\documentclass[a4paper,12pt,twoside]{amsart}
\usepackage{amsmath,amsthm,amsfonts}

\newtheorem{thm}{Theorem}[section]

\newtheorem{cor}[thm]{Corollary}
\newtheorem{lem}[thm]{Lemma}

\newtheorem{prop}[thm]{Proposition}
\newtheorem{proposition}[thm]{Proposition}
\newtheorem{definition}[thm]{Definition}
\newtheorem{defn}[thm]{Definition}

\newtheorem{assumption}[thm]{Assumption}
\newtheorem{nota}[thm]{Notation}
\newtheorem{notas}[thm]{Notations}

\theoremstyle{definition}

\newtheorem{exs}[thm]{Examples}
\newtheorem{rem}[thm]{Remark}
\newtheorem{rems}[thm]{Remarks}

\newcommand{\R}{{\mathbb{R}}}
\newcommand{\A}{{\mathbb{A}}}
\newcommand{\E}{{\mathbb{E}}}
\newcommand{\Q}{{\mathbb{Q}}}
\newcommand{\C}{{\mathbb{C}}}

\newcommand{\Z}{{\mathbb{Z}}}
\newcommand{\N}{{\mathbb{N}}}
\newcommand{\T}{{\mathbb{T}}}
\newcommand{\PP}{{\mathbb{P}}}

\newcommand{\Cal}{\mathcal}
\newcommand{\cal}{\mathcal}
\newcommand{\cH}{\mathcal{H}}

\def \mod {{\rm \ mod \,} }

\def \diag {{\rm diag} \,}
\def \Var {{\rm Var}}
\def\e{{\rm e}}

\def\j{{\underline j}} \def\k{{\underline k}}
 \def\el{{\underline \ell}}
\def\t{{\underline t}} \def\n {{\underline n}}
\def\0{{\underline 0}} \def\m {{\underline m}}
\def\a{{\underline a}} 
\def\f{{\underline f}} \def\p {{\underline p}}
\def\u{{\underline u}} \def\r{{\underline r}}
\def\v{{\underline v}} \def\w{{\underline w}}
\def\s{{\underline s}} \def\x{{\underline x}}
\def\z{{\underline z}} \def\q{{\underline q}}

\def\for{\text{ for }}
\def\and{\text{ and }}

\def\stm0{{\setminus \{\0\}}}

\def\eop{\qed}
\def\Proof {\vskip -3mm {{\it Proof}. }}
\def\proof {\vskip -3mm {{\it Proof}. }}
\def\Sect {{Sect. \hskip -1mm}}

\def \leb {{\rm Leb}}
\def \Log {{\rm Log}}
\def \Card {{\rm Card}}

\begin{document}
\baselineskip 15pt
\parindent=0mm

\title[CLT for random walks of
commuting endomorphisms] {CLT for random walks of commuting \\
endomorphisms on compact abelian groups}

\date{26 October 2014}
\author{Guy Cohen and Jean-Pierre Conze}
\address{Guy Cohen, \hfill \break Dept. of Electrical Engineering,
\hfill \break Ben-Gurion University, Israel}
\email{guycohen@ee.bgu.ac.il}
\address{Jean-Pierre Conze,
\hfill \break IRMAR, CNRS UMR 6625, \hfill \break University of
Rennes I, Campus de Beaulieu, 35042 Rennes Cedex, France}
\email{conze@univ-rennes1.fr}


\begin{abstract} Let $\Cal S$ be an abelian group of automorphisms of
a probability space $(X, {\Cal A}, \mu)$ with a finite system of
generators $(A_1, ..., A_d)$. Let $A^{\el}$ denote $A_1^{\ell_1} ...
A_d^{\ell_d}$, for ${\el}= (\ell_1, ..., \ell_d)$. If $(Z_k)$ is a
random walk on $\Z^d$, one can study the asymptotic
distribution of the sums $\sum_{k=0}^{n-1} \, f \circ A^{\,{Z_k(\omega)}}$
and $\sum_{\el \in \Z^d} \PP(Z_n= \el) \, A^\el f$, for a function $f$ on $X$.

In particular, given a random walk on commuting matrices in $SL(\rho, \Z)$ or in
${\Cal M}^*(\rho, \Z)$ acting on the torus $\T^\rho$, $\rho \geq 1$,
what is the asymptotic distribution of the associated ergodic
sums along the random walk for a smooth function on $\T^\rho$ after normalization?

In this paper,  we prove a central
limit theorem when $X$ is a compact abelian connected group $G$
endowed with its Haar measure (e.g. a torus or a connected extension of a torus),
$\Cal S$ a totally ergodic $d$-dimensional group of commuting algebraic
automorphisms of $G$ and $f$ a regular function on $G$. The proof is based
on the cumulant method and on preliminary results on the spectral properties
of the action of $\Cal S$, on random walks and on the variance of the associated ergodic sums.
\end{abstract}

\maketitle \tableofcontents

\section*{\bf Introduction}

Let $\Cal S$ be a finitely generated abelian group and let $(T^s, \,
s \in \Cal S)$ be a measure preserving action of $\Cal S$ on a
probability space $(X, {\Cal A}, \mu)$. A probability distribution
$(p_s, \, s \in \Cal S)$ on $\Cal S$ defines a random walk (r.w.)
$(Z_n)$ on $\Cal S$. We obtain a Markov chain on $X$ (and a random
walk on the orbits of the action of $\Cal S$) by defining the
transition probability from $x \in X$ to $T^s x$ as $p_s$. The
Markov operator $P$ of the corresponding Markov chain is $Pf(x) =
\sum_{s \in \Cal S} \, p_s f(T^s x)$. Limit theorems have been
investigated for the associated random process. For instance,
conditions on $f \in L_0^2(\mu)$ are given in \cite{DerLin07} for a
``quenched" central limit theorem for the ergodic sums along the r.w.
$\sum_{k=0}^{n-1} f(T^{Z_k(\omega)} x$) when $\Cal S$ is $\Z$. ``Quenched",
there, is understood for a.e. $x \in X$ and in law with respect to $\omega$.
For other limit theorems in the quenched setting, see for instance \cite{CuMe14}
and references therein.

Another point of view, in the study of random dynamical systems, is to prove limiting laws
for a.e. fixed $\omega$ (cf. \cite{KiLi06}). This is our framework:
for a fixed $\omega$, we consider the asymptotic distribution of the above ergodic
sums after normalization with respect to $\mu$.

More precisely, let $(T_1, ..., T_d)$ be a system of generators in ${\Cal S}$. Every
$T \in {\Cal S}$ can be represented as $T = T^{\el}= T_1^{\ell_1} ...
T_d^{\ell_d}$, for ${\el}= (\ell_1, ..., \ell_d) \in {\N}^d$. (The elements of $\Z^d$ are
underlined to distinguish them from the scalars.) For a
function $f$ on $X$ and $T \in \Cal S$, $Tf$ denotes the composition
$f \circ T$. The map $f \to Tf$ defines an isometry on ${\Cal H} =
L^2_0(\mu)$, the space of square integrable functions $f$ on $X$
such that $\mu(f) =0$.

Given a r.w. $W = (Z_n)$ on $\Z^d$, we study the asymptotic distribution
(for a fixed $\omega$ and with respect to the measure
$\mu$ on $X$) of $\sum_{k=0}^{n-1} \, T^{Z_k(\omega)}f, {n \to \infty}$,
after normalization. We consider also iterates of the Markov
operator $P$ introduced above (called below barycenter operator) and
the asymptotic distribution of $(P^n f)_{n \geq 1}$ after normalization.

Other sums for the random field $(Tf, T \in {\Cal S})$,  $f \in L^2_0(\mu)$,
can be considered. For instance, if $(D_n) \subset \N^d$ is an increasing sequence of
domains, a question is the asymptotic normality of $|D_n|^{-\frac12} \sum_{{\el} \, \in D_n} \,
T^{\,{\el}} f$ and of the multidimensional ``periodogram"
$|D_n|^{-\frac12} \sum_{\el \, \in D_n} \, e^{2\pi i \langle \el,
\theta \rangle} \, T^{{\el}} f, \ \theta \in \R^d$. All above
questions can be formulated in the framework of what is called below
``summation sequences" and their associated kernel.

The specific model which is studied here is the following one: $(X,
{\Cal A}, \mu)$ will be a compact abelian connected group $G$
endowed with its Borel $\sigma$-algebra ${\Cal A}$ and its Haar
measure $\mu$ and $\Cal S$ will be a semigroup of commuting
algebraic endomorphisms of $G$. This extends the classical situation
of a CLT for a single endomorphism (R.~Fortet (1940) and M.~Kac
(1946) for $G = \T^1$, V. Leonov (1960) for a general ergodic
endomorphism of $G$).

The main result (Theorem \ref{main}) is a quenched CLT for
($\sum_{k=0}^{n-1} \, T^{Z_k(\omega)}f)$, when $f$ is a regular
function on $G$, $(T^\el)_{\el \in \Z^d}$ a $\Z^d$-action on $G$ by
automorphisms and $(Z_k)$ a r.w. on $\Z^d$. After reduction of the
r.w., we examine three different cases: a) moment of order 2,
centering, dimension 1; b) moment of order 2, centering, dimension
2; c) transient case. This covers all cases when there is a moment
of order 2. As usual, the transient case is the easiest. The
recurrent case requires to study the self-intersections of the r.w.

The paper is organized as follows. In \Sect \ref{summSeqSect} we
discuss methods of summation and recall some facts on the spectral
analysis of a finitely generated group of unitary operators with
Lebesgue spectrum. The notion of ``regular" summation sequence
allows a unified treatment of different cases, including ergodic
sums over sets or sums along random sequences generated by random
walks. In \Sect\ref{rwsummSec} we prove auxiliary results on the
regularity of quenched and barycenter summation sequences defined by
a r.w. on $\Z^d$. These first two sections apply to a general
$\Z^d$-action with Lebesgue spectrum.

In \Sect \ref{endoGroup}, the model of multidimensional actions by
endomorphisms or automorphisms on a compact abelian group $G$ is
presented. We recall how to construct totally ergodic toral
$\Z^d$-actions by automorphisms and we give explicit examples. The link
between the regularity of a function $f$ on $G$ and its spectral
density $\varphi_f$ is established, with a specific treatment for
tori for which the required regularity for $f$ is weaker. A
sufficient condition for the CLT on a compact abelian connected group $G$ is
given in \Sect \ref{CLTGroup}.

The results of the first sections are applied in \Sect
\ref{applRW}. A quenched CLT for the ergodic sums
of regular functions $f$ along a r.w. on $G$ by commuting automorphisms
is shown when $G$ is connected. For a transient r.w. the variance, which is
related to a measure on $\T^d$ with an absolutely continuous part,
does not vanish. Another summation method, iteration of
barycenters, yields a polynomial decay of the iterates for regular
functions and a CLT for the normalized sums, the nullity of the variance
being characterized in terms of coboundary. An appendix (\Sect \ref{selfInt})
is devoted to the self-intersections properties of a r.w. In \Sect \ref{cumulAppend},
we recall the cumulant method used in the proof of the CLT.

To conclude this introduction, let us mention that in a previous
work \cite{CohCo13} we extended to ergodic sums of
multidimensional actions by endomorphisms the CLT proved by T.~Fukuyama
and B.~Petit \cite{FuPe01} for coprime integers acting on
the circle. After completing it, we were informed of the results of
M. Levin \cite{Lev13} showing the CLT for ergodic sums over
rectangles for actions by endomorphisms on tori. The proof of the
CLT for sums over rectangles is based in both approaches, as well as
in \cite{FuPe01}, on results on $S$-units which imply mixing of all
orders for connected groups (K. Schmidt and T. Ward
\cite{SchWar93}). The cumulant method used here is also based on
this property of mixing of all orders.
\section{\bf Summation sequences } \label{summSeqSect}

Let us consider first the general framework of an {\it abelian
finitely generated group} ${\Cal S}$ (isomorphic to $\Z^d$) of unitary
operators on a Hilbert space ${\Cal H}$. There exists a system of
independent generators $(T_1, ..., T_d)$ in ${\Cal S}$ and each
element of $\Cal S$ can be written in a unique way as $T^{\el}=
T_1^{\ell_1} ... T_d^{\ell_d}, \text{ with }{\el}= (\ell_1, ...,
\ell_d) \in \Z^d$.

Our main example (Sect. \ref{endoGroup}) will be given by groups of
commuting automorphisms (or extensions of semigroups of commuting
surjective endomorphisms) of a compact abelian group $G$. They act by
composition on ${\Cal H} = L_0^2(G, \mu)$ (with $\mu$ the Haar
measure of $G$) and yield examples of $\Z^d$-actions by unitary
operators.

Given $f \in \Cal H$, there are various choices of summation
sequences for the random field $(T^{\el} f, \el \in \Z^d)$. In the first
part of this section, we discuss this point and the behavior of the
kernels associated to summation sequences. The second part
of the section is devoted to the spectral analysis of
$d$-dimensional commuting actions with Lebesgue spectrum.

\subsection{\bf Summation sequences, kernels, examples} \label{summatEx}
\begin{definition} \label{sumSeq} {\rm We call {\it summation sequence}
a sequence $(R_n)_{n \geq 1}$ of functions from $\Z^d$ to $\R^+$
with $0 < \sum_{\el \in \Z^d} R_n(\el) < +\infty, \ \forall n \geq
1$. Given ${\Cal S} = \{T^\el, \el \in \Z^d\}$ and $f \in {\Cal H}$, the
associated sums are $\sum_{\el \in \Z^d} R_n(\el) \, T^\el f$.

Let $\check R_n(\el) := R_n(-\el)$. The normalized nonnegative
kernel $\tilde R_n$ (with $\|\tilde R_n\|_{L^1(\T^d)} = 1$) is
\begin{eqnarray}
\tilde R_n(\t) = {|\sum_{\el \in \Z^d} R_n(\el) e^{2\pi i \langle
\el, \t \rangle}|^2 \over \sum_{\el \in \Z^d} |R_n(\el)|^2} =
\sum_{\el \in \Z^d} {(R_n*\check R_n)(\el) \over (R_n*\check
R_n)(\0) } \, e^{2\pi i \langle \el, \t \rangle}, \ \t \in \T^d.
\label{chch}
\end{eqnarray}}
\end{definition}
\begin{definition} \label{regularProcess} {\rm We say that $(R_n)$
is $\zeta$-{\it regular}, if $(\tilde R_n)_{n \geq 1}$ weakly
converges to a probability measure $\zeta$ on $\T^d$, i.e.,
$\int_{\T^d} \tilde R_n \,\varphi \, d\t \underset{n \to \infty}
{\longrightarrow}\int_{\T^d} \varphi \, d\zeta$ for every continuous
$\varphi$ on $\T^d$.

The existence of the limit: $\displaystyle L(\p) = \lim_{n \to
\infty} \int \tilde R_n(\t) \, e^{-2\pi i \langle \p, \t \rangle} \
d\t = \lim_n {(R_n*\check R_n)(\p) \over (R_n*\check R_n)(\0)}$, for
all $\p \in \Z^d$, is equivalent to $\zeta$-regularity with $\hat
\zeta (\p) = L(\p)$. Note that $\tilde R_n$ and $\zeta$ are even.}
\end{definition}
\begin{exs} \label{exples1} {\it Summation over sets} \end{exs}
\vskip -3mm A class of summation sequences is given by summation
over sets. If $(D_{n})_{n \geq 1}$ is a sequence of finite subsets
of $\Z^d$ and $R_n = 1_{D_n}$, we get the ``ergodic" sums
$\sum_{{\el} \in D_n} \, T^\el f$. The simplest choice for $(D_n)$
is an increasing family of $d$-dimensional rectangles. The sequence
$(R_n) = (1_{D_n})$ is $\delta_0$-regular if and only if $(D_n)$
satisfies
\begin{eqnarray}
\lim_{n\to\infty} |D_n|^{-1} |(D_n+{\p}) \, \cap \, D_n| =1, \
\forall {\p}\in \Z^d \ (\text{F\o{}lner condition)}.
\label{FolnCond}
\end{eqnarray}
a) (Rectangles) For ${\underline N} =(N_1, ..., N_d) \in \N^d$ and
$D_{\underline N} := \{{\el}\in \N^d:\ \ell_i \le N_i, \, i=1, ...,
d\}$, the kernels are the $d$-dimensional Fej\'er kernels
$K_{\underline{N}}(t_1,..., t_d) = K_{N_1}(t_1)\cdots K_{N_d}(t_d)$
on $\T^d$, where $K_N(t)$ is the one-dimensional Fej\'er kernel
${1\over N} ({\sin \pi N t \over \sin \pi t })^2$.

b) A family of examples satisfying (\ref{FolnCond}) can be obtained
by taking a non-empty domain $D \subset \R^d$ with ``smooth"
boundary and finite area and putting $D_n=\lambda_n D\cap \Z^d$,
where $(\lambda_n)$ is an increasing sequence of real numbers
tending to $+\infty$.

c) If $(D_\n)$ is a (non F\o{}lner) sequence of domains in $\Z^d$
such that $\lim_\n \frac{|(D_n+{\p}) \cap D_n|} {|D_n|}=0, \forall
{\p}\not={{0}}$, the associated kernel $(\tilde R_n)$ is
$\zeta$-regular, with $\zeta$ the uniform measure on $\T^d$.

\begin{exs} \label{exples2} {\it Sequential summation} \end{exs}
\vskip -3mm Let $(\x_k)_{k \geq 0}$ be a sequence in $\Z^d$. Putting
$\z_n:= \sum_{k=0}^{n-1} \x_k$, $n \geq 1$, and $R_n(\el) =
\sum_{k=0}^{n-1} 1_{\z_k = \el}$ for $\el \in \Z^d$, we get a
summation sequence of sequential type. Here, $\sum_{\el \in \Z^d}
R_n(\el) = n$.

The associated sums are $\sum_{\el \in \Z^d} R_n(\el) \, T^\el f =
\sum_{k=0}^{n-1} T^{\z_k} f$. Let us define
\begin{eqnarray}
&&v_{n} = v_{n, \0} = \#\{0 \leq k', k < n: \ \z_k= \z_{k'}\}
= \sum_{\el \in \Z^d} R_n^2(\el),   \label{defVnGen0} \\
&&v_{n, \p} := \#\{0 \leq k', k < n: \ \z_k - \z_{k'} = \p \} =
\sum_{0 \leq k', k < n} 1_{z_k - z_{k'} = \p}.  \label{defVnGenEl}
\end{eqnarray}
The quantity $v_n$ is the number of ``self-intersections" of
$(\z_k)_{k \geq 1}$ up to time $n$. We have
\begin{eqnarray}
\tilde R_n(\t) = \sum_{\p \in \Z^d} {v_{n,\p} \over v_{n, \0}} \, e^{2\pi i
\langle \p, \t \rangle}, \ \t \in \T^d, \label{chch1}
\end{eqnarray}
and the $\zeta$-regularity of $(\z_k)_{k \geq 1}$ is equivalent to
$\displaystyle \lim_n {v_{n,\p} \over v_{n, \0}} = \hat \zeta(\p),
\forall \p \in \Z^d$. Remark that $v_{n} - v_{n, p} \geq 0$ by
Cauchy-Schwarz inequality, since $v_{n,p}=\sum_{\ell} R_n(\ell)
R_n(\ell+p)$.

In \Sect 2, we consider random sequential summation sequences
defined by random walks.

\vskip 3mm \subsubsection{\bf Reduction to the ``genuine" dimension}
\label{reduc}

\

A reduction to the smallest lattice containing the support of the
summation sequence is convenient and can be done using Minkowski
lemma.

Below $(e_1, ..., e_d)$ is the canonical basis of $\R^d$. For $r <
d$, $E_r$ (resp. ${\Cal E}_{r}$) is the lattice over $\Z$ (resp.
the vector space) generated by $(e_1, ..., e_{r})$ and $F_{d -r}$
the sublattice of $\Z^d$ generated by $(e_{r+1}, ..., e_d)$.

The set of non singular $d \times d$-matrices with coefficients in
$\Z$ is denoted by ${\Cal M}^*(d, \Z)$.
\begin{lem} \label{ImLatt} a) Let $L$ be a sublattice of $\,\Z^d$. If
$r \geq 1$ is the dimension of the $\R$-vector space $\Cal L$
spanned by $L$, there exists $C \in SL(d, \Z)$ such that $\Cal L = C
{\Cal E}_{r}$. \hfill \break b) If $r = d$, there exists $B \in
{\Cal M}^*(d, \Z)$, such that $L = B \Z^d$ and $\Card (\Z^{d} / L) =
|\det(B)|$.
\end{lem}
\Proof a) The proof is by induction on $r$. By definition of $\Cal
L$ we can find a non zero vector $\v = (v_1, ..., v_d)$  in $\Cal L$
with integral coordinates such that $gcd(v_1, ..., v_d) =1$.

Since a vector in $\Z^d$ is extendable to a basis of $\Z^d$ if and
only if the gcd of its coordinates is 1, we can construct an
integral matrix $C_1$ such that $\det C_1 = 1$ with $\v$ as first
column.

For $\w\in L$, let $\tilde \w := \langle \w, \w \rangle \v - \langle
\v, \w\rangle \v$. The space $\Cal L \cap \v^\perp$ is generated by
the lattice $\{\tilde \w, \w \in L\}$ and has dimension $r-1$. By
the induction hypothesis, there is an integral matrix $C_2$ with
$\det C_2 = 1$, such that $C_2 e_r = e_r$ and $C_2 {\Cal E}_{r-1} =
C_1^{-1} (\Cal L \cap \v^\perp)$. Taking $C = C_1C_2$, we have $C
{\Cal E}_{r-1} = C_1C_2{\Cal E}_{r-1} = \Cal L \cap \v^\perp$, $C
e_r = C_1 C_2 e_r = C_1 e_r = \v \in \Cal L$.

b) The existence of $B$ is equivalent to the existence in $L$ of a
set $\{\f_1,\ldots, \f_{d}\}$ of linearly independent vectors
generating $L$. Let us construct such a set.

The set $\{x_1: \x = (x_1, ..., x_d) \in L\}$ is an ideal, hence the
set of multiples of some integer $a$. We can assume that $a \not =
0$ (otherwise we permute the coordinates). Let $\f =(f_1, ..., f_d)
\in L$ be such that $f_1 = a$. The set $\{\x - a^{-1} x_1 \f, \x \in
L\}$ is a sublattice $L_0$ of $L$ with $0$ as first coordinate of
every element. It is of dimension $d-1$ over $\R$, hence, by
induction hypothesis, generated by linearly independent vectors
$\{\f_1,\ldots, \f_{d-1}\}$ in $L_0$. Clearly $\{\f_1, \ldots,
\f_{d-1}, \f\}$ is a set of linearly independent vectors generating
$L$.

Let $\Cal R$ be a set of representatives of the classes of $\Z^d$
modulo $B \Z^d$ and let $K$ be the unit cube in $\R^d$. Since if $L$
is a cofinite lattice in $\Z^d$, two fundamental domains $F_1, F_2$
in $\R^d$ for $L$ have the same measure: $\leb(F_1) = \leb(F_2)$.
Taking $F_1= B K$ and $F_2 = \bigcup_{\r \in \Cal R} (K + \r)$, this
implies $|\det(B)| = \leb (B K) = \leb \Bigl(\bigcup_{\r \in \Cal R}
(K + \r) \Bigr )= \Card (\Cal R)$.\eop

Let us consider the lattice $L$ in $\Z^d$ generated by the elements
$\el \in \Z^d$ such that $R_n(\el) >0$ for some $n$. For instance,
for a sequential summation generated by $(\x_k)_{k \geq 0}$, the
lattice $L$ is the group generated in $\Z^d$ by the vectors $\x_k$.

By the previous lemma, after a change of basis given by a matrix in
$SL(d, \Z)$ and reduction, we can assume without loss of generality
that $d$ is the {\it genuine dimension}, i.e., that the vector space
$\Cal L$ generated by $L$ over $\R$ has dimension $d$.

\subsection{\bf Spectral analysis, Lebesgue spectrum} \label{LebSpec}

\

Recall that,  if ${\Cal S}$ is an abelian group of unitary
operators, for every $f \in \cH$ there is a positive finite measure
$\nu_f$ on $\T^d$, the spectral measure of $f$, with Fourier
coefficients $\hat \nu_f({\el}) = \langle T^{\el} f, f \rangle$,
$\el \in \Z^d$. When $\nu_f$ is absolutely continuous, its density
is denoted by $\varphi_f$.

{\it In what follows, we assume that $\Cal S$ (isomorphic to $\Z^d$)
has the Lebesgue spectrum property for its action on ${\Cal H}$}, i.e.,
there exists a closed subspace ${\Cal K}_0$ such that $\{T^\el {\Cal
K}_0, \, \el \in \Z^d\}$ is a family of pairwise orthogonal
subspaces spanning a dense subspace in ${\Cal H}$. If $(\psi_j)_{j
\in \Cal J}$ is an orthonormal basis of ${\Cal K}_0$, $\{T^\el
\psi_j, j \in \Cal J, \, \el \in \Z^d \}$ is an orthonormal basis of
${\Cal H}$.

In \Sect \ref{endoGroup}, for algebraic automorphisms of compact
abelian groups, the characters will provide a natural basis.

Let $\cH_j$ be the closed subspace generated by $(T^\el \psi_j)_{\el
\in \Z^d}$ and $f_j$ the orthogonal projection of $f$ on $\cH_j$. We
have $\nu_f = \sum_j \nu_{f_j}$. For $j \, \in \Cal J$, we denote by
$\gamma_j$ an everywhere finite square integrable function on $\T^d$
with Fourier coefficients $a_{j, \, \el} := \langle f, T^{\el}
\psi_j \rangle$. A version of the density of the spectral measure
corresponding to $f_j$ is $|\gamma_j|^2$.

{\it By orthogonality of the subspaces $\cH_j$, it follows that, for
every $f \in \cH$, $\nu_f$ has a density $\varphi_f$ in $L^1(d\t)$
which reads $\varphi_f(\t) = \sum_{j \in \Cal J} |\gamma_j(\t)|^2$.}

Changing the system of independent generators induces composition of
the spectral density by an automorphism acting on $\T^d$. After the
reduction performed in \ref{reduc} (Lemma \ref{ImLatt}) we obtain an
action with Lebesgue spectrum with, possibly, a smaller dimension.

We have: $\int_{\T^d} \sum_{j \in \Cal J} |\gamma_j(\t)|^2 \ d\t =
\sum_{j \, \in \Cal J} \sum_{\el \in \Z^d} |a_{j, \el}|^2 =
\int_{\T^d} \varphi_f(\t) \ d\t =\|f\|_2^2 < \infty$. For $\t$ in
$\T^d$, let $M_\t f$ in ${\Cal K}_0$ be defined, when the series
converges in $\Cal H$, by: $M_\t f := \sum_j \gamma_j(\t) \,
\psi_j$.

Under the condition $ \sum_{j\in \Cal
J}\big(\sum_{\el}|a_{j,\el}|\big)^2 < +\infty$, $M_\t f$ is defined
for every $\t$, the function $\t \to \|M_\t\|_2^2$ is a continuous
version of $\varphi_f$. For a general $f \in \Cal H$, it is defined
for $\t$ in a set of full measure in $\T^d$ (cf. \cite{CohCo12}).

{\bf Variance for summation sequences}

Suppose that $(R_n)$ is a $\zeta$-regular summation sequence. Let
$f$ be in $L_0^2(\mu)$ with a continuous spectral density
$\varphi_f$. By the spectral theorem, we have for $\theta \in \T^d$:
\begin{eqnarray}
(\sum_\el \, R_n^2(\el))^{-1} \|\sum_{\el \in \Z^d} \, R_n(\el) \,
e^{2\pi i \langle \el, \theta \rangle} \, T^\el f\|_2^2 = (\tilde
R_n * \varphi_f)(\theta) \underset{n \to \infty} {\longrightarrow}
(\zeta * \varphi_f)(\theta). \label{spectrThm}
\end{eqnarray}
{\it Examples.} 1) If $(D_n)$ is a F\o{}lner sequence of sets in
$\Z^d$, then (for $\theta=0$) $\zeta = \delta_0$ and we obtain the
usual asymptotic variance $\sigma^2(f) = \varphi_f(\0)$. More
generally the usual asymptotic variance at $\theta$ for the
``rotated" ergodic sums $R_n^\theta f$ is
\begin{eqnarray*}
\lim_n |D_n|^{-1} \, \|\sum_{{\el} \, \in D_n} \, e^{2\pi i \langle
\el, \theta \rangle} T^{{\el}} f\|_2^2 = \lim_n(\tilde R_n *
\varphi_f)(\theta) = (\delta_0 * \varphi_f)(\theta) =
\varphi_f(\theta).
\end{eqnarray*}
When $(D_n)$ is a sequence of $d$-dimensional cubes, by
the Fej\'er-Lebesgue theorem, for every $f$ in $\Cal H$, since $\varphi_f \in L^1(\T^d)$,
for a.e. $\theta$, the variance at $\theta$ exists and is equal to
$\varphi_f(\theta)$.

2) Let $(\x_k)_{k \geq 0}$ be a sequence in $\Z^d$ and $\z_n =
\sum_{k=0}^{n-1} \x_k$. If $(\z_n)$ is $\zeta$-regular, then (cf.
(\ref{chch1})): $\displaystyle \lim_n v_n^{-1} \, \|\sum_{k=0}^{n-1}
\, T^{{\z_k}} f\|_2^2 = \lim_n \int \tilde R_n \varphi_f \, d\t =
\zeta(\varphi_f)$.

\section{\bf Random walks} \label{rwsummSec}
\subsection{\bf Random summation sequences (rws) defined by random walks} \label{rwsummSubSec}

\

For $d \geq 1$ and $J$ a set of indices, let $\Sigma := \{\el_j, j
\in J\}$ be a set of vectors in $\Z^d$ and $(p_j, j \in J)$ a
probability vector with $p_j > 0, \forall j \in J$.

Let $\nu$ denote the probability distribution with support $\Sigma$
on $\Z^d$ defined by $\nu(\el) = p_j$ for $\el = \el_j$, $j \in J$.
The euclidian norm of $\el \in \Z^d$ is denoted by $|\el|$.

Let $(X_k)_{k \in \Z}$ be a sequence of i.i.d. $\Z^d$-valued random
variables with common distribution $\nu$, i.e., $\PP(X_k = \el_j) =
p_j$, $j \in J$. The associated {\it random walk} $W = (Z_n)$ in
$\Z^d$ (starting from $\0$) is defined by $Z_0 := \0$, $Z_n := X_0
+... + X_{n-1}$, $n \geq 1$.

The characteristic function $\Psi$ of $X_0$ (computed with a
coefficient $2\pi$) is
\begin{eqnarray}
&&\Psi(\t) = \E[\e^{2\pi i\langle X_0, \t \rangle}] = \sum_{\el \in
\Z^d} \PP(X_0 = \el) \, e^{2 \pi i \langle \el, \t\rangle} = \sum_{j
\in J} p_j e^{2 \pi i \langle \el_j, \t\rangle}, \, \t \in \T^d.
\label{charact} \\
&&\text{It satisfies: } \PP(Z_{n} = \k) = \int_{\T^d} \Psi^{n}(\t)
e^{-2\pi i\langle \k, \t \rangle} \ d\t. \label{RwFour}
\end{eqnarray}
If $(\Omega, \PP)$ denotes the product space $((\Z^d)^\Z,
\nu^{\otimes \Z})$ with coordinates $(X_n)_{n \in \Z}$ and $\tau$
the shift, then $W$ can be viewed as the cocycle generated by the
function $\omega \to X_0(\omega)$ under the action of $\tau$, i.e.,
$Z_n = \sum_{k=0}^{n-1} X_0 \circ \tau^k, n \geq 1$.

Let $T_1, ..., T_d$ be $d$ commuting unitary operators on a Hilbert
space $\Cal H$ generating a representation of $\Z^d$ in $\Cal H$
with Lebesgue spectrum. Let $(Z_n)_{n \geq 1} = ((Z_n^1, ...,
Z_n^d))_{n \geq 1}$ be a r.w. with values in $\Z^d$. We consider the
random sequence of unitary operators $(T^{Z_n})_{n \geq 1}$, with
$T^{Z_n} = T_1^{Z_n^1}...T_d^{Z_n^d}$. For $f \in \Cal H$, the
``quenched" process (for $\omega$ fixed) of the ergodic sums along
the random walk is
\begin{eqnarray}
&\sum_{k=0}^{n-1} \, T^{Z_k(\omega)}f = \sum_{\el \in \Z^d} \,
R_n(\omega, \el) \, T^{\el}f, \text{ with } R_n(\omega, \el) =
\sum_{k=0}^{n-1} 1_{Z_k(\omega) = \el}, \, n \geq 1.
\label{randSumseq}
\end{eqnarray}
The sequence of ``local times" $(R_n(\omega, \el))$ will be called
{\it a random walk sequential summation} (abbreviated in ``r.w.
summation" or in ``rws"). It satisfies $\sum_\el R_n(\omega,\el) =
n$.

We consider also the (Markovian) barycenter operator $P: f \to
\sum_{j \in J} p_j \ T^{\el_j} f = \E_\PP \, [T^{X_0(.)}f]$ and its
powers:
\begin{eqnarray}
&&P^n f = \E_\PP [T^{Z_n(.)} f]) = \sum_{\el \in \Z^d} R_n(\el) \
T^\el f, \text{ with } R_n(\el) = \PP(Z_n = \el). \label{barycproc}
\end{eqnarray}
When the operators $T_j$ are given by measure preserving
transformations of a probability space $(X, \Cal A, \mu)$, we call
{\it quenched limit theorem} a limit theorem for (\ref{randSumseq})
w.r.t. $\mu$. If the limit is taken w.r.t. $\PP \times \mu$, it is
called {\it  annealed limit theorem}.

\begin{notas} {\rm The random versions of (\ref{defVnGen0})
and (\ref{defVnGenEl}) are
\begin{eqnarray}
&&V_n(\omega) := \#\{0 \leq k', k < n: \ Z_k(\omega) =
Z_{k'}(\omega)
\} = \sum_{\el \in \Z^d} R_n^2(\omega, \el),\label{defVn} \\
&&V_{n, \p}(\omega) := \#\{0 \leq k', k < n: \ Z_k(\omega) -
Z_{k'}(\omega) = \p \} = \sum_{0 \leq k', k < n} 1_{Z_k(\omega) -
Z_{k'}(\omega) = \p}.
\end{eqnarray}
}\end{notas} $V_n(\omega) = V_{n, \0}(\omega)$ is the number of
``self-intersections" of $W$. $V_{n, \p}(\omega)$ is a sum of
ergodic sums.

We have:
\begin{eqnarray}
&&V_{n, \p}(\omega) = n1_{\p= \0} + \sum_{k = 1}^{n-1} \, \sum_{j =
0}^{n-k-1} (1_{Z_k (\tau^j \omega) = \p} + 1_{Z_k (\tau^j \omega) =
-\p});\label{VnpErgSum} \\
&&\E V_{n, \p} = n1_{\p= 0} + \sum_{k = 1}^{n-1} \,(n-k) \, [\PP(Z_k
= \p) + \PP(Z_k = -\p)]. \label{EspVnpErgSum}
\end{eqnarray}

\vskip 2mm \goodbreak {\bf Variance for quenched processes}

 For a r.w.\ summation $R_n(\omega, \el) = \sum_{k=0}^{n-1}
1_{Z_k(\omega) =\el}$, let $\underline R_n^\omega(t)$ and
$\widetilde R_n^\omega(t)$ denote, respectively, the non normalized
and normalized kernels (cf. (\ref{chch}))
\begin{eqnarray}
&&\underline R_n^\omega(\t) = |\sum_{k=0}^{n-1} e^{2\pi i \langle
Z_k(\omega), \t \rangle} |^2 = |\sum_{\el \in \Z^d} R_n(\omega, \el)
e^{2\pi i \langle \el, t \rangle}|^2, \ \widetilde R_n^\omega(\t) =
\underline R_n^\omega(\t) / \int_{\T^d} \underline R_n^\omega d\t.
\label{kernDefOmNorm}
\end{eqnarray}
Recall that, if $\varphi_f$ is the spectral density of $f \in \Cal
H$ for the action of $\Z^d$, we have $\displaystyle \|\sum_{k=1}^{n}
T^{Z_k(\omega)} f\|^2 = \int_{\T^d} \underline R_n^\omega(t) \,
\varphi_f(\t) \, d\t$ by the spectral theorem (cf.
(\ref{spectrThm})) and
\begin{eqnarray}
&&\int_{\T^d} \frac1n \underline R_n^\omega(t) \,e^{-2\pi i \langle
\p, \t \rangle} d\t = \frac1n V_{n, \p}(\omega) = 1_{\p= \0} +
\sum_{k = 1}^{n-1} \, \frac1n\sum_{j = 0}^{n-k-1} [1_{Z_k (\tau^j
\omega) = \p} + 1_{Z_k (\tau^j \omega) = - \p}]. \label{FourKern}
\end{eqnarray}
{\it A rws defined by a r.w. $(Z_n)$ is said to be $\zeta$-regular
if $(\tilde R_n^\omega)_{n \geq 1}$ weakly converges to a
probability measure $\zeta$ (not depending on $\omega$) on $\T^d$ for a.e. $\omega$.}

This is equivalent to $\lim_{n \to \infty} {V_{n, \p}(\omega) \over
V_{n}(\omega)} = \hat \zeta(\p)$, for a.e. $\omega$ and all $\p \in
\Z^d$ (see Definition \ref{regularProcess} and Example
\ref{exples2}). After auxiliary results on r.w.'s and
self-intersections, we will show that every rws is $\zeta$-regular
for some measure $\zeta$.

\subsection{\bf Auxiliary results on random walks}
\label{prelimRW}

\

We recall now some general results on random walks in $\Z^d$ (cf.
\cite{Sp64}). Most of them are classical. Nevertheless, since we do
not assume strict aperiodicity and in order to explicit the constants
in the limit theorems, we include reminders and some proofs.
We start with definitions and preliminary results.

\subsubsection{\bf Reduced form of a random walk \label{redRW}}

\

Let $L(W)$ (or simply $L$) be the sublattice of $\Z^d$ generated by
the support $\Sigma =\{\el_j, j \in J\}$ of $\nu$. From now on
(without loss of generality, as shown by Lemma \ref{ImLatt}), {\it
we assume that $L(W)$ is cofinite in $\Z^d$ and we say that $W$ is
reduced} (in other words ``genuinely $d$-dimensional", cf.
\cite{Sp64}). Therefore the vector space $\Cal L$ generated by $L$
is $\R^d$ and there is $B \in {\Cal M}^*(d, \Z)$ such that $L = B
\Z^d$. The rank $d$ will be sometimes denoted by $d(W)$.

Let $D = D(W)$ be the sublattice of $\Z^d$ generated by $\{\el_j -
\el_{j'}, j, j' \in J\}$, $\Cal D$ the vector subspace of $\R^d$
generated by $D$ and $\Cal D^\perp$ its orthogonal supplementary in
$\R^d$. We denote $\dim(\Cal D)$ by $d_0(W)$. With the notations of
\ref{reduc}, by Lemma \ref{ImLatt}, there is $B_1 \in {\Cal M}^*(d,
\Z)$ such that $D = B_1 \Z^d$ if $d_0(W) = d$ and $D = B_1 F_{d -1}$
if $d_0(W) = d-1$.
\begin{lem} \label{DLatt} a) The quotient group $L(W)/D(W)$ is
cyclic (each $\ell_j$ is a generator of the group). It is finite if
and only if $d_0(W) = d(W)$. \hfill \break b) If $W$ is reduced,
then $d_0(W) = d(W)$ or $d(W)-1$. \hfill \break c) If $W$ has a
moment of order 1 and is centered, then $d_0(W) = d(W)$ and $D$ and
$L$ generate the same vector subspace.
\end{lem}
\proof The point {\it a)} is clear, since $\el_j = \el_1 \, \mod D,
\, \forall j \in J$. For {\it b)}, suppose $\Cal D^\perp \not =
\{\0\}$. There is $v_0 \in \Cal D^\perp$ such that $\langle v_0,
\el_1 \rangle = 1$. If $v $ is in $\Cal D^\perp$, then $\langle v -
\langle v, \el_1 \rangle \, v_0, \el_j\rangle = 0$, which implies $v
= \langle v, \el_1 \rangle \, v_0$, since $\Cal L = \R^d$.

{\it c)} If $v \in \Cal D^\perp$, then $\langle v, \el_j \rangle =
\langle v, \el_1 \rangle = \sum_i p_i \langle v, \el_1 \rangle =
\langle v, \sum_i p_i \el_i \rangle = 0, \forall j$; hence $v = \0$.
\eop

{\bf Recurrence/transience.} Recall that a r.w. $W = (Z_n)$ is
recurrent if and only if $\sum_{n=1}^\infty \PP(Z_n = \0) = +
\infty$. If $\sum p_j |\el_j| < \infty$, then $W$ is transient if
$\sum p_j \el_j \not = 0$ and, if $d=1$, recurrent if $\sum p_j
\el_j = 0$. These last two results are special cases of a general
result for cocycle over an ergodic dynamical system. A r.w. $W$ is
transient or recurrent according as $\Re e ({1 \over 1- \Psi})$ is
integrable on the $d$-dimensional unit cube or not (\cite{Sp64}).

The local limit theorem (LLT) (cf. Theorem \ref{LatLocThm}) implies
that a (reduced) r.w. $W$ with finite variance is recurrent if and
only if it is centered and $d(W) = 1$ or 2.

For further references, let us introduce the following condition for
a r.w. $W$:
\begin{eqnarray}
&&W \text{\it is reduced, has a moment of order 2 and is centered.} \label{hypR2}
\end{eqnarray}
If $W$ (reduced) has a moment of order 2, either it is centered with
$d(W) \leq 2$ (hence recurrent) or transient. In the first case
$d_0(W) = d$, in the second $d_0(W) = d$ or $d-1$.

{\bf Annulators of $L$ and $D$ in $\T^d$.}  The values of $\t \in
\T^d$ such that $\Psi(\t) = 1$ or $|\Psi(\t)| = 1$ play an important
role. They are characterized in the following lemma, where $B$ and
$B_1$ are the matrices introduced at the beginning of \ref{redRW}.
Recall that the annulator of a sublattice $L$ in $\Z^d$ is the
closed subgroup $\{\t \in \T^d: e^{2\pi i \langle \r, \t \rangle} =
1, \forall \r\ \in L \}$.
\begin{nota} \label{defdgamma} {\rm
We denote by $\Gamma$ (resp. $\Gamma_1$) the annulator of $L$ (resp. $D$),
by $d\gamma$ (resp. $d\gamma_1$) the Haar probability measure of the group
$\Gamma$ (resp. $\Gamma_1$).}
\end{nota}
\begin{lem} \label{modulus1} 1) We have $\Gamma = \{\t \in \T^d: \Psi(\t) = 1\} = (B^t)^{-1}
\Z^d \mod \Z^d$. The group $\Gamma$ is finite and $\Card(\Gamma) =
|\det B|$. \hfill \break 2a) We have $\Gamma_1 = \{\t \in \T^d:
|\Psi(\t)| = 1\} = \widetilde {\Cal D^\perp} + (B_1^t)^{-1} \Z^d
\mod \Z^d$ with $\widetilde {\Cal D^\perp} := {\Cal D^\perp} /
\Z^d$, which is either trivial or a 1-dimensional torus in $\T^d$.
The quotient $\Gamma_1 / \widetilde {\Cal D^\perp}$ is finite and
$a_0(W) := \Card (\Gamma_1 / \widetilde {\Cal D^\perp}) = |\det
B_1|$. \hfill \break 2b) If $\p \in D + n\el_1$, then the function
$F_{n, \p}(\t) := \Psi(\t)^n \ e^{- 2\pi i \langle \p, \t\rangle}$
is invariant by translation by elements of $\Gamma_1$. Its integral
is 0 if $\p \not \in D + n\el_1$. We have $|F_{n, \p}(\t) | < 1,
\forall \t \not \in \Gamma_1$, and $F_{n, \p}(\t) = 1, \forall \t
\in \Gamma_1$.
\end{lem}
\Proof We prove 2). The proof of 1) is analogous. We have
$|\Psi(\t)| \leq 1$ and by strict convexity, $|\Psi(\t)| = 1$ if and
only if $e^{2\pi i \langle \el_j, \t \rangle} = e^{2\pi i \langle
\el_1, \t \rangle}, \, \forall j \in J$, i.e., if and only if $\t
\in \Gamma_1$.

Recall that $D = B_1 \Z^d$ or $B_1 F_{d -1}$. Let us treat the
second case. By orthogonality, $\Cal D^\perp = (B_1^t)^{-1} {\Cal
E}_1$. If $\t \in \Gamma_1$, $e^{2\pi i \langle \r, \t \rangle} = 1,
\, \forall \r \in D$, so that, for $j= 2, ..., d$, $e^{2\pi i
\langle B_1 e_j, \t \rangle} = 1$, hence $\langle e_j, B_1^{t} \t
\rangle \in \Z$. It follows that the $d - 1$ last coordinates of
$B_1^{t} \t$ are integers. Therefore, $B_1^{t} \t \in {\Cal E}_1$
mod $\Z^d$, i.e., $\t \in (B_1^t)^{-1} {\Cal E}_1 + \Gamma_0 = \Cal
D^\perp + \Gamma_0$ mod $\Z^d$, where $\Gamma_0 := (B_1^t)^{-1} \Z^d
\mod \Z^d$ and $\widetilde {\Cal D^\perp} = {\Cal D^\perp}$ modulo
$\Z^d$ is a closed 1-dimensional subtorus of $\T^d$.

By Lemma \ref{ImLatt}.b, we have $\Card(\Gamma_0) = \Card (\Z^{d} /
B_1^t \Z^{d}) = |\det(B_1)|$. For a (reduced) centered r.w.,
$\Gamma_0 = \Gamma_1$ and $\Card (\Gamma_1) = |\det(B_1)|$.

2b) Since $F_{n, \p}(\t +\t_0) = e^{i \langle n \el_1 - \p, \,
\t_0\rangle} F_{n, \p}(\t)$ for $\t_0 \in \Gamma_1$, $F_{n,\p}$ is
invariant by all $\t_0 \in \Gamma_1$ if $\p \in D + n\el_1$, and its
integral is 0 if $\p \not \in D + n\el_1$.\eop

For a reduced r.w. $W$, there are $(\r_j, j \in J)$ in $\Z^d$ such
that $\el_j = B \, \r_j$ and the lattice generated by $(\r_j, j \in
J)$ is $\Z^d$. The r.w. is said to be {\it aperiodic} if $L = \Z^d$,
{\it strongly aperiodic} if $D = \Z^d$. If we replace $W$ by the
r.w. $W'$ defined by the r.v. $X_0'$ such that $\PP(X_0' = r_j) =
p_j$, we obtain an aperiodic r.w. This allows in several proofs to
assume aperiodicity without loss of generality when the r.w. is reduced
(for example see Theorem \ref{recd12} in \Sect \ref{selfInt}).

Strong aperiodicity is equivalent to $|\Psi(\t)| < 1$ for $\t \not =
\0 \mod \Z^d$ (cf. \cite{Sp64} and Lemma \ref{modulus1}). It is also
equivalent to: for every $\el\in \Z^d$, the additive group generated
by $\Sigma + \el$ is $\Z^d$. It implies $d_0(W) = d(W)$ and
$a_0(W)=1$. Observe that, contrary to aperiodicity, it is not always
possible to reduce proofs to the strictly aperiodic case.

A r.w. is ``deterministic" if $\PP(X_0 = \el_0) = 1$ for some
$\el_0$, so that $|\Psi(\t)| \equiv 1$ in this case.

\vskip 2mm {\bf Quadratic form}
\begin{lem} \label{quadrF} Let $W$ satisfy (\ref{hypR2}).
Let $Q$ be the quadratic form $Q(\u) = \Var (\langle X_0, \u
\rangle) = \sum_{j \in J} p_j \langle \el_j, \u\rangle^2$ and
$\Lambda$ the corresponding symmetric matrix. Then $Q$ is definite
positive. If $\Card (J)= d + 1$, then $\det(\Lambda) = c \,
\prod_{j=1}^{\Card (J)} p_j$, where $c$ does not depend on the
$p_j$'s.
\end{lem}
\Proof If $Q(\u)=0$, then $\langle X_0,\u\rangle$ is a.e. constant,
i.e., $\langle \el_i, \u \rangle= \langle\el_j, \u \rangle$ for all
$i,j\in J$; hence $\u = 0$, since (\ref{hypR2}) implies $\Cal D =
\R^d$.

Now, if $d' = \Card(J)$ is finite, let $J'$ be the set of indices
$\{2, ..., d'\}$. The quadratic form $q(\u) := \sum_{j \in J'} p_{j}
u_j^2 - (\sum_{j \in J'} p_{j} u_j)^2, \, \u \in \R^{d'-1}$, is
defined by the symmetric matrix $A = D_{d'} M$, where $D_{d'} :=
\diag(p_2, ..., p_{d'})$ and $M := \left(
\begin{matrix} 1 - p_2 & -p_3 & . & -p_{d'} \cr -p_2 & 1- p_3 & . &
-p_{d'} \cr . & . & . & . \cr -p_2 & -p_3& . & 1 -p_{d'}
\end{matrix} \right)$.

By subtracting line from line in the matrix $M$, we find $\det(M) =
1 - \sum_{j \in J'} p_j$ and therefore $\det(A) = \prod_{j=1}^{d'}
p_j$. The quadratic form $q$ is positive definite since
\begin{eqnarray*}
\sum_{j \in J'} p_j u_j^2 - (\sum_{j \in J'} p_j u_j)^2 - \sum_{2
\leq j' < j \leq d'} p_j p_{j'} (u_j - u_{j'})^2 = (1- \sum_{j \in
J'} p_j) \sum_{j \in J'} p_j \u_j^2.
\end{eqnarray*}
Let $U$ be the map from $\R^{d}$ to $\R^{d'-1}$: $\v \to U \v$,
where $U \v$ is the vector with coordinates $\langle \el_j - \el_1,
\v\rangle$, $j \in J'$. The
quadratic form $Q$ can be written $Q(\u) = \sum_{j \in J'} p_j
\langle \el_j - \el_1, \u\rangle^2 - (\sum_{j \in J'} p_j \langle
\el_j - \el_1, \u \rangle)^2 = q(U \u)$. Since $\Cal D = \R^d$, $U$
is injective, hence an isomorphism if and only if $\dim \Cal D = d =
d' -1$.

If $d = d' -1$, the determinant of $\Lambda=U^t A U$ is $c \,
\prod_{j=1}^{d+1} p_j$. The integer $c$ does not depend on the
probability vector $(p_j)_{j \in J}$. We have $c =\det(\tilde U)^2$ where, for
$j = 1, ..., d$, the matrix $\tilde U$ representing $U$ has as
$j$-{th} row the $d$ coordinates of the vector $\el_{(j+1)}-\el_1$. \eop

\goodbreak
\subsubsection{\bf Local limit theorem} \

Let $W = (Z_n)$ be a reduced random walk in $\Z^d$ associated to the
distribution $\PP(X_0 = \el_j) = p_j, j \in J$, with a finite second
moment. Recall that $a_0(W)$ and $\Lambda$ are defined in Lemmas
\ref{modulus1} and \ref{quadrF}. The local limit theorem (LLT) gives
an equivalent of $\PP(Z_n = \k)$ when $n$ tends to infinity:
\begin{thm} \label{LatLocThm} Suppose $W$ centered.
If $\k \not \in D + n\el_1$, then
$\PP(Z_n = \k) = 0$. If $\k \in D + n\el_1$, then we have
with $\sup_\k |\varepsilon_n(\k)| = o(1)$:
\begin{eqnarray}
(2 \pi n)^{{d\over 2}} \, \, \PP(Z_n = \k) = a_0(W) \, \det
(\Lambda)^{-\frac12} \, e^{- {1 \over 2n}\langle \Lambda^{-1} \, \k,
\, \k \rangle} + \varepsilon_n(\k). \label{errLLT}
\end{eqnarray}
\end{thm}
\begin{rems} \label{arithProgr}
1) If $(Z_n)$ is strongly aperiodic centered with finite second
moment, then $\lim_{n\to\infty} \, [(2\pi n)^{\frac{d}2} \,
\PP(Z_n=\k)] = \det(\Lambda)^{-\frac12}$ (\cite{Sp64}, P.10). A
version of the LLT in the centered case, extending a result of
P\'olya, was proved by van Kampen and Wintner (1939) in
\cite{vKWi39} without the strong aperiodicity assumption.

2) By Lemma \ref{quadrF}, if $J$ is finite and $\Card(J)= d + 1$,
then $\det(\Lambda) = c \prod_{j=1}^{\Card(J)} p_j$, where $c$ is an integer
$\geq 1$ not depending on the $p_j$'s.

3) It follows from (\ref{errLLT}) that, for every $\k$, there is
$N(\k)$ such that, for $n \geq N(\k)$, $\PP(Z_n = \k) > 0$ if and
only if $\k \in D + n\el_1$. It implies also $\sup_{\k \in \Z^d}
\PP(Z_n = \k) = O(n^{- {d \over 2}})$.

4) The condition $\k \in D + n\el_1$ reads $n\el_1 = \k$ mod $D$.
Since here $d_0(W) = d(W)$, $L/D$ is a cyclic finite group by Lemma
\ref{DLatt} and the set $E(\k) := \{n \geq 1: \k \in D + n\el_1\}$
for a given $\k$ is an arithmetic progression.

5) Let us mention a general version of the LLT which can be proved
under the assumption of finite second moment for a centered or non
centered r.w.

The quadratic form $Q$ reads now: $Q(\u) = \sum_{j \in J} p_j
\langle \el_j, \u\rangle^2 - (\sum_{j \in J} p_j \langle \el_j, \u
\rangle)^2$. If $\Lambda$ is the self-adjoint operator such that
$Q(\u) = \langle \Lambda \u, \u \rangle$, then $\Cal D$ is invariant
by $\Lambda$ and the restriction of $Q$ to $\Cal D$ is definite
positive. We denote by $\Lambda_0$ the restriction of $\Lambda$ to
$\Cal D$.
\begin{thm} \label{LLT2} For $\k \in \Z^d$, let $z_{n,\k} :=
{\k \over \sqrt n} - \sqrt n \, \E(X_0) = {\k \over \sqrt n} - \sqrt
n \sum_j p_j \el_j$. If $\k \not \in D + n\el_1$, then $\PP(Z_n =
\k) = 0$. If $\k \in D + n\el_1$, then, with $\sup_\k
|\varepsilon_n(\k)| = o(1)$:
\begin{eqnarray}
(2 \pi n)^{{d_0(W) \over 2}} \, \PP(Z_n = \k) = a_0(W) \, \det
(\Lambda_0)^{-\frac12} \, e^{- \frac12\langle
\Lambda_0^{-1}z_{n,\k}, \, z_{n,\k} \rangle} + \varepsilon_n(\k).
\label{errLLTGen}
\end{eqnarray}
If the moment of order 3 is finite, then
$\sup_\k |\varepsilon_n(\k)| = O(n^{-\frac12})$.
\end{thm}
In the non centered case, when $z_{n,\k}$ is not too large, for
example when $e^{- \langle \Lambda_0^{-1}z_{n,\k}, \, z_{n,\k}
\rangle} \geq 2 \sup_\k |\varepsilon_n(\k)|$, Equation
(\ref{errLLTGen}) gives an information on how the r.w. spreads
around the drift. With a finite third moment, it shows that, for any
constant $C$, for $n$ large, the r.w. takes with a positive
probability the values $\k$ such that $\k \in D + n \el_1$ which
belong to a ball of radius $C \sqrt n$ centered at the drift $n
\E(X_0)$.
\end{rems}
\subsubsection{\bf Upper bound for $\Phi_n(\omega) := \sup_{\el \in
\Z^d} R_n(\omega,\el) = \sup_{\el \in \Z^d} \sum_{k=0}^{n-1}
1_{Z_k(\omega) = \el}$.}

\

For the quenched CLT, we need some results on the local times and on
the number of self-intersections of a r.w. $W$.

\begin{prop} \label{BolthausenProp} (cf. \cite{Bo89}) a) If
the r.w. $W$ has a moment of order 2, then
\begin{eqnarray}
&&\for d= 1: \ \Phi_n(\omega) = o(n^{\frac12 + \varepsilon}); \ \for
d= 2: \ \Phi_n(\omega) = o(n^\varepsilon), \forall \varepsilon > 0.
\label{IneqepsilRn2}
\end{eqnarray}
b) If $W$ is transient with a moment of order $\eta$ for some $\eta
> 0$, then $\Phi_n(\omega) = o(n^\varepsilon), \forall \varepsilon >
0$.
\end{prop}
\Proof Let $\Cal A_n^{r}:= \{\el \in \Z^d: |\el| \le n^{r}\}$. If
$\E(\|X_0\|^{\eta}) < \infty$ for $\eta > 0$, then for $r\eta
> 1$, $\sum_{k=1}^\infty \PP(|X_k| > k^r) \leq (\sum_{k=1}^\infty
k^{-r\eta}) \E(\|X_0\|^{\eta}) < \infty$. By Borel-Cantelli lemma,
it follows $|X_k| \leq k^r$ a.s. for $k$ big enough. Therefore there
is $N(\omega)$ a.e. finite such that: $|X_0 + ... + X_{n-1}| \leq
n^{r+1}$ for $n > N(\omega)$. Hence $\sup_{\el \in \Z^d}
R_n^m(\omega,\el) = \sup_{\el \in \Cal A_n^{r+1}} R_n^m(\omega,\el)$
for $n \geq N(\omega)$.

{\it a)} We take $r=1$. For all $m \geq 1$ and for constants $C_m,
C_m'$ independent of $\el$, we have:
\begin{eqnarray}
&& \E[R_n^m(.,\el)] = \E[\sum_{k=0}^{n-1} 1_{Z_k = \el}]^m \le C_m
n^{m/2}, \for d=1, \and \le C_m' (\Log n)^m, \for d=2.
\label{majERnm}
\end{eqnarray}
To show (\ref{majERnm}), it suffices to bound $\sum_{0 \leq k_1 <
k_2 < ... < k_m < n} \PP(Z_{k_1} = \el, Z_{k_2} = \el, ..., Z_{k_m}
= \el)$. By independence and stationarity, with $F_n(k) =
1_{[0,n-1](k)} \PP(Z_k = \0)$ the preceding sum is
$$\sum_{0 \leq k_1 < k_2 < ... < k_m < n} \PP(Z_{k_1} = \el)\, \PP(Z_{k_2- k_1} =
\0) \, \cdots \, \PP(Z_{k_m- k_{m-1}} = \0) \leq C \sum_k (F_n * F_n
* ...* F_n)(k),$$
hence bounded by $C \bigl(\sum_k F_n(k)\bigr)^m$ and (\ref{majERnm})
follows from the bound given by the LLT.

By (\ref{majERnm}), $\E[\sup_{\el \in \Cal A_n^{2}}
R_n^m(\omega,\el)]$ is less than $n^2 \sup_{\el}
\E[R_n^m(\omega,\el)] \le C_m n^2 \cdot n^{m/2}$ for $d=1$ and less
than $n^4 \sup_{\el} \E[R_n^m(\omega,\el)] \le C_m n^4 (\Log \,
n)^m$ for $d=2$. Hence, for all $\varepsilon > 0$,
\begin{eqnarray*}
&&\for d=1, \ \ \sum_{n=1}^\infty \E[n^{-(\frac12+ \varepsilon)}
\sup_{\el \in \Cal A_n^{d}} R_n(\omega,\el)]^m \le C_m
\sum_{n=1}^\infty n^{2+m/2} / n^{(\frac12+\varepsilon) m} <\infty,
\text{ if } m > 3 /\varepsilon.\\
&&\for d=2, \ \ \sum_{n=1}^\infty \E[ n^{- \varepsilon}\, \sup_{\el
\in \Cal A_n^{d}} R_n(\omega,\el)]^m\, \le C_m \sum_{n=1}^\infty n^4
(\Log n)^m / n^{\varepsilon m}<\infty, \text{ if } m >
5/\varepsilon.
\end{eqnarray*}
Since $\sup_{\el \in \Z^d} R_n^m(\omega,\el) = \sup_{\el \in \Cal
A_n^{2}} R_n^m(\omega,\el)$ for $n \geq N(\varepsilon)$, this proves
{\it a)}.

{\it b)} According to $\sum_{k=0}^\infty \PP(Z_k = \p) < +\infty$,
using the same method as above, it can be shown that there are
constants $C_m$ and $M$ such that $\E[R_n^m(.,\el)]\le C_m M^{m}$.
 From the existence of a moment of order $\eta >0$, there is $r$ and
$N(\omega) < +\infty$ a.e. such that $\sup_{\el \in \Z^d}
R_n^m(\omega,\el) = \sup_{\el \in \Cal A_n^{r+1}} R_n^m(\omega,\el)$
for $n \geq N(\omega)$.

In view of $\E[\sup_{\el \in \Cal A_n^{r}} R_n^m(\omega,\el)] \leq
n^{d(r+1)} \sup_{\el} \E[R_n^m(\omega,\el)] \le C_m n^{d(r+1)} \,
M^{m}$, for every $\varepsilon >0$ there is $m$ such that
$\sum_{n=1}^\infty \E[ n^{- \varepsilon}\, \sup_{\el \in \Cal
A_n^{r}} R_n(\omega,\el)]^m < \infty$. This implies {\it b)}. \eop

\subsubsection{\bf Self-intersections} \label{subsubSelf}

\

{\bf 1) Recurrent case}

The proof of the following result is postponed to \Sect
\ref{selfInt}.
\begin{thm} \label{recd12} Let $W$ satisfy (\ref{hypR2}) and let $\p$ be in
$L(W)$. If $d(W)= 1$ or 2, $\displaystyle \lim_n {V_{n, \p}(\omega)
\over V_{n}(\omega)} = 1$ a.e. If $d(W)=2$, $V_{n, \p}$ satisfies a
SLLN: $\lim_n V_{n, \p}(\omega) / \E V_{n, \p} = 1$ a.e.
\end{thm}

{\bf 2) Transient case}

In the transient case (without moment assumptions), a general
argument is available:
\begin{lem} \label{ergLem} If $(\Omega, \PP, \tau)$ is an ergodic
dynamical system and $(f_k)_{k \geq 1}$ a sequence of functions in
$L^1(\Omega, \PP)$ such that $\sum_{k \geq 1} \|f_k\|_r < \infty$,
for some $r > 1$, then
\begin{eqnarray}
\lim_n \frac1n \sum_{k=1}^{n-1} \sum_{\ell=0}^{n-k-1} f_k(\tau^\ell
\omega) = \sum_{k=1}^\infty \int f_k \, d\PP, \text{ for a.e. }
\omega. \label{ergLemFormula}
\end{eqnarray}
\end{lem}
\Proof \ We can assume $f_k \geq 0$. For the maximal function
$\tilde f_{k}(\omega):= \sup_{n \geq 1} \frac1n \sum_{\ell=0}^{n-1}
\, f_{k}(\tau^\ell \omega)$, by the ergodic maximal lemma, there is
a finite constant $C_r$ such that $\|\tilde f_{k}\|_r \leq C_r
\|f_{k}\|_r$. Therefore $\sum_{k=1}^\infty \, \tilde f_{k} \in
L^r(\Omega, \PP)$; hence $\sum_{k=1}^\infty \, \tilde f_{k} (\omega)
< +\infty$, for a.e. $\omega$.

Let $\omega$ be such that $\sum_{k=1}^\infty \, \tilde f_{k}
(\omega) < +\infty$. For $\varepsilon > 0$, there is $L$ such that
$\sum_{k > L} \int f_{k} \, d\PP \leq \sum_{k > L} \|f_{k}\|_r \leq
\varepsilon$ and $\sum_{k > L} \tilde f_{k}(\omega) \leq
\varepsilon$; hence, uniformly in $n$, $\frac1n \sum_{k= L +1}^ n
\sum_{\ell=0}^{n-k-1} f_k(\tau^\ell \omega) \leq \varepsilon$. By
the ergodic theorem, we have $\lim_n \frac1n \sum_{k=1}^L \,
\sum_{\ell =0}^{n-k-1} f_k(\tau^\ell \omega) = \sum_{k=1}^L \, \int
f_k \, d\PP$. Therefore, for $n$ big enough: $|\frac1n \sum_{k=1}^n
\, \sum_{\ell =0}^{n-k-1} f_k(\tau^\ell \omega) - \sum_{k=1}^\infty
\, \int f_k \, d\PP | \leq 2 \varepsilon$. \eop

Recall that $\Psi(\t) = \E [e^{2\pi i \langle X_0, \t\rangle}], \ \t
\in \T^d$. We use the notation:
\begin{eqnarray}
w(\t) := {1 - |\Psi(\t)|^2 \over |1 - \Psi(\t)|^2}, \ c_w =
\int_{\T^d} w \, d\t, \text{ when } w \text{ is integrable}.
\label{defWt}
\end{eqnarray}
Outside the finite group $\Gamma$, $w$ is well defined, nonnegative
and $w(\t) = 0$ only on $\Gamma_1 \setminus \Gamma$. Hence it is
positive for a.e. $\t$, excepted when the r.w. is ``deterministic".

\begin{proposition} \label{renewProp} (\cite{Sp64})
Let $W =(Z_n)$ be a transient random walk in $\Z^d$.

a) The function $w$ is integrable on $\T^d$ and there is a
nonnegative constant $K$ such that, if $W$ is aperiodic,
\begin{eqnarray}
&&I(\p) := 1_{\p= \0} + \sum_{k = 1}^{\infty} \, [\PP(Z_k = \p) +
\PP(Z_k = - \p)] = \int_{\T^d} \cos 2\pi \langle \p, \t\rangle \,
w(\t) \, d\t + K. \label{renewval}
\end{eqnarray}
For a general reduced r.w. W, the sums $I(\p)$ are the
Fourier coefficients of the measure $w d\t + K d\gamma$.

b) If $d > 1$, then $K = 0$; if $d = 1$ and $m(W):=\sum_{\ell \in
\Z} \, \PP(X_0 = \ell) \, |\ell| < + \infty$, then $K = |\sum_{\ell
\in \Z} \, \PP(X_0 = \ell) \, \ell|^{-1}$; if $d = 1$ and $m(W) =
+\infty$, then $K =0$.
\end{proposition}
\proof For completeness we give a proof of a), following Spitzer (\S
9, P2 in \cite{Sp64}). For b) a more difficult result without
assumption on the moments need to be used (see Spitzer \S 24, P5,
P6, P8, T2 in \cite{Sp64}).

Let $(Z_n)$ be a transient, reduced, aperiodic r.w. Since $\sum_{k =
1}^{\infty} \, \PP(Z_k = \p) < +\infty$, we have
\begin{eqnarray*}
&&\sum_{k = 0}^{\infty} \, \lambda^k \PP(Z_k = \p) = \int_{\T^d}
e^{2\pi i \langle \p, \t\rangle} \, {1 \over 1 - \lambda \Psi(\t)}
\, d\t, \, \forall \lambda \in [0, 1[, \\
&&\sum_{k = 0}^{\infty} \, [\PP(Z_k = \p) + \PP(Z_k = - \p)] = 2
\lim_{\lambda \uparrow 1} \int_{\T^d} \cos 2\pi \langle \p,
\t\rangle \, \Re e ({1 \over 1 - \lambda \Psi(\t)}) \, d\t;
\end{eqnarray*}
therefore, $I(\p) = \lim_{\lambda \uparrow 1} \int_{\T^d} \cos 2\pi
\langle \p, \t\rangle \, w_\lambda(\t)  \, d\t$, with $\displaystyle
w_\lambda(\t) := {1 - \lambda^2 |\Psi(\t)|^2 \over |1 - \lambda
\Psi(\t)|^2}$.

Taking $\p = \0$ in the previous formula, we deduce, by Fatou's
lemma, the integrability of $w$ on $\T^d$ and with a nonnegative
constant $K$ the equality:
\begin{eqnarray*}
&&I(\0) = 1 + 2 \sum_{k = 1}^{\infty} \, \PP(Z_k = \0) =
\lim_{\lambda \uparrow 1} \int_{\T^d} \, w_\lambda(\t) \, d\t =
\int_{\T^d} w(\t) \, d\t  + K = c_w +K,
\end{eqnarray*}
since $w_\lambda(\t) \geq 0$, $\lim_{\lambda \uparrow 1}
w_\lambda(\t) = {1 - |\Psi(\t)|^2 \over |1 - \Psi(\t)|^2} = w(\t)$,
for $\t \in \T^d \stm0$.

By aperiodicity of $W$, $\Psi(\t) \not = 1$ for $\t \in U_\eta^c$, where
$U_\eta^c$ is the complementary in $\T^d$ of the
ball $U_\eta$ of radius $\eta$ centered at $\0$, for
$\eta > 0$. This implies
$\sup_{\t \in U_\eta^c} \sup_{\lambda < 1} w_\lambda(\t) < +\infty$.

Therefore, we get: $\displaystyle \lim_{\lambda \uparrow 1}
\int_{U_\eta^c} \cos 2\pi \langle \p, \t\rangle \,  w_\lambda(\t) \,
d\t = \int_{U_\eta^c} \cos 2\pi \langle \p, \t\rangle \,  w(\t) \,
d\t$, hence:
\begin{eqnarray}
I(\p) = \int_{U_\eta^c} \cos 2\pi \langle \p, \t\rangle \, w(\t)\,
d\t + \lim_{\lambda \uparrow 1} \int_{U_\eta} \cos 2\pi \langle \p,
\t\rangle \, w_\lambda(\t) \, d\t, \, \forall \eta
>0. \label{IpLim}\end{eqnarray}
Let $\varepsilon > 0$. By positivity of $w$, we have, for
$\eta(\varepsilon)$ small enough:
\begin{eqnarray*}
&&(1 - \varepsilon) \int_{U_{\eta(\varepsilon)}} \, w_\lambda(\t) \,
d\t \leq \int_{U_{\eta(\varepsilon)}} \cos 2\pi \langle \p,
\t\rangle \, w_\lambda(\t) \, d\t \leq  (1+\varepsilon)
\int_{U_{\eta(\varepsilon)}} \,  w_\lambda(\t) \, d\t;\\
&&\text{hence, using (\ref{IpLim}): } (1 - \varepsilon)
\int_{U_{\eta(\varepsilon)}} \, w_\lambda(\t) \, d\t -
\int_{U_{\eta(\varepsilon)}} \cos 2\pi \langle \p, \t\rangle \,
w(\t)\, d\t  \\&& \leq I(\p) - \int_{\T^d} \cos 2\pi \langle \p,
\t\rangle \, w(\t)\, d\t \leq (1 + \varepsilon)
\int_{U_{\eta(\varepsilon)}} \, w_\lambda(\t) \, d\t -
\int_{U_{\eta(\varepsilon)}} \cos 2\pi \langle \p, \t\rangle \,
w(\t)\, d\t.
\end{eqnarray*}
For $\varepsilon$ small enough, $\int_{U_{\eta(\varepsilon)}} \cos
2\pi \langle \p, \t\rangle \, w(\t)\, d\t$ can be made arbitrary
small, since $w$ is integrable, as well as $\varepsilon
\sup_{\lambda < 1} \int_{U_\eta} w_\lambda \, d\t$, since
$\sup_{\lambda < 1} \int_{\T^d} w_\lambda \, d\t < \infty$.
Therefore the value of $I(\p) - \int_{\T^d} \cos 2\pi \langle \p,
\t\rangle \, w(\t)\, d\t$ does not depend on $\p$, hence is equal to
$K$. This shows (\ref{renewval}). The Riemann-Lebesgue Lemma implies
$K = \lim_{|\p| \to +\infty} \sum_{k = 1}^{\infty} \, \PP(Z_k = \pm
\p)$.

In the aperiodic case, by (\ref{renewval}) the sums $I(\p)$ are the
Fourier coefficients of the measure $w d\t + K \delta_0$. In the
general case, by replacing $\el_j$ by $\r_j$ with $\el_j = B
 \r_j$ and $\t$ by $B^{t} \t$ (cf. Subsection \ref{prelimRW}), we get an
aperiodic r.w. and we apply the previous result. The push forward of
the measure $\delta_0$ by $B^t$ is the measure $d\gamma$, which
shows a).\eop

\goodbreak
\subsection{\bf $\zeta$-regularity, variance}
\subsubsection{\bf Regularity of summation sequences defined by r.w.}

\

We have $V_{n, \p} = 0$, for $\p \not \in L(W)$. If $W$ satisfies
(\ref{hypR2}), by the LLT (Theorem \ref{LatLocThm}), for $\p \in
L(W)$, with constants $C_i$ independent of $\p$:
\begin{eqnarray}
&&\E V_{n, \p} \sim C_1 \, n^{\frac32} \for d =1, \ \E V_{n, \p}
\sim C_2 \, n \, \Log \, n \for d = 2, \ \E V_{n, \p} \sim C_d \, n
\for d >2. \label{EVnp}
\end{eqnarray}
It follows $\lim_n \E V_{n, \p} / \E V_{n, \0} = 1, \forall \p \in
L(W)$. Let us check (\ref{EVnp}) for $d=2$. Let $C$ be the constant
$a_0(W) \, \det (\Lambda)^{-\frac12} \, (2 \pi)^{{-1}}$. Recall that
$L/D$ is a finite cyclic group (each $\ell_j$ is a generator of the
cyclic group $L/D$, cf. Lemma \ref{DLatt}).

By (\ref{errLLT}), for every $\varepsilon > 0$, there is
$K_\varepsilon$ such that:
\begin{eqnarray*}
&&|\sum_{n=1}^N \PP(Z_n = \p) - C \, \sum_{n=1}^N  {1_{n\el_1 = \p \mod D} \over n}| \\
&&\leq  C \, \sum_{n=1}^N \ 1_{n\el_1 = \p \mod D}{|e^{- {1 \over
2n}\langle \Lambda^{-1} \, \p, \, \p \rangle} -1 +
\varepsilon_n(\p)| \over n} \leq C (K_\varepsilon + \varepsilon \,
\Log \, N).
\end{eqnarray*}

Since $(\Log N )^{-1} \, \sum_{n=1}^N {1_{n\el_1 = \p \mod D} \over
n} \to (\Card \, L/D)^{-1}$ we have:
\begin{eqnarray}
&&  \lim_N \, [(\Log N )^{-1} \, \sum_{n=1}^N \PP(Z_n = \p)] =  C \,
\lim_N \, [(\Log N )^{-1} \, \sum_{n=1}^N  {1_{n\el_1 = \p \mod D}
\over n}] = {C \over \Card \, L/D}. \label{limPotLogp}
\end{eqnarray}

 From (\ref{limPotLogp}) and (\ref{EspVnpErgSum}), we have in
dimension 2, for every $\p \in L$:
\begin{eqnarray*}
&&{\E V_{N, \p} \over N \Log \, N} = {1_{\p= 0} + \sum_{k = 1}^{N-1}
\,(1-{k \over N}) \, [\PP(Z_k = \p) + \PP(Z_k = -\p)] \over \Log \,
N} \to {2 \, C \over \Card \, L/D}.
\end{eqnarray*}

Recall that $\zeta$-regularity of a rws $(R_n(\el, \omega))_{n \geq
1}$ is equivalent to $\displaystyle \lim_n {V_{n,\p}(\omega) \over
V_{n,\0}(\omega)} = \hat \zeta (\p), \forall \p \in \Z^d$, for a.e.
$\omega$.

\begin{thm} \label{transrecSpec} 1) If the rws is defined by a r.w.
$W$ which satisfies (\ref{hypR2}) with $d(W) \leq 2$ (hence centered
and recurrent), then it is $\zeta$-regular with $\zeta = d\gamma$ (=
$\delta_0$ if $W$ is strictly aperiodic). \hfill \break - For $d(W)
=1$, the normalization $V_{n,\0}(\omega)$ depends on $\omega$.
\hfill \break - For $d(W) =2$, the normalization is $C \, n \, \Log
\, n$, with $C = \pi^{-1} a_0(W) \, \det(\Lambda)^{-\frac12}$.

2) If the rws is defined by a transient r.w., then it is
$\zeta$-regular and the normalization is $C \, n$, with $C = 1 + 2
\sum_{k = 1}^{\infty} \, \PP(Z_k = \0)$.

If $d = 1$, then $C = c_w +K$ and $d\zeta(\t) = (c_w + K)^{-1}\,
(w(\t) \, d\t + d\gamma(t))$, where $c_w$ is defined in
(\ref{defWt}) and $K$ is the constant given by Proposition
\ref{renewProp} ($K \not = 0$ if $m(W)$ is finite). If $d \geq 2$,
then  $C = c_w$ and $d\zeta(\t) = c_w^{-1} \, w(\t) \, d\t$.
\end{thm} \proof 1) The recurrent case follows from Theorem \ref{recd12}.

As $V_n(\omega)/ \E V_n \to 1$, we choose the constant $C$ such that
$\lim_n \E V_n / C n \Log n = 1$, i.e., by the LLT: $C = 2 \lim_n
(\sum_{k = 1}^{n-1} \, \PP(Z_k = \0) /\Log n) = \pi^{-1} a_0(W) \,
\det (\Lambda)^{-\frac12}$.

2) In the transient case, the Fourier coefficients of the kernel
$\frac1n \underline R_n^\omega$ given by (\ref{FourKern}) converge
by Lemma \ref{ergLem} to the finite sum of the series $1_{\p= \0} +
\sum_{k = 1}^{\infty} \, [\PP(Z_k = \p) + \PP(Z_k = - \p)]$, which
are the Fourier coefficients of the measure $w d\t + K d\gamma$
according to Proposition \ref{renewProp}. \eop

\begin{rems} {\rm 1) With Condition (\ref{hypR2}), if $\Card(J)
= d+1$, then $\det(\Lambda) = c \, \prod p_j$, where $c$ is an integer
$\geq 1$ depending only on the $\el_j$'s (cf. Lemma \ref{quadrF}). In the
strictly aperiodic case, we obtain $C = \pi^{-1} \, (c \prod p_j)^{-\frac12}$.

2) For $d \not = 1$, the asymptotic variance w.r.t. $\PP \times \mu$
is the same as the quenched variance.

3) For a reduced r.w., $\zeta$ is a discrete measure only in the
centered case when $d \leq 2$, and in the non centered case when $d=
1$ and the r.w. is deterministic.

4) As an example, let us consider a r.w. $W$ on $\Z$ defined by
$\PP(X_k = \ell_j) = p_j$, with $\ell_1 = 1$ and $0 < p_1 < 1$.
Suppose that $W$ is transient. Let $S$ be a map with Lebesgue
spectrum. The random sequence $(S^{Z_k(\omega)})$ can be viewed as a
random commutative perturbation of the iterates of $S$. If
$\varphi_f$ is continuous, the asymptotic variance of ${1 \over
\sqrt n} \|\sum_{k=0}^{n-1} f(S^{Z_k(\omega)}.)\|_2$ is $\int
\varphi_f \, d\zeta$. By Theorem \ref{transrecSpec}, the variance is
$\not = 0$ if $f \not \equiv 0$. A CLT for the quenched ergodic sums
can be shown when $S$ is of hyperbolic type.} \end{rems}

\subsubsection{\bf Regularity of barycenter summations} \label{barycprocsubsec}

\

Let $W = (Z_n)$ be a reduced r.w. in $\Z^d$ and let $\tilde W =
(\tilde Z_n)$ be the symmetrized r.w. If $P$ is the barycenter
operator defined as in (\ref{barycproc}) by $Pf = \sum_j p_j
T^{\el_j} f$, we have
$$P^n f = \sum_{\el \in \Z^d} R_n(\el) \ T^\el f, \text{ with }
R_n(\el) = \PP(Z_n = \el)$$ and by (\ref{chch}) $(R_n * \check
R_n)(\el) = \PP(\tilde Z_n = \el)$. We have $d_0(\tilde W) =
d(W)-\delta = d_0(W)$ with $\delta = 0$ or 1. The behaviour of $P^n$
is given by the LLT applied to $\tilde W$. When $\delta = 0$,
$\Gamma_1$ is finite and its probability Haar measure $d\gamma_1$ is
a discrete measure. If $\delta = 1$, $d\gamma_1$ is the product of a
discrete measure by the uniform measure on a circle. The kernel
$\tilde R_n$ is $\sum_\el {\PP(\tilde Z_n = \el)\over \PP(\tilde Z_n
= 0)} \, e^{2\pi i \langle \el , \, \t\rangle}$ and we have
\begin{eqnarray*} \|P^n f\|_2^2 &&= \int_X |\E(f(A^{Z_n(.)}
x))|^2 \, d\mu(x) = \sum_{\el \in \Z^d} \PP(\tilde Z_n = \el) \,
\hat \varphi_f(\el).
\end{eqnarray*}
\begin{proposition} \label{regBaryc} The kernel $(\tilde R_n)$ (with a normalization
factor $\sum_\el R_n(\el)^2$ of order $n^{d-\delta}$) weakly
converges to the measure $d\gamma_1$.
\end{proposition}
\Proof The characteristic function of $\tilde X_0$ is $|\Psi(\t)|^2
= |\E(e^{2\pi i \langle X_0, \t\rangle})|^2$. The LLT given by
Theorem \ref{LLT2} ($\tilde W$ is centered but non necessarily
reduced) implies
$$\PP(\tilde Z_n = \el) = (2 \pi n)^{-{d_0(\tilde W) \over 2}} \, [a_0 \, \det
(\tilde\Lambda_0)^{-\frac12} \, e^{-{1 \over 2 n} \langle
\tilde\Lambda_0^{-1} \el, \, \el \rangle} + \varepsilon_n(\el)].$$

For the r.w. $(\tilde Z_n)$, the condition $\k \in D + n\el_1$
reduces to $\k \in D$, since $\0$ belongs to the support of the
symmetrized distribution. Therefore, for the symmetric r.w., we
obtain:
\begin{eqnarray}
\lim_n {\int |\Psi(\t)|^{2n} \, e^{-2\pi i \langle \p, \t \rangle}
d\t \over \int |\Psi(\t)|^{2n} \, \, d\t}= \lim_n {\PP(\tilde Z_n =
\p) \over \PP(\tilde Z_n = \0)} \to 1_D(\p), \forall \p \in \Z^d.
\label{ratioSymm}
\end{eqnarray}

Remark that when $W$ is strictly aperiodic, one can show that
$\bigl(|\Psi(\t)|^{2n}/\int |\Psi(\t)|^{2n} d\t\bigr)_{n \geq 1}$ is
an approximation of identity. This gives a direct proof of
(\ref{ratioSymm}).

Since $1_D$ is the Fourier transform of $d\gamma_1$ viewed as a
measure on the torus $\T^d$, the weak convergence of $(\tilde R_n)$
to $d\gamma_1$ follows from (\ref{ratioSymm}). \eop

\vskip 3mm
\section{\bf $\Z^d$-actions by commuting endomorphisms on
$G$} \label{endoGroup}

\subsection{\bf Endomorphisms of a compact abelian group $G$}
\label{ssectPrelim}

\

Let $G$ be a compact abelian group with Haar measure $\mu$. The
group of characters of $G$ is denoted by $\hat G$ or $H$ and the set
of non trivial characters by $\hat G^*$ or $H^*$. The Fourier
coefficients of a function $f$ in $L^1(G, \mu)$ are $c_f(\chi) :=
\int_G \, \overline \chi \, f \, d\mu$, $\chi \in \hat G$.

Every surjective endomorphism $B$ of $G$ defines a measure
preserving transformation on $(G, \mu)$ and a dual injective
endomorphism on $\hat G$. For simplicity, we use the same notation
for the actions on $G$ and on $\hat G$. If $f$ is function on $G$,
$B f$ stands for $f \circ B$.

We consider a semigroup ${\Cal S}$ of surjective commuting
endomorphisms of $G$, for example the semigroup generated by
commuting matrices on a torus. It will be useful to extend it to a group of
automorphisms acting on a (possibly) bigger group $\tilde G$ (a
solenoidal group when $G$ is a torus). Let us recall
briefly the construction.
\begin{lem} \label{embedd} There is a smallest compact abelian group
$\tilde G$ (connected, if $G$ is connected) such that $G$ is a
factor of $\tilde G$ and ${\Cal S}$ is embedded in a group $\tilde
{\Cal S}$ of automorphisms of $\tilde G$.\end{lem} \Proof On the set
$\{(\chi, A), \, \chi \in H, A \in \Cal S\}$ we consider the law
$(\chi, A) + (\chi', A') = (A' \chi + A \chi', AA')$). Let $\tilde
H$ be the quotient by the equivalence relation $\Cal R$ defined by
$(\chi, A) \, \Cal R \, (\chi', A')$ if and only if $A'\chi =
A\chi'$.

The transitivity of the relation $\Cal R$ follows from the
injectivity of the dual action of each $A \in \Cal S$. The map $\chi
\in \hat G \ \to \ (\chi, Id) / \Cal R$ is injective. The elements
$A \in \Cal S$ act on $\tilde H$ by $(\chi, B) / \Cal R \ \to \
(A\chi, B) / \Cal R$. The equivalence classes are stable by this
action. We can identify $\Cal S$ and its image. For $A \in \Cal S$,
the automorphism $(\chi, B) / \Cal R \ \to \ (\chi, AB) / \Cal R $
is the inverse of $(\chi, B) / \Cal R \ \to \ (A\chi, B) / \Cal R $.

Let $\tilde G$ be the compact abelian group (with
Haar measure $\tilde \mu$) dual of the group $\tilde H$ endowed with the discrete
topology. The group $H = \hat G$ is isomorphic to a subgroup of
$\tilde H$ and $G$ is a factor of $\tilde G$. We obtain an embedding
of ${\Cal S}$ in a group $\tilde {\Cal S}$ of automorphisms of
$\tilde H$ and, by duality, in a group of automorphisms of $\tilde
G$. If $\hat G$ is torsion free, then $\tilde H$ is also torsion
free and its dual $\tilde G$ is a connected compact abelian group.
\eop

If $(A_1, ..., A_d)$ are $d$ algebraically independent automorphisms
of $\tilde G$ generating $\tilde {\Cal S}$, then $\tilde {\Cal S} =
\{A^\el = A_1^{\ell_1}... A_d^{\ell_d}, \el = (\ell_1, ..., \ell_d)
\in \Z^d\}$ and $\tilde {\Cal S}$ is isomorphic to $\Z^d$ if it is
torsion free. The corresponding $\Z^d$-action on $(\tilde G, \tilde \mu)$
is denoted $\A$.
\begin{assumption} \label{hypo1} {\it In what follows, we consider a set
$(B_j, j \in J)$ of commuting surjective endomorphisms of $G$ and
the generated semigroup $\Cal S$. Denoting by $\tilde {\Cal S}$ the
extension of $\Cal S$ to a group of automorphisms of $\tilde G$ (if
the $B_j$'s are invertible, then $\tilde G$ is just $G$), we assume
that $\tilde {\Cal S}$ is torsion-free and non trivial, hence has a
system of $d \in [1, +\infty]$ algebraically independent generators
$A_1, ..., A_d$ (not necessarily in $\Cal S$). We suppose $d$
finite. In other words we consider a set ($B_j, j \in J)$ of
endomorphisms such that $B_j = A^{\el_j}$, where $A_1, ..., A_d$ are
$d$ algebraically independent commuting automorphisms of $\tilde
G$.}
\end{assumption}
We begin with some spectral results. The measure preserving
$\Z^d$-action $\A$ is said to be {\it totally ergodic} if $A^\el$ on
$(\tilde G, \tilde \mu)$ is ergodic for every $\el \in \Z^d \stm0$.

One easily show that total ergodicity is equivalent to: $A^\el \chi \not = \chi$
for any non trivial character $\chi$ and $\el \not = \underline 0$
(free dual $\Z^d$-action on $H^*$), to the Lebesgue spectrum property
for $\tilde {\Cal S}$ acting on $(\tilde G, \tilde \mu)$, as well as to
$2$-mixing, i.e., $\lim_{\|\el\| \to \infty} \mu(C_1 \cap A^{-\el}
\, C_2)= \mu(C_1) \, \mu(C_2)$, for all Borel sets $C_1, C_2$ of $G$.

For a semigroup $\Cal S = \{A^\el, \el \in (\Z^+)^d\}$, total ergodicity of the generated group
is equivalent to the property (expressed on the dual of $G$): $A^\el \chi
\not = A^{\el'} \chi$ for $\chi \in H^*$ and $\el \not = \el'$.

Examples for $G = \T^\rho$ are given in Subsection \ref{torEx}. What
we call ``totally ergodic" for a general compact abelian group $G$
is also called ``partially hyperbolic" for the torus.

\vskip 2mm {\bf Spectral density $\varphi_f$ of a regular function $f$}

For $f$ in $L^2(\tilde G)$, we have $\hat \varphi_f(\n) = \langle f,
A^\n f\rangle = \sum_{\chi \in \tilde H} \ c_f(A^\n \chi) \,
\overline{c_f(\chi)}, \n \in \Z^d$. Observe that, using the
projection $\Pi$ from $\tilde G$ to $G$, a function $f$ on the group
$G$ can be lifted to $\tilde G$ and a character $\chi \in \hat G$
viewed as a character on $\tilde G$ via composition by $\Pi$.
Putting $\tilde f(x) = f(\Pi x)$, we have $\tilde f = f\circ \Pi =
\sum_{\chi \in \hat G} c_{f}(\chi) \, \chi \circ \Pi$, so that the
only non zero Fourier coefficients of $\tilde f$ correspond to
characters on $G$.

If $f$ is a function defined on $G$, the Fourier analysis can be
done for the $\Z^d$-action $\A$ in the group $\tilde G$, but
expressed in terms of the Fourier coefficients of $f$ computed in
$G$. For the spectral density of $f$, viewed as the spectral density
of $\tilde f$ for the $\Z^d$-action $\A$ on $\tilde G$, one checks
that $\hat \varphi_f(\n) = \sum_{\chi \in \hat G} \ c_f(A^\n \chi)
\, \overline{c_f(\chi)}, \n \in \Z^d$, where $c_f(A^\n \chi)= 0$ if
$A^\n \chi \not \in \hat G$.

For the action by endomorphisms of compact abelian groups, the
family $(\psi_j)$ (cf. notations in Subsection \ref{LebSpec}) is
$(\chi_j)$. We have: \hfill \break $a_{j,\el} = \langle f, \ A^\el
\chi_j\rangle = c_f(A^{\el}\, \chi_j), \ \gamma_j(\t) = \sum_{\el
\in \Z^d} c_f(A^\el \, \chi_j) \, e^{2 \pi i \langle \el, \t
\rangle}, \ M_\t f = \sum_{j\in J_0} \gamma_j(\t) \, \chi_j$, where
$J_0$ denotes a {\it section} of the action of $\tilde {\Cal S}$ on
$\tilde H^*$, i.e., a subset $\{\chi_j\}_{j \in J_0} \subset \tilde
H^* = \tilde H \stm0$ such that every $\chi \in \tilde H^*$ can be
written in a unique way as $\chi = A_1^{n_1}... A_d^{n_d} \,
\chi_j$, with $j \in J_0$ and $(n_1, ..., n_d) \in \Z^d$.

Therefore, $M_\t f$ is defined for every $\t$, if $\sum_{j \in
J_0}\big(\sum_{\el \in \Z^d} |c_f(A^{\el}\, \chi_j)|^2) < \infty$.

We denote by $AC_0(G)$ the class of real functions on $G$ with {\it
absolutely convergent Fourier series} and $\mu(f) =0$, endowed with
the norm: $\|f \|_c := \sum_{\chi \in \hat G} |c_f(\chi)| <
+\infty$.

\begin{proposition} \label{condAbsConv} If $f$ is in $AC_0(G)$,
then $\sum_{\el \in \Z^d} |\langle A^{\el}f, f\rangle| < \infty$ and
the spectral density $\varphi_f$ of $f$ is continuous. For any
subset $\Cal{E}$ of $\hat G$, if $f_1(x) = \sum_{\chi \in \Cal{E}}
c_f(\chi) \, \chi$, then
\begin{eqnarray}
\|\varphi_{f-f_1}\|_\infty \leq \|f-f_1\|_c^2. \label{densspectAC}
\end{eqnarray}
\end{proposition}
\proof (We work in $\tilde G$, but for simplicity we do not write
$\tilde~$). Since every $\chi \in H^*$ can be written in an unique
way as $\chi = A^{\el}\, \chi_j$, with $j \in J_0$ and $\el \in
\Z^d$, we have: $\sum_{\el \in \Z^d} |c_f(A^{\el}\, \chi_j)| \leq
\sum_{j \in J_0} \sum_{\el \in \Z^d} |c_f(A^{\el}\, \chi_j)| =
\sum_{\chi \in \tilde H^*} |c_f(\chi)| = \|f \|_{c}$.

Then, for every $j \in J_0$, the series defining $\gamma_j$ is
uniformly converging, $\gamma_j$ is continuous and $\sum_{j \in J_0}
\|\gamma_j\|_\infty \leq \sum_{j \in J_0} \sum_{\el \in \Z^d}
|c_f(A^{\el}\, \chi_j)| = \|f \|_{c}$.

The spectral density of $f$ is $\varphi_f(\t) = \sum_{j \in J_0} \,
|\sum_{\el \in \Z^d \stm0} c_f(A^\el \chi_j) \, e^{2\pi i
\langle\el, \t\rangle} \, |^2$. The function $\sum_{j \in J_0}
|\gamma_j|^2$ is a continuous version of the spectral density
$\varphi_f$ and $\|\varphi_f\|_\infty$ is bounded by
\begin{eqnarray*}
&&\sum_{j \in J_0} \, (\sum_{\el \in \Z^d} |c_f(A^\el \chi_j)|)^2
\leq \sum_{j \in J_0} \, (\sum_{\el \in \Z^d} |c_f(A^\el \chi_j)|)
(\sum_{\chi\in \hat G} |c_f(\chi)|) \leq (\sum_{\chi\in \hat G}
 \, |c_f(\chi)|)^2 = \|f \|_{c}^2.
\end{eqnarray*}
Inequality (\ref{densspectAC}) follows by replacing $f$ by $f -
f_1$. \eop

\subsection{\bf The torus case: $G = \T^\rho$, examples of $\Z^d$-actions}
\label{torEx}

\

Every $B$ in the semigroup ${\Cal M}^*(\rho, \Z)$ of non singular
$\rho \times \rho$ matrices with coefficients in $\Z$ defines a
surjective endomorphism of $\T^\rho$ and a measure preserving
transformation on $(\T^\rho, \mu)$. It defines also a dual
endomorphism of the group of characters $H = \widehat {\T^\rho}$
identified with $\Z^\rho$ (action by the transposed of $B$). Since
we compose commuting matrices, for simplicity we do not write the
transposition. When $B$ is in the group $GL(\rho, \Z)$ of matrices
with coefficients in $\Z$ and determinant $\pm 1$, it defines an
automorphism of $\T^\rho$.

For the torus, the construction in Lemma \ref{embedd} reduces to the
following. Let $B_j, j \in J$, be matrices in ${\Cal M}^*(\rho, \Z)$
and $q_j=|\det(B_j)|$. Suppose for simplicity $J$ finite. Then
$\tilde G$ is the compact group dual of the discrete group $\tilde H
:=\{\prod_{j} q_j^{\ell_j} \, \k, \k \in \Z^\rho, \ell_j \in \Z^{-}
\}$, $\Z^\rho$ is a subgroup of $\tilde H$ and $\T^\rho$ is a factor
of $\tilde G$.

It is well known that $A \in {\Cal M}^*(\rho, \Z)$ acts ergodically on
$(\T^\rho, \mu)$ if and only if $A$ has no eigenvalue
root of unity. A $\Z^d$-action $(A^\el, \el \in \Z^d)$ on
$(\T^\rho, \mu)$ is totally ergodic
if and only if it is free on $\Z^\rho \ \stm0$, or equivalently if $A^\el$ has
no eigenvalue root of unity if $\el \not = \0$.
\begin{lem} \label{irreduc1} Let $M \in {\Cal M}^*(\rho, \Z)$ be
a matrix with irreducible (over $\Q$) characteristic polynomial $P$.
If $\{B_j, j \in J\}$ are $d$ matrices in ${\Cal M}^*(\rho, \Z)$
commuting with $M$, they generate a commutative semigroup $\Cal S$
of endomorphisms on $\T^\rho$ which is totally ergodic if and only
if $B^\el \not = B^{\el'}$, for $\el \not = \el'$ (where $B^\el :=
B_1^{\ell_1}.... B_d^{\ell_d})$. The $\Z^d$-action extending $\Cal
S$ is the product of a totally ergodic $\Z^{d'}$-action, with $d\,' \leq d$
by an action of finite order.
\end{lem}
\proof Since $P$ is irreducible, the eigenvalues of $M$ are
distinct. It follows that the matrices $B_j$ are simultaneously
diagonalizable on $\C$, hence are pairwise commuting. Now suppose
that there are $\el \in \Z^d \stm0$ and $v \in \Z^\rho \setminus
\{\0\}$ such that $B^\el v = v$. Let $\Cal E$ be the subspace of
$\R^\rho$ generated by $v$ and its images by $M$. The restriction of
$B^\el$ to $\Cal E$ is the identity. $\Cal E$ is $M$-invariant, the
characteristic polynomial of the restriction of $M$ to $\Cal E$ has
rational coefficients and factorizes $P$. By the assumption of
irreducibility over $\Q$, this implies $\Cal E = \R^\rho$. Therefore
$B^\el$ is the identity.

Let $\Cal K$ be the kernel of the homomorphism $h: \el \to B^\el$
and $\tilde h$ the quotient of $h$ in $\Z^d / \Cal K$.
The finitely generated group $\Z^d / \Cal K$ is isomorphic to $\Cal U \oplus
\Cal T$, with $\Cal U$ isomorphic to $\Z^{d'}$ for $d' \in [0, d]$,
the restriction of $\tilde h$ to $\Cal U$ a totally ergodic $\Z^{d'}$-action
and $\Cal T$ a finite group.\eop

\vskip 3mm {\bf Examples of $\Z^d$-actions by automorphisms}

In general proving total ergodicity and computing explicit
independent generators is difficult. This may be easier with
endomorphisms. For example, let $(q_j, j \in J)$ be integers $> 1$
and let $x \to q_j x \text{ mod } 1$ be the corresponding
endomorphisms acting on $\T^1$. They generate a semigroup embedded
in a group acting giving a totally ergodic $\Z^d$-action on an
extension of $\T^1$, where $d \in [1, +\infty]$ is the dimension of
the vector space over $\Q$ generated by $\Log \, q_j, j \in J$. The
construction extends to $\rho \times \rho$ matrices $B_j$, $j \in
J$, such that $|\det(B_j)| > 1$ by replacing $q_j$ by $|\det(B_j)|$,
under the condition of Lemma \ref{irreduc1}. Without irreducibility
condition, the action of commuting  $\rho \times \rho$ matrices
$B_j$, when the numbers $\Log \,|\det(B_j)|$ are linearly
independent over $\Q$, extends to a totally ergodic $\Z^d$-action
with $d = \Card(J)$ on an extension of $\T^\rho$.

On the contrary, for automorphisms of $G= \T^\rho$, it can be
difficult to compute independent generators of the generated group
for $\rho > 3$ or 4. We would like to discuss this point and recall
some facts (see in particular \cite{KKS02} and \cite{DaKa10}).

The construction of $\Z^d$-action by automorphisms on $\T^\rho$ is
related to the group of units in number fields (cf. \cite{KKS02}).
To simplify let us consider a matrix $M \in GL(\rho, \Z)$ with an
irreducible (over $\Q$) characteristic polynomial $P$. The elements
of the centralizer of $M$ in $GL(\rho, \R)$ are simultaneously
diagonalizable. The centralizer of $M$ in $\Cal M(\rho, \Q)$ can be
identified with the ring of polynomials in $M$ with rational
coefficients modulo the principal ideal generated by the polynomial
$P$ and hence with the field $\Q(\lambda)$, where $\lambda$ is an
eigenvalue of $M$, by the map $p(A) \to p(\lambda)$ with $p \in
\Q[X]$.

By Dirichlet's theorem, if $P$ has $d_1$ real roots and $d_2$ pairs
of complex conjugate roots, there are $d_1 + d_2-1$ fundamental
units in the group of units in the ring of integers in the field
$K(P)$ associated to $P$. The centralizer $\Cal{C}(M)$ of $M$ in
$GL(\rho, \Z)$ provides a totally ergodic $\Z^{d_1+d_2 -1}$-action
by automorphisms on $\T^\rho$ (up to a product by a finite cyclic
group consisting of roots of unity). The computation of the number
of real roots of $P$ gives the dimension $d = d_1+d_2 -1$ of the
$\Z^d$ free action on $\T^\rho$ generated by $\Cal{C}(M)$.

Nevertheless it can be difficult to compute elements with
determinant $\pm 1$ in $\Cal{C}(M)$. The explicit computation of
fundamental units (hence of independent generators) relies on an
algorithm (see H. Cohen's book \cite{Co93}) which is in practice
limited to low dimensions.

\vskip 3mm \goodbreak {\bf Examples for $\T^3$}

Let $P(x) = -x^3 +qx +n$ be a polynomial with coefficients in $\Z$,
irreducible over $\Q$.

Let $M = \left(
\begin{matrix} 0 & 1 & 0 \cr 0 & 0 & 1 \cr n & q & 0 \cr
\end{matrix} \right)$ be its companion matrix. Let $\lambda$ be a root
of $P$. If the field $K(P)$ is listed in a table giving the
characteristics of the first cubic real fields (see \cite{Co93},
\cite{SteRud76}), we find a pair of fundamental units for the group
of units in the ring of integers in $K(P)$ of the form
$P_1(\lambda)$, $P_2(\lambda)$, with $P_1, P_2 \in \Z[X]$. The
matrices $A_1= P_1(M)$ and $A_2= P_2(M)$ provide elements of
$\Cal{C}(M)$ giving a totally ergodic $\Z^2$-action on $\T^3$ by
automorphisms.

\vskip 3mm \goodbreak {\it 1) Explicit examples (from the table in
\cite{SteRud76})}

a) Let us consider the polynomial $P(x) = - x^3 + 12x + 10$ and its
companion matrix $M = \left( \begin{matrix} 0 & 1 & 0 \cr 0 & 0 & 1
\cr 10 & 12 & 0 \cr
\end{matrix} \right)$. Let $\lambda$ be a root of $P$. The table gives
a pair of fundamental units: $P_1(\lambda) = \lambda^2 -3\lambda -3,
\ P_2(\lambda) = -\lambda^2 +\lambda +11$. Let $A_1, A_2$ be the
matrices

$A_1 = P_1(M)= \left(
\begin{matrix} -3 & -3 & 1 \cr 10 & 9 & -3 \cr -30 & -26 & 9 \cr
\end{matrix} \right), \ \ A_2 = P_2(M) = \left(
\begin{matrix} 11 & 1 & -1 \cr -10 & -1 & 1 \cr 10 & 2 & -1 \cr
\end{matrix}
\right).$ \vskip 3mm b) Consider now the polynomial $P(x) = - x^3 +
9x + 2$ and its companion matrix $M' = \left(
\begin{matrix} 0 & 1 & 0 \cr 0 & 0 & 1 \cr 2 & 9 & 0 \cr
\end{matrix}
\right).$ Let $\lambda$ be a root of $P$. The table gives a pair of
fundamental units for the algebraic group associated to $P$:
$P_1(\lambda) = 85\lambda^2 - 245\lambda -59, \ P_2(\lambda) = - 18
\lambda^2 + 4 \lambda + 161$.

Take $A_1 = P_1(M')= \left(
\begin{matrix} -59 & -245 & 85 \cr 170 & 706 & -245 \cr -490 & -2035 & 706 \cr
\end{matrix} \right), \ \ A_2 = P_2(M') = \left(
\begin{matrix} 161 & 4 & -18 \cr -36 & -1 & 4 \cr 8 & 0 & -1 \cr
\end{matrix}
\right)$.

In both cases a) and b), the matrices $A_1$ and $A_2$ are in $GL(3,
\Z)$ and generate a totally ergodic actions of $\Z^2$ by
automorphisms on $\T^3$.

{\it 2) A simple example on $\T^4$}

If $P(x) = x^4 + a x^3 + b x^2 + a x + 1$, the polynomial $P$ has
two real roots: $\lambda_0, \lambda_0^{-1}$ and two complex
conjugate roots of modulus 1: $\lambda_1, \overline \lambda_1$.

Let $\sigma_j= \lambda_j + \overline \lambda_j$, $j=0,1$. They are
roots of $Z^2 -a Z +b-2 =0$.

If the conditions: $a^2 -4b+8 > 0, a > 4, b > 2, 2a > b+2$ are
satisfied (i.e., $2 < b < 2a -2, \ a > 4$, since $2a-2 \leq
\frac14a^2 +2$), then $\lambda_0, \lambda_0^{-1}$ are solutions of
$\lambda^2 - \sigma_0 \lambda +1 = 0$, and $\lambda_1, \overline
\lambda_1$ are solutions of $\lambda^2 - \sigma_1 \lambda +1 = 0$,
where
$$\sigma_0 = -\frac12 a - \frac12 \sqrt{a^2 -4b +8},
\ \sigma_1 = -\frac12 a + \frac12 \sqrt{a^2 -4b +8}.$$ The
polynomial $P$ is not factorizable over $\Q$. Indeed, suppose that
$P = P_1 P_2$ with $P_1$, $P_2$ with rational coefficients and
degree $\geq 1$. Since the roots of $P$ are irrational, the degrees
of $P_1$ and $P_2$ are 2. Necessarily their roots are, say,
$\lambda_1, \overline \lambda_1$ for $P_1$, $\lambda_0,
\lambda_0^{-1}$ for $P_2$. The sum $\lambda_1 + \overline
\lambda_1$, root of $Z^2 -a Z +b-2 =0$, is not rational and the
coefficients of $P_1$ are not rational. Let us take $A := \left(
\begin{matrix} 0 & 1 & 0 & 0 \cr 0 & 0 & 1 & 0 \cr 0 & 0 & 0 & 1 \cr
-1 & -a & -b & -a \cr
\end{matrix}
\right), B = A+I.$ From the irreducibility over $\Q$, it follows
that, if there is a non zero fixed integral vector for $A^k
B^{\ell}$, where $k, \ell$ are in $\Z$, then we have $A^k B^{\ell} =
Id$. This implies: $\lambda_1^k \, (\lambda_1 -1)^\ell = 1$, hence,
since we have $|\lambda_1| = 1$, it follows $|\lambda_1 + 1| = 1$,
i.e. $\lambda$ is also solution of $z^2 - z + 1 = 0$, which is not
true.

\vskip 2mm An example is $P(x) = x^4 + 5x^3 +7 x^2 + 5x +1$. If $A$
is the companion matrix, then the matrices $A$ and $B= A+I$ generate
a $\Z^2$-totally ergodic action on $\T^4$.

\vskip 3mm {\it 3) Construction by blocks}

Let $M_1, M_2$ be two ergodic matrices respectively of dimension
$d_1$ and $d_2$. Let $(p_i, q_i)$, $i=1,2$, be two pairs of integers
such that $p_1q_2 - p_2 q_1 \not = 0$. On the torus $\T^{d_1+d_2}$
we obtain a $\Z^2$-totally ergodic action by taking $A_1, A_2$ of
the following form: $A_1 = \left( \begin{matrix} M_1^{p_1} & 0 \\ 0
& M_2^{q_1} \cr
\end{matrix}
\right), \ A_2 = \left( \begin{matrix} M_1^{p_2} & 0 \\ 0 &
M_2^{q_2} \cr
\end{matrix}
\right).$ Indeed, if there exists $v = \left( \begin{matrix} v_1 \\
v_2\cr
\end{matrix} \right) \in \Z^{d_1 + d_2} \setminus \{0\}$ invariant by $A_1^n
A_2^\ell$, then $M_1^{np_1 + \ell p_2} v_1 = v_1, \ M_2^{nq_1 + \ell
q_2} v_2 = v_2$, which implies $np_1 + \ell p_2 = 0$, $nq_1 + \ell
q_2 = 0$; hence $n = \ell = 0$.

This method gives explicit free $\Z^2$-actions on $\T^4$. The same
method gives free $\Z^3$-actions on $\T^5$. As discuss
above, it is more difficult to explicit examples of full dimension,
i.e., with 3 independent generators on $\T^4$, or with 4 independent
generators on $\T^5$.
\subsection{\bf Spectral densities and Fourier series for tori}

\

The continuity of $\varphi_f$ for a general compact abelian group
follows from the absolute convergence of the Fourier series of $f$,
hence, for $G = \T^\rho$ from the condition:
\begin{eqnarray}
|c_f(\k)| = O(\|k\|^{-\beta}), \text{ with } \beta > \rho.
\label{regFour2b2}
\end{eqnarray}
Condition (\ref{regFour2b2}) implies a rate of approximation of $f$
by the partial sums of its Fourier series, but for the torus a
weaker regularity condition on $f$ can be used. First, let us recall
a result on the approximation of functions by trigonometric
polynomials.

For $f\in L^2(\T^\rho)$, the Fourier partial sums of $f$ over squares
of sides $N$ are denoted by $s_{N}(f)$. The {\it integral modulus of
continuity} of $f$ is defined as
$$\omega_2(\delta, f)=\sup_{|\tau_1|\le
\delta,\ldots,|\tau_\rho|\le \delta} \|f(.+\tau_1, \cdots, .+\tau_\rho)-
f\|_{L^2(\T^\rho)}.$$
\begin{lem} The following condition on the modulus of continuity
\begin{eqnarray}
&&\exists \ \alpha > d \text{ and } \, C(f)<+\infty \text{
such that } \ \omega_{2}(\delta, f) \leq C(f) \, (\ln {1 \over
\delta})^{-\alpha}, \forall \delta > 0, \label{5.11}
\end{eqnarray}
implies, for a constant $R(f)$:
\begin{eqnarray}
&&\|f - s_{N}(f)\|_2 \leq R(f) \, (\ln N)^{-\alpha}, \text{ with }
\alpha > d. \label{ineq5.23b2}
\end{eqnarray}
\end{lem}
\proof Let $K_{N_1,\ldots,N_\rho}$ be the $\rho$-dimensional Fej\'er
kernel, for $N_1,..., N_\rho \geq 1$, and let $J_{N_1, \ldots,N_\rho}
(t_1,\cdots,t_\rho)=K_{N_1,\ldots, N_\rho}^2(t_1,\cdots, t_\rho)
/\|K_{N_1,\ldots,N_\rho}\|_{L^2(\T^\rho)}^2$ be the $\rho$-dimensional
Jackson's kernel. Using the moment inequalities
$\int_0^{\frac12}t^k\, J_{N}(t) \, dt = O(N^{-k}), \ \forall N \geq
1, \, k=0,1,2$, satisfied by the 1-dimensional Jackson's kernel,
we obtain that there exists a positive constant $C_\rho$ such that,
for every $f \in L^2(\T^\rho)$, $\|J_{N,\ldots,N}*
f-f\|_2\le C_\rho\,\omega_2(\frac1{N}, f), \forall N \ge1$.

It follows: $\|f-s_{N}(f)\|_2 \le C_{\rho}\, \omega_2(\frac1{N}, f), \
\forall f\in L^2(\T^\rho), \ \forall N \ge 1$, since
$\|f-s_{N}(f)\|_2\leq \|f-P\|_2$ for every trigonometric polynomial
$P$ in $\rho$ variables of degree at most $N \times \cdots\times N$.
Hence (\ref{ineq5.23b2}) follows from (\ref{5.11}). \eop

\vskip 3mm
The required regularity in the next theorem is weaker than
in Proposition \ref{condAbsConv}. The proof is like that of the
analogous result in \cite{Leo60c}. It uses the lemma below due to D.
Damjanovi\'c and A. Katok \cite{DaKa10}, extended by M. Levin
\cite{Lev13} to endomorphisms.
\begin{lem} \cite{DaKa10}  \label{DaKa}
If $(A^\n, \n \in \Z^d)$ is a totally ergodic $\Z^d$-action on
$\T^\rho$ by automorphisms, there are $\tau > 0$ and $C > 0$, such
that for all $(\n, \k) \in \Z^d \times (\Z^\rho \stm0)$ for which
$A^\n \k\ \in \Z^\rho$.
\begin{eqnarray}
\|A^\n \k\| \geq C e^{\tau \|\n\|} \|\k\|^{-\rho}. \label{DaKaMin}
\end{eqnarray}
\end{lem}
\begin{thm} \label{ThmRatedeCor} Let $\el \to A^\el$ be a totally ergodic
$d$-dimensional action by commuting endomorphisms on $\T^\rho$. Let
$f$ be in $L_0^2(\T^\rho)$ satisfying the regularity condition
(\ref{5.11}) or more generally (\ref{ineq5.23b2}). Let $f_1(x) :=
\sum_{\n \in \Cal{E}_1} c_\n (f) e^{2\pi i\langle \n , x \rangle}$,
where $\Cal{E}_1$ is a subset of $\Z^\rho$. Then there are finite
constants $B(f), C(f)$ depending only on $R(f)$ such that
\begin{eqnarray}
|\langle A^{\el} f_1, f_1\rangle| \leq B(f) \|f_1\|_2
\|\el\|^{-\alpha}, \ \forall \el \not = \0, \label{ineq5.24}
\end{eqnarray}
the spectral density is continuous, $\sum_{\el \in \Z^d} |\langle
A^{\el}f, f\rangle| < \infty$ and $\|\varphi_{f - f_1}\|_\infty \leq
C(f) \|f-f_1\|_2$.
\end{thm}
\proof (Recall the convention $c_{A^\el \k}(f) = 0$ if $A^\el\k \not
\in \Z^\rho$.) It suffices to prove the result for $f$ since, by
setting $c_f(n)=0$ outside $\Cal{E}_1$, we obtain (\ref{ineq5.24})
with the same constant $B(f)$ as shown by the proof. Let $\lambda,
b, h$ be such that $1 < \lambda < e^{\tau}, \ 1 < b <
\lambda^{{1\over \rho}}, \ h:=\lambda b^{-\rho} > 1$ (where $\tau$
is given by Lemma \ref{DaKa}). We have for $\el \in \Z^d$:
\begin{eqnarray}
&&\langle A^{\el} f, f\rangle = \sum_{\k \, \in \Z^\rho} c_\k(f) \,
\overline c_{A^\el\k}(f) = \sum_{\|\k\| < b^{\|\el\|}} \ + \
\sum_{\|\k\| \geq b^{\|\el\|}} = (A) + (B). \label{ineq5.25}
\end{eqnarray}
 From Inequality (\ref{DaKaMin}) of Lemma \ref{DaKa}, we deduce
that, if $\|\k\| < b^{\|\el\|}$, then $\|A^\el\k\| \geq D
\lambda^{\|\el\|} \, \|\k\|^{-\rho} \geq D \lambda^{\|\el\|} \,
b^{-\rho \|\el\|} = D h^{\|\el\|}, \ \el \not = \0$. It follows, for
the sum (1):
\begin{eqnarray*}
|(A)| \leq (\sum_{\|\k\| < b^{\|\el\|}} |c_\k(f)|^2)^\frac12 \
(\sum_{\|\k\| < b^{\|\el\|}} |c_{A^\el\k}(f)|^2)^\frac12 \leq
\|f\|_2 \,\sum_{\|\m\|
> Dh^{\|\el\|}} |c_{\m}(f)|^2.
\end{eqnarray*}

By Parseval inequality and (\ref{ineq5.23b2}), there is a finite
constant $B_1(f)$ such that, for $\el \not = 0$:
\begin{eqnarray}
&& (\sum_{\|\m\| > Dh^{\|\el\|}} |c_{\m}(f)|^2)^\frac12 \leq \|f-
s_{[Dh^{\|\el\|}], ..., [Dh^{\|\el|}]}(f)\|_2 \leq {R(f) \over (\ln
[Dh^{\|\el\|}])^\alpha} \leq B_1(f) \|\el\|^{-\alpha}.
\label{ineq5.28}
\end{eqnarray}
 From the previous inequalities, it follows: $|(A)| \leq B_1(f)
\|f\|_2 {\|\el\|}^{-\alpha}, \forall \, {\|\el\|} \not = 0$.
Analogously, we obtain $|\sum_{\|\k\| \geq b^{\|\el\|}} c_\k(f) \,
\overline c_{A^\el\k}(f)| \leq B_2(f) \|f\|_2 {\|\el\|}^{-\alpha},
\, \el \not = 0$ for $(B)$ in (\ref{ineq5.25}).

Taking $B(f)= B_1(f) + B_2(f)$, (\ref{ineq5.24}) follows from
(\ref{ineq5.25}) and implies the last statements. \eop

\vskip 5mm \section{\bf CLT for summation sequences of endomorphisms
on $G$} \label{CLTGroup}

\vskip 2mm \subsection{\bf Criterium for the CLT on a compact
abelian connected group $G$}

\

For $\Z^d$-dynamical systems satisfying the $K$-property, a
martingale-type property can be used to obtain a CLT.
(For martingale methods applied to $d$-dimensional random fields,
see for example \cite{Go09}, \cite{VoWa14}.) For $\Z^d$-action by automorphisms, the $K$-property is
equivalent to mixing of all orders (cf. \cite{SchWar93}) for zero-dimensional compact abelian groups
but does not hold for instance for the model that we consider on tori. In the
absence of $K$-property, we will use for abelian semigroups of
endomorphisms of connected compact abelian groups the method of
mixing of all orders applied by Leonov to a single ergodic
automorphism.

\vskip 2mm {\bf Mixing actions by endomorphisms ($G$ connected)}
\label{mixAuto}

The proof of the CLT given by Leonov in \cite{Leo60b} for a single
ergodic endomorphism $A$ of a compact abelian group $G$ is based on
the computation of the moments of the ergodic sums $S_nf$ when $f$
is a trigonometric polynomial. It uses the fact that $A$ is mixing
of all orders, which follows from the $K$-property for the
$\Z$-action of a single ergodic automorphism (\cite{Roh61}). For
$\Z^d$-actions by automorphisms on compact abelian groups, mixing of
all orders is not always true (cf. \cite{Le78}, \cite{Sc95}), but it
is satisfied for actions on connected compact abelian groups
(Theorem \ref{r-mixing} below) and the method of moments can be
used.

In 1992, W. Philip in \cite{Ph94} and K. Schmidt and T. Ward in
\cite{SchWar93} applied results on the number of solutions of
$S$-units equations (see (\cite{Schl90, EvScSch02}) to endomorphisms
or automorphisms of compact abelian groups.
\begin{thm} \label{r-mixing} (\cite [Corollary 3.3] {SchWar93})
Every 2-mixing $\Z^d$-action by automorphisms on a compact connected
abelian group G is mixing of all orders.\end{thm} With the notations
of Lemma \ref{embedd}, if ${\Cal S}$ is a totally ergodic semigroup
of endomorphisms on a compact connected abelian group $G$, then its
extension $\tilde {\Cal S}$ to a group of automorphisms of $\tilde
G$ is mixing of all orders by Theorem \ref{r-mixing}.

{\it From now on, we consider a totally ergodic $\Z^d$-action $\el
\to A^\el$ by commuting automorphisms on $G$ (or on the extension
$\tilde G$, but we will not write $\tilde \ $) which is mixing of
all orders (an assumption satisfied when $G$ is connected by Theorem
\ref{r-mixing}).}

Let $(R_n)_{n \geq 1}$ be a summation sequence on $\Z^d$. We use the
notations and the results of the appendix (\Sect \ref{cumulAppend})
on cumulants. For $f \in L^2(G)$, we put $\sigma_n(f) := \|\sum_\el
R_n(\el) \, A^\el f \|_2$ and assume $\sigma_n^2(f) \not = 0$, for
$n$ big.
\begin{lem} \label{tclPol0} For a trigonometric polynomial $f$ with zero mean,
the condition
\begin{eqnarray}
&& \sum_\el \prod_{k= 1}^r R_n(\el + {\j}_k) = o(\sigma_n^{r}(f)),
\forall \{\j_1, ..., \j_r\} \in \Z^d, \forall r \geq 3, \text{
implies} \label{cumulCN}
\end{eqnarray}
\begin{eqnarray}
\sigma_n(f)^{-1} \, \sum_\el R_n(\el) \, A^\el f \
\overset{distrib} {\underset{n \to \infty} \longrightarrow} \ \Cal
N(0,1). \label{CLTPol0}
\end{eqnarray}
\end{lem}
\proof Let $(\chi_\k, \k \in \Lambda)$ be a finite set of characters
on $G$, $\chi_{0}$ the trivial character. If $f = \sum_{\k \in
\Lambda} c_\k(f) \, \chi_\k$, the moments of the process $(f(A^\n
.))_{\n \in \Z^d}$ are
\begin{eqnarray*}
&&m_f(\n_1, ..., \n_r) = \int f(A^{\n_1} x) ... f(A^{\n_r} x) \ dx =
\sum_{\k_1, ..., \k_r \in \Lambda} c_{k_1} ... c_{k_r} 1_{A^{\n_1}
\chi_{\k_1} ... A^{\n_r} \chi_{\k_r} = \, \chi_{0}}.
\end{eqnarray*}
 For $r$ fixed, the function $(\k_1, ..., \k_r) \to m_f(\k_1, ...,
\k_r)$ takes a finite number of values, since $m_f$ is a sum with
coefficients 0 or 1 of the products $c_{k_1} ... c_{k_r}$ with $k_j$
in a finite set. The cumulants of given order according to Identity
(\ref{cumFormu1}) take also a finite number of values.

Therefore, since mixing of all orders implies $\underset{\max_{i,j}
\|\el _i - \el _j\| \to \infty} \lim \, s_f(\el _1,..., \el _r) = 0$
(cf. Notation (\ref{notasf})) by Lemma \ref{skTozeroLem}, there is
$M_r$ such that $s_f(\el _1, ..., \el _r) = 0$ if $\max_{i,j} \|\el
_i - \el_j\| > M_r$.

We apply Theorem \ref{Leonv} (cf. appendix). Let us check
(\ref{smallCumul}), i.e.,
\begin{eqnarray*}
\sum_{(\el_1, ..., \el_r)
\, \in (\Z^d)^r} \, c(X_{\el_1}, ... , X_{\el_r}) \, R_n(\el_1) ...
R_n(\el_r) = o(\|Y^n\|_2^r), \forall r \geq 3.
\end{eqnarray*}
Using (\ref{cumLin0}), we obtain
\begin{eqnarray*}
&&|\sum_{\el_1, ..., \el_r} s_f(\el _1, ..., \el _r) \, R_n(\el
_1)...R_n(\el _r)| = |\sum_{\max_{i,j} \|\el _i - \el_j\| \leq M_r}
s_f(\el_1, \el _2, ..., \el _r) \, R_n(\el _1)...R_n(\el _r)| \\
&&\ \ \leq \sum_{\el} \sum_{\|\j _2\|, ..., \|\j_r\| \leq M_r, \,
j_1 =\0} |s_f(\el, \el+ \j_2, ..., \el+\j_r)| \, \prod_{k=1}^r
R_n(\el+\j_k)\\ && \ \ = \sum_\el \sum_{\|\j _2\|, ..., \|\j_r\|
\leq M_r, \, j_1 =\0} |s_f(\j_1,\j_2, ..., \j_r)| \ \prod_{k=1}^r
R_n(\el+j_k).
\end{eqnarray*}
The right hand side is less than $C \sum_\el \sum_{\|\j _2\|, ...,
\|\j_r\| \leq M_r, \, j_1 =\0} \ \prod_{k=1}^r R_n(\el+j_k)$.
Therefore (\ref{cumulCN}) implies (\ref{smallCumul}). \eop
\begin{thm} \label{tclRegKer} Let $(R_n)_{n \geq 1}$ be a summation
sequence on $\Z^d$ which is $\zeta$-regular (cf. Definition
\ref{regularProcess}). Let $f$ be a function in $AC_0(G)$ with
spectral density $\varphi_f$. The condition
\begin{eqnarray} && \bigl(\sup_\el R_n(\el)\bigr)^{r-1} \, \sum_\el R_n(\el)
= o((\sum_{\el \in \Z^d} \, R_n(\el)^2)^{r/2}), \, \text{for every }
r \geq 3, \label{cumulCN3}
\end{eqnarray}
implies (with the convention that the limiting distribution
is $\delta_0$ if $\sigma^2(f) = 0$)
\begin{eqnarray}
(\sum_{\el \in \Z^d} \, R_n(\el)^2)^{-\frac12} \, \sum_{\el \in
\Z^d} R_n(\el) f(A^\el .) \overset{distr} {\underset{n \to \infty}
\longrightarrow } \Cal N(0, \zeta(\varphi_f)). \label{cvgce1}
\end{eqnarray}
\end{thm}
\proof We use (\ref{densspectAC}) and the $\zeta$-regularity of
$(R_n)$: for $g$ in $AC_0(G)$,
$$(\sum_{\el \in \Z^d} \, R_n(\el)^2)^{-1} \, \|\sum_{\el \in \Z^d} R_n(\el)
\, A^\el g \|_2^2 = \int_{\T^d} \, \tilde R_n \, \varphi_g \, dt
\underset{n \to \infty} \to \zeta(\varphi_g).$$ Let $(\Cal E_s)_{s
\geq 1}$ be an increasing sequence of finite sets in $\hat G$ with
union $\hat G \setminus \{0\}$ and let $f_s(x) := \sum_{\chi \in
\Cal E_s} \, c_f(\chi) \, \chi$ be the trigonometric polynomial
obtained by restriction of the Fourier series of $f$ to $\Cal E_s$.
Let us consider the processes defined respectively by
\begin{eqnarray*}
U_n^{s} := (\sum_{\el \in \Z^d} \, R_n(\el)^2)^{-\frac12} \,
\sum_{\el \in \Z^d} R_n(\el) f_s(A^\el .), \ U_n := (\sum_{\el \in
\Z^d} \, R_n(\el)^2)^{-\frac12} \, \sum_{\el \in \Z^d} R_n(\el)
f(A^\el .).
\end{eqnarray*}
We can suppose $\zeta(\varphi_f) > 0$, since otherwise the limiting
distribution is $\delta_0$. By Proposition \ref{condAbsConv} we have
$\zeta(\varphi_{f - f_s}) \leq \|f-f_s\|_c$. It follows
$\zeta(\varphi_{f_s}) \not = 0$ for $s$ big enough. We can apply
Lemma \ref{tclPol0} to the trigonometric polynomials $f_{s}$, since
$\sigma_n^2(f_s) \sim (\sum_{\el \in \Z^d} \, R_n(\el)^2) \,
\zeta(\varphi_{f_s})$ with $\zeta(\varphi_{f_s}) > 0$ and since
Condition (\ref{cumulCN3}) implies Condition (\ref{cumulCN}) in
Lemma \ref{tclPol0}.

It follows: ${U_n^{s} \overset{distr} {\underset{n \to \infty}
\longrightarrow} \Cal N(0,\zeta(\varphi_{f_s}))}$ for every $s$.
Moreover, since
\begin{eqnarray*}
\lim_n \int |U_n^s - U_n|_2^2 \ d\mu &=& \lim_n \int_{\T^d} \,
\tilde R_n \,\varphi_{f-f_s} \, d \t = \zeta(\varphi_{f - f_s}) \leq
\|f-f_s\|_c,
\end{eqnarray*}
we have $\limsup_n \mu[|U_n^s - U_n| > \varepsilon] \leq
\varepsilon^{-2} \limsup_n \int |U_n^s - U_n|_2^2 \ d\mu \underset{s
\to \infty} \to 0$ for every $\varepsilon > 0$.

Therefore the condition $\lim_s \limsup_n \mu[|U_n^s - U_n| >
\varepsilon] = 0$, $\forall \varepsilon > 0$, is satisfied and the
conclusion $U_n \overset{distr} {\underset{n \to \infty}
\longrightarrow } \Cal N(0, \zeta(\varphi_{f}))$ follows from
Theorem 3.2 in \cite{Bill99}. \eop

With the notations of \Sect \ref{summSeqSect}, Theorem
\ref{tclRegKer} implies for sequential summations:
\begin{cor} \label{tclSeqSum}
If $(\x_n)$ is a sequence in $\Z^d$ such that  $\z_n =
\sum_{k=0}^{n-1} \x_k$ is $\zeta$-regular, then the convergence
$v_n^{-\frac12} \, \sum_{k=0}^{n-1} f(A^{z_k} .) \overset{distr}
{\underset{n \to \infty} \longrightarrow } \Cal N(0,
\zeta(\varphi_f))$ follows from the condition
\begin{eqnarray}
&&n \, \bigl(\sup_\el \sum_{k=0}^{n-1} 1_{z_k = \el}\bigr)^{r-1} =
o(v_n^{r/2}), \, \forall  r \geq 3. \label{cumulCN6}
\end{eqnarray}
\end{cor}

Before considering random walks, let us apply Theorem
\ref{tclRegKer} to summation over sets:
\begin{cor} \label{tclFoln} Let $(D_n)_{n \geq 1}$ be a F\o{}lner
sequence of sets in $\N^d$ and let $f$ be in $AC_0(G)$. We have
$\sigma^2(f) = \lim_n \|\sum_{\el \in D_n} \, A^\el f\|_2^2/|D_n| =
\varphi_f(0)$ and
$$|D_n|^{-\frac12} \, \sum_{\el \in D_n} \, A^\el f(.)
\overset{distr} {\underset{n \to \infty} \longrightarrow } \Cal
N(0,\sigma^2(f)).$$
\end{cor} \proof The sequence $R_n(\el) = 1_{D_n}(\el)$ is $\zeta$-regular,
with $\zeta = \delta_0$. Suppose that $\varphi_f(0) \not = 0$. We
have $\sigma_n^2(f) \sim |D_n| \, \varphi_f(0)$ and $R_n(\el +
{\j}_k) = 0$ or 1. Therefore Condition (\ref{cumulCN3}) holds and
the result follows from Theorem \ref{tclRegKer}. \eop

\begin{rems} 1) The previous result is valid for the rotated sums: for $f$ in
$AC_0(G)$, for every $\theta$,
\begin{eqnarray}
\sigma_\theta^2(f) = \varphi_f(\theta), \ |D_n|^{-\frac12} \,
\sum_{\el \in D_n} e^{2 \pi i \langle\el, \theta\rangle} f(A^\el .)
\overset{distr} {\underset{n \to \infty} \longrightarrow } \Cal
N(0,\sigma_\theta^2(f)). \label{CLTRot1}
\end{eqnarray}
2) When $G = \T^d$, in view of Theorem \ref{ThmRatedeCor}, the
conclusions of Theorem \ref{tclRegKer} and Corollary \ref{tclFoln}
are valid under the weaker regularity assumption (\ref{ineq5.23b2}),
in particular under a logarithmic H\"olderian regularity (condition
(\ref{5.11})).

3) As mentioned in the introduction, the result of Corollary
\ref{tclFoln} for the sums over $d$-dimensional rectangles and
regular functions was obtained by M. Levin (\cite{Lev13}).
\end{rems}
{\it A CLT for the rotated sums for a.e. $\theta$ without regularity
assumption}

When $(D_n)$ is a sequence of $d$-dimensional cubes in $\Z^d$, the
following CLT for the rotated sums holds for a.e. $\theta$ without
regularity assumptions on $f$.
\begin{thm} \label{rotatedCLT} Let $(D_n)_{n \geq 1}$ be a sequence of cubes
in $\Z^d$. For $f$ in $L^2(G)$, we have for a.e. $\theta \in \T^d$:
$\sigma_\theta^2(f) = \varphi_f(\theta)$ and $|D_n|^{-\frac12} \,
\sum_{\el \in D_n} e^{2 \pi i \langle\el, \theta\rangle} \, A^\el
f(.) \overset{distr} {\underset{n \to \infty} \longrightarrow } \Cal
N(0,\sigma_\theta^2(f))$.
\end{thm}
\proof As in \cite{CohCo12}, we use the relation $\lim_n
|D_n|^{-\frac12} \, \| \sum_{\el \in D_n} \, e^{2\pi i \langle \el,
\theta \rangle} \, T^\el(f -M_\theta f)\|_2 = 0$, which is
satisfied, for any $f \in L^2(G)$, for $\theta$ in a set of full
measure (depending on $f$).\eop

\section{\bf Application to r.w. of commuting endomorphisms on $G$}
\label{applRW}

Now we apply the previous sections to random walks of commuting
endomorphisms or automorphisms on a compact abelian group $G$.

Let us consider a family $(B_j, j \in J)$ of commuting endomorphisms
of $G$ (extended if necessary to automorphisms of $\tilde G$) and a
probability vector $\nu = (p_j, j \in J)$ such that $p_j > 0,
\forall j$. These data define a random walk $U_n := Y_0...Y_{n-1}$
where $(Y_k)_{k \geq 0}$ are i.i.d. r.v. with common distribution
$\PP(Y_k = B_j) = p_j$.

If the $B_j$'s are given a priori, we can try to express $U_n$ via a
r.w. on $\Z^d$ for some $d$. This means that we have to find
algebraically independent generators $A_1, ..., A_d$ for some $d$,
in order to express the $B_j$'s as $B_j = A^{\el_j}$, for $\el_j \in
\Z^d$. Another approach is to start from a totally ergodic $\Z^d$-action
$\A$ on $G$ (or $\tilde G$) and from a r.w. $W$ on $\Z^d$ and to
transfer $W$ to a r.w. on the group of automorphisms of $\tilde G$,
thus getting the $B_j$'s a posteriori.

For example, when $G$ is a
torus, we can use the method described in Subsection \ref{torEx},
which gives an explicit construction of algebraically independent
generators of a totally ergodic $\Z^d$-action by invertible matrices
on a torus.

With $(\Omega, \PP) = ((\Z^d)^\Z, \nu^{\otimes \Z})$ and $(X_n)$ the
sequence of coordinates maps of $\Omega$, we obtain on $(\Omega
\times \tilde G, \PP \times \tilde \mu)$ a dynamical system
$(\omega, x) \to (\tau \omega, A^{X_0(\omega)} x)$. The iterates are
$(\omega, x) \to (\tau^n \omega, A^{X_0(\omega)+ ...+
X_{n-1}(\omega)} x) = (\tau^n \omega, A^{Z_n(\omega)} x)$, $n \geq
1$.

The random walk $U_n= Y_0...Y_{n-1}$ can be expressed as
$(A^{Z_n})$, where $(Z_n)$ is the r.w. in $\Z^d$ (which can be and
will be assumed reduced) with distribution $\PP(X_0 = \el_j) = p_j$.
\vskip 3mm \subsection{\bf Random walks and quenched CLT}

\

Recall that the measure $d\gamma$ is defined in Definition
\ref{defdgamma}, $w(\t) = {1 - |\Psi(\t)|^2 \over |1 - \Psi(\t)|^2}$
and $V_n(\omega) = \#\{0 \leq k', k < n: \ Z_k(\omega) =
Z_{k'}(\omega) \}$.

\begin{thm} \label{main} Let $W$ be a reduced centered r.w.
Let $\el \to A^\el$ be a totally ergodic $\Z^d$-action by
automorphisms on $G$. Let $f$ be in $AC_0(G)$ with spectral density
$\varphi_f$.

I) Suppose $W$ with a finite moment of order 2 and centered.  \hfill
\break a) If $d =1$, then, for a.e. $\omega$, $\displaystyle
V_n(\omega) ^{-\frac12} \sum_{k=0}^{n-1} \, A^{Z_k(\omega)} f(.) \
\overset{distr} {\underset{n \to \infty} \longrightarrow } \ \Cal
N(0, \gamma(\varphi_f))$.

b) If $d =2$, then, for a.e. $\omega$, $\displaystyle (C n \Log
n)^{-\frac12} \sum_{k=0}^{n-1} \, A^{Z_k(\omega)} f(.) \
\overset{distr} {\underset{n \to \infty} \longrightarrow } \ \Cal
N(0, \gamma(\varphi_f))$,  with $C = \pi^{-1} a_0(W) \,
\det(\Lambda)^{-\frac12}$.

II) If $W$ is transient with finite moment of order $\eta$ for some
positive $\eta$, then, for a.e. $\omega$, $\displaystyle (C
n)^{-\frac12} \sum_{k=0}^{n-1} \, A^{Z_k(\omega)} f(.) \
\overset{distr} {\underset{n \to \infty} \longrightarrow } \ \Cal
N(0, \zeta(\varphi_f))$. We have for $\zeta$ the following cases:
\hfill \break if $d = 1$, then $C = c_w +K$ and $d\zeta(\t) = (c_w +
K)^{-1}\, (w(\t) \, d\t + d\gamma(t))$, (cf. notations of Theorem
\ref{transrecSpec}) ($K \not = 0$ if $m(W)$ is finite); \hfill
\break if $d \geq 2$, then $C = c_w$ and $d\zeta(\t) = c_w^{-1} \,
w(\t) \, d\t$ (the absolutely continuous part is non trivial if and
only if the random walk is non deterministic).
\end{thm} \proof Theorem \ref{transrecSpec} gives the
$\zeta$-regularity for the r.w. summation $(R_n(\omega, \el))_{n
\geq 1} = (\sum_{k=0}^{n-1} 1_{Z_k(\omega) = \el})_{n \geq 1}$ and
the expression of $\zeta$. This is the first step (variance). It
remains to check Condition (\ref{cumulCN6}) of Corollary
\ref{tclSeqSum}, which reads here:
\begin{eqnarray}
&&n \, \bigl(\sup_\el \sum_{k=0}^{n-1} 1_{Z_k(\omega) =
\el}\bigr)^{r-1} = o(V_n(\omega)^{r/2}), \, \text{ for every } r
\geq 3. \label{cumulCN7}
\end{eqnarray}
a) For the recurrent 1-dimensional case, since $\sum_\el
R_n^2(\omega, \el) \geq C n^\frac32 / \Log\Log \, n$ for a.e.
$\omega$ by (\ref{minVn1}), to have (\ref{cumulCN7}) it suffices
that $n \, (\sup_\el R_n(\omega, \el))^{r-1} = o\bigl(({n^\frac32
\over \Log\Log \, n})^{r/2}\bigr), \forall r \geq 3$, i.e.,
$$\sup_\el R_n(\omega, \el) = o(n^{{3r - 4 \over 4r - 4}} \,
(\Log\Log \, n)^{-{r \over 2r -2}}), \forall r \geq 3.$$ The
condition is satisfied, since the exponent ${3r - 4 \over 4r - 4}$
is bigger than $\frac58$ and $\sup_\el R_n(\el) \leq n^{\frac12 +
\varepsilon}$.

b) For the recurrent 2-dimensional case, for a.e. $\sum_\el
R_n^2(\omega, \el) \sim \E \sum_\el R_n^2(., \el) \sim C n \Log \,
n$.

We need: $n \, \bigl(\sup_\el R_n(\el)\bigr)^{r-1} = o((n \Log \,
n)^{r/2}), \forall r \geq 3$, i.e.,
\begin{eqnarray*}
\sup_\el R_n(\el) = o(n^{{r -2 \over 2r-2}} \ (\Log \, n)^{{r\over
2r -2}}) , \forall r \geq 3.
\end{eqnarray*}
The above condition is satisfied, since the exponent $n^{{r -2 \over
2r-2}}$ increases from $\frac14$ to $\frac12$ when $r$ varies from 3
to $+\infty$ and by Lemma \ref{BolthausenProp} $\sup_\el R_n(\el) =
o(n^\varepsilon), \forall \varepsilon > 0$.

II) For the transient case, we have the same estimation without the
logarithmic factor. Excepted in the deterministic case, the variance
is $>0$ unless $f \equiv 0$ a.e.\eop
\begin{rem}
If $(\Phi_n(\omega, x))_{n \geq 1}$ is a process depending on two
variables $x$ and $\omega$ such that, for a normalization by
$\sigma_n$ independent of $\omega$, $\int e^{it \sigma_n^{-1}
\Phi_n(\omega, x)} \, d\mu(x) \ \to e^{-t^2/2}$, for a.e. $\omega$,
then the CLT holds for $\Phi_n$ w.r.t. $d\mu \times d\PP$ holds,
since by the dominated convergence theorem: $\int_\Omega
\bigl(\int_X e^{it \sigma_n^{-1} \Phi_n(\omega, x)} \, d\mu(x))\,
d\PP(\omega) \to e^{-t^2/2}$.

It follows that, for a r.w. $W$ which is transient or with finite
variance, centering and $d(W) =2$, the annealed version of the CLT
in Theorem \ref{main}.b is satisfied. This result can be viewed as a
``toral" version of Bolthausen theorem in \cite{Bo89}.

For $W$ with finite variance, centering and $d(W) = 1$, the theorem
of Kesten and Spitzer in \cite{KeSp79} for a r.w. in random scenery
gives an (annealed) convergence toward a distribution which is not
the normal law. For the quenched process in their model or in the
toral model of Theorem \ref{main}, convergence toward a normal law
holds, but with a normalization depending on $\omega$. Let us
mention that, in the toral model, the annealed theorem analogous to
the result of \cite{KeSp79} for the recurrent 1-dimensional r.w.
holds (S. Le Borgne, personal communication). Let us mention other
results of quenched type like for instance in \cite{GuPlPoi13}.
\end{rem}

{\bf An example}

Let us give an explicit example on $\T^3$.
Consider the centered random walk on $\Z^2$ with distribution $\nu$
supported on $\el_1 = (2, 1)$, $\el_2 = (1, -2)$, $\el_3 =(-3, 1)$,
such that $\PP(X_0 = \el_j) = \frac13$, for $j= 1, 2, 3$. With the
notation of Subsection \ref{prelimRW}, here $D$ is the sublattice
generated by $\{(1, 3), (4, -3)\}$ and we have $\Cal D^\perp =
\{\0\}$ and $a_0 = 15$.

Let $A_1$ and $A_2$ be the commuting matrices computed in Subsection
\ref{torEx} (example 1.a). Let the matrices $B_j, j= 1, 2, 3,$ be
defined by $B_1 = A_1^{2} A_2 = \left(
\begin{matrix} 29 & 23 & -8 \cr -80 & -67 & 23 \cr 230 & 196 & -67 \cr
\end{matrix} \right),$ \\
$B_2 = A_1 A_2^{-2}= \left(
\begin{matrix} -13 & -11 & 4 \cr 40 & 35 & -11 \cr -110 & -92 & 35 \cr
\end{matrix} \right),$
$B_3 = A_1^{-3} A_2= \left(
\begin{matrix} 107 & 16 & -7 \cr -70 & 23 & 16 \cr 160 & 122 & 23 \cr
\end{matrix} \right)$.

The random walk on $\T^3$ such that with equal probability we move
from $x \in \T^3$ to $B_j x$, $j= 1, 2, 3$, gives rise to a random
process$(U_n(\omega) x)_{n \geq 1}$ on the torus such that for a
regular function $f$ the limiting distribution of $((n \Log
n)^{-\frac12} \sum_{k = 0}^{n-1} f(V_k(\omega) .))$ is a normal law,
with variance $\sigma_f(0) = c \, \varphi_f(0)$, where the constant
$c$ is given by the LLT.

\subsection{\bf Powers of barycenter operators}

\

We show now that the iterates of the barycenter operators satisfy
the condition of Theorem \ref{tclRegKer}. Let ($B_j, j \in J)$ be a set of
endomorphisms satisfying Assumption \ref{hypo1}, i.e.,
such that $B_j = A^{\el_j}$, $\el_j \in \Z^d$, where $A_1, ..., A_d$ are
$d$ algebraically independent commuting automorphisms of $\tilde G$.
We suppose that the corresponding  $\Z^d$-action $\A$ on
$(\tilde G, \tilde \mu)$ is totally ergodic.

Let $P$ be the barycenter operator $Pf(x) := \sum_{j \in J} \, p_j
\, f(B_j x)$. The associated random walk $(Z_n)$ on $\Z^d$ is defined by $\PP(X_0
= \ell_j) = p_j$, for $j \in J$. The lattice $L(\tilde W)$ generated
by the support of the distribution of $\tilde W$ coincides with
$D(W)$. The lattices $D(W)$ and $D(\tilde W)$ are the same. We use
the LLT (Theorem \ref{LLT2}) for $\tilde W$ (with the exponent
$d(\tilde W) = d_0(W) = d(W)$ or $d_0(W) -1$ and $\Lambda$ replaced
by $\Lambda_0$).

We suppose that $D(W)$ is not trivial, so that $W$ is not
deterministic and $d_0 = d_0(W) \geq 1$. The measure $d\gamma_1$ was
defined in Notation \ref{defdgamma} (see also
\ref{barycprocsubsec}).

\begin{thm} \label{tclBaryc0} Let $f$ be a function in $AC_0(G)$ with spectral
density $\varphi_f$. Then, for a constant $C$ depending on the
random walk, we have w.r.t. the Haar measure on $G$:
$$C \, n^{d_0\over 4} \, P^n f \overset{distr} {\underset{n \to \infty}
\longrightarrow} \Cal N(0, \sigma_P^2(f)), \ with \ \sigma_P^2(f) =
\int \, \varphi_{f} \, d\gamma_1.$$
In particular, if $B_j = A_j, j= 1, ..., d$, then
$$(4\pi)^{{d_0 \over 4}} \, (p_1 ... p_{d})^{\frac14} \, n^{d_0\over 4}
\, P^n f \overset{distr} {\underset{n \to \infty} \longrightarrow
} \Cal N(0, \sigma_P^2(f)), \ with \ \sigma_P^2(f) = \int_{\T^1} \,
\varphi_{f}(u, u, ..., u) \, du.$$
\end{thm}
\proof We apply Proposition \ref{regBaryc}. Here $R_n(\el) = \PP(Z_n
= \el)$, $\sum_{\el \in \Z^d} R_n(\el) = 1$ and by Theorem
\ref{LatLocThm}, $\sup_\el R_n(\el) = O(n^{-d_0/2})$. If
$\sigma_P^2(f) \not = 0$, the variance $\sigma_n^2(f)$ is
asymptotically like $\sum_{\el \in \Z^d} \, R_n(\el)^2  = \PP(\tilde
Z_n = \0) \sim C_0 n^{-d_0/2}$, for a constant $C_0 > 0$.

Here Condition (\ref{cumulCN3}) reads $\bigl(\sup_\el
R_n(\el)\bigr)^{r-1} = o(n^{- rd_0/4}), \, \forall r \geq 3$. It is
satisfied, since $n^{-(r-1) \, d_0/2} = o(n^{- rd_0/4}), \, \forall
r \geq 3$. We conclude by Theorem \ref{tclRegKer}. \eop

{\bf Examples and remarks.}

1) Let $A_1, A_2$ be two commuting matrices with coefficients in
$\cal{M}(\rho,\Z)$ generating a totally ergodic $\Z^2$-action on
$\T^\rho$, $\rho \geq 3$. Let $P$ be the barycenter operator $Pf(x)
:= p_1 f(A_1 x) + p_2 f(A_2 x)$, $p_1, p_2 > 0, p_1 + p_2 = 1$.

The symmetrized r.w. $(\tilde Z_n)$ is 1-dimensional and the support
of the distribution of $\tilde X_0$ is $\{(0,0), \, (1,-1), (-1,1)
\}$. If $\varphi_f$ is continuous, then $\lim_{n\to\infty} \sqrt{4
p_1 p_2 \pi n} \, \|P^n f\|_2^2 = \int_\T \, \varphi_{f}(u,u) \,
du$. If $f$ satisfies the regularity condition (\ref{5.11}) on
$\T^\rho$, we have, with $\sigma_P^2(f) = \int_\T \,
\varphi_{f}(u,u) \, du$, $(4 p_1 p_2 \pi)^\frac14 \, n^\frac14 \,
P^n f \overset{distr} {\underset{n \to \infty} \longrightarrow }
\Cal N(0,\sigma_P^2(f))$.

If $f$ is in $AC_0(G)$, then $\sigma_P(f) = 0$ if and only if
$\varphi_{f}(u,u) = 0$, for every $u \in \T^1$. In particular, if
$f$ is not a mixed coboundary (cf. Theorem \ref{coboundChar}), then
$\sigma_P(f) \not = 0$ and the rate of convergence of $\|P^nf\|_2$
to 0 is the polynomial rate given by Theorem \ref{tclBaryc0}.

2) Suppose now that $A_1, A_2$ are automorphisms. We can take, for
examples the matrices computed in Subsection \ref{torEx}. Let $P$ be
defined by $Pf(x) := p_1 f(A_1 x) + p_2 f(A_2 x) + p_3 f(A_1^{-1} x)
+ p_4 f(A_2^{-1} x)$, $p_j > 0$, $\sum_j p_j = 1$.

We have analogous result, excepted that here $L(\tilde W) = D(\tilde
W) = D(W) = \Z \, (1, 1) \oplus \Z \, (-1, 1)$. The measure
$d\gamma_1$ is the barycenter of $\delta_{(0, 0)}$ and
$\delta_{(\frac12, \frac12)}$.

3) Observe that the previous example for the quenched limit theorem
used for its computation the tables of units in a field number. As
remarked in \Sect \ref{endoGroup}, for barycenters, it is easier to
give examples which are endomorphisms. The maps on the torus are not
necessarily invertible, but the analysis of the process uses a
symmetrized r.w. on $\Z^d$. The difficulty is then to compute the
dimension of the associated r.w. unless it is given by primality
conditions of the determinants. Let us give a simple example.

Let $P$ on $L_0^2(\T^1)$ be defined by: $(Pf)(x) = \frac18 f(2x) +
\frac28 f(3x) + \frac18 f(5x) + \frac38 f(6x) + \frac18 f(15x)$. The
behavior of the iterates $P^n$ is given by the LLT applied to the
symmetrized r.w. $(\tilde Z_n)$ (strictly aperiodic in $\Z^3$) with
distribution supported on $\tilde \Sigma = \{(0, 0, 0), \pm (1, 0,
0), \pm (0, 1, 0), \pm (0, 0, 1), \pm (1, -1, 0), \pm (1, 0, -1),
\pm (0, 1, -1), \pm (1, 1, -1), \\ \pm (1, -1, -1)\}$. If $f$ is
regular on $\T^1$, the process $(n^{-3/4} P^nf(.))_{n \geq 1}$
converges in distribution to a normal law (non degenerate $f$ is not
a mixed coboundary).

4) For $\nu$ a discrete measure on the semigroup $\Cal T$ of
commuting endomorphisms of $G$, let us consider a barycenter of the
form $Pf(x) = \sum_{T \in \Cal T} \, \nu(T) f(Tx)$. When there is a
finite moment of order 2 and $d(W) < + \infty$, the decay of $P^n f$
is of order $n^{-{d(W) \over 4}}$ when $\varphi_f$ is continuous and
$\varphi_f(0) \not = 0$. A question is to estimate the decay when
$\nu$ has an infinite support and $d(W)$ is infinite and to study
the asymptotic distribution of the normalized iterates (if there is
a normalization).

 For example if $Pf(x) = \sum_{q \in \Cal P} \nu(q) f(qx)$, where
$\Cal P$ is the set of prime numbers and $(\nu(q), q \in \Cal P)$ a
probability vector with $\nu(q) > 0$ for every prime $q$, what is
the decay to 0 of $\|P^nf\|_2$, when $f$ is H\"olderian on the
circle?

When  $d_0(W)$ is infinite, a partial result is that the decay is
faster than $Cn^{-r}$, for every $r \geq 1$. This follows from the
following observation:

Let $P_1$ and $P_2$ be commuting contractions of $L^2(G)$ such that
$\|P_1^nf\|_2 \leq Mn^{-r}$. Let $\alpha \in ]0, 1], \beta = 1 -
\alpha$. Then we have: $\|(\alpha P_1 + \beta P_2)^n f) \|_2 \leq
\sum_{k= 0}^n {n\choose k} \alpha^k \beta^{n-k} \|P_1^k f\|_2.$
There is $c < 1$ such that $\sum_{k \leq {n \over 2\alpha}}
{n\choose k} \alpha^k \beta^{n-k} \leq c^n$; therefore: $\sum_{k=
0}^n {n\choose k} \alpha^k \beta^{n-k} \|P_1^k f\|_2 \leq M ({n
\over 2\alpha})^{-r} + c^n \leq M' n^{-r}$.

5) The case of commutative or amenable actions strongly differs from
the case of non amenable actions for which a ``spectral gap
property" is available in certain cases (\cite{FuSh99}), which
implies a quenched CLT theorem (cf. (\cite{CoLeb11}).

For actions by algebraic {\it non commuting} automorphisms $B_j,
j=1, ..., d$, on the torus, the existence of a spectral gap for $P$
of the form $Pf(x) := \sum_{j \in J} \, p_j \, f(B_j x)$ is related
to the fact that the generated group has no factor torus on which it
is virtually abelian (cf. \cite{BeGu11}). In general, a question is
to split the action into a (virtually) abelian action on the
$L^2$-space of a factor torus and a supplementary subspace with a
spectral gap, in order to get a full description of the quenched and
the barycenter processes.

\vskip 3mm {\bf Coboundary characterization}

For a $\Z^d$-action $\Cal S$ with Lebesgue spectrum by automorphisms
on $\T^\rho$ and $f$ on $\T^\rho$, let us give a characterization
for $\varphi_f(0) = 0$ in terms of coboundaries.

For the dual action of $\Z^d$ on $\Z^\rho$, we construct the section
$J_0$ introduced before Proposition \ref{condAbsConv} in the
following way. For a fixed $\el$, the set $\{A^\k \el, \k \in
\Z^d\}$ is discrete and $\lim_{\|\k\| \to \infty} \|A^\k \el\| =
+\infty$. Therefore we can choose an element $\j$ in each class
modulo the action of $\Cal S$ on $\Z^\rho$ which achieves the
minimum of the norm. By this choice, we have
\begin{eqnarray}
\|\j\| \leq \|A^\k \j\|, \, \forall \j \in J_0, \, \k \in \Z^d.
\label{majJK0}
\end{eqnarray}

\begin{thm} \label{coboundChar} If $|c_f(\k)| = O(\|k\|^{-\beta}),
\text{ with } \beta > \rho$, then $\varphi_f(0) = 0$ if and only if
$f$ satisfies the following mixed coboundary condition: there are
continuous functions $u_i$, $1 \leq i \leq d$ such that
\begin{eqnarray}
f= \sum_{i=1}^d (I- A_i) u_i. \label{cobound12}
\end{eqnarray}
\end{thm}
\proof Let $\varepsilon \in \, ]0, \beta - \rho[$. If $\delta =
(\beta - \rho - \varepsilon) / (\beta (1 + \rho))$, we have $\delta
\beta\rho -\beta(1-\delta) = - (\rho + \varepsilon)$. There is a
constant $C_1$ such that $\|\k\|^d \, e^{-\delta \beta \tau \|\k\|}
\leq C_1, \, \forall \k \in \Z^d$. According to (\ref{DaKaMin}), we
have $|c_f(A^\k \j)| \leq C \|A^\k \j\|^\beta \leq C e^{-\beta \tau
\|\k\|} \, \|\j\|^{\beta \rho}$; hence
\begin{eqnarray}
e^{\delta \beta \tau \|\k\|} \, |c_f(A^\k j)|^\delta \leq C
\|\j\|^{\delta \beta\rho}. \label{majkj1}
\end{eqnarray}
For every $\el \in \Z^\rho \stm0$, there is a unique $(k, j) \in
\Z^d \times J_0$ such that $A^\k \j = \el$; hence by Inequality
(\ref{majJK0}):
\begin{eqnarray*}
&&\sum_{\j \in J_0}\sum_{\k \in \Z^d} \|\k\|^d \, |c_f(A^\k \j)| =
\sum_{\j \in J_0} \sum_{\k \in \Z^d} \|\k\|^d \, |c_f(A^\k \j)|^\delta
|c_f(A^\k \j)|^{1-\delta} \\
&& \ \ \leq C_1 \sum_{\j \in J_0}\sum_{\k \in \Z^d} e^{\delta \beta
\tau \|\k\|} \, |c_f(A^\k \j)|^\delta \, |c_f(A^\k \j)|^{1-\delta} \
\leq C_2 \sum_{\j \in J_0}\sum_{\k \in \Z^d} \|\j\|^{\delta
\beta\rho} \, |c_f(A^\k \j)|^{1-\delta}.
\end{eqnarray*}
According to (\ref{majkj1}) and (\ref{majJK0}), the right hand side is less than
\begin{eqnarray*}
&&C_2 \sum_{\j \in J_0}\sum_{\k \in \Z^d} \|A^\k
\j\|^{\delta \beta\rho} \, |c_f(A^\k \j)|^{1-\delta}
 = C_2 \sum_{\el \in \Z^\rho \stm0} \|\el\|^{\delta \beta\rho}
\, |c_f(\el)|^{1-\delta} \\
&&\leq C_3 \sum_{\el \in \Z^\rho \stm0}
\|\el\|^{\delta \beta\rho -\beta(1-\delta)} \leq C_3 \sum_{\el \in
\Z^\rho \stm0} \|\el\|^{-(\rho + \varepsilon)}< +\infty.
\end{eqnarray*}
The sufficient condition for (\ref{cobound12}) given in
\cite{CohCo13} reads here: $\sum_{j \in J_0} \sum_{\k \in \Z^d} (1
+\|\k\|^d) \, |c_f(A^\k j)| < \infty$. This condition holds by the
previous inequality. Since here the functions involved in the proof
of the coboundary characterization are characters, hence continuous
and uniformly bounded, the functions $u_i$ in (\ref{cobound12}) are
continuous. \eop

\section{\bf Appendix I: self-intersections of a centered
r.w.} \label{selfInt}

In this appendix, we prove Theorem \ref{recd12}. {\it As already
noticed, we can assume aperiodicity in the proofs.}

\vskip 3mm \subsection{\bf d=1: a.s. convergence of $V_{n, p} / V_{n,
0}$}

\

We need the following lemmas.
\begin{lem} \label{min1d} If $W$ is a 1-dimensional r.w. with finite variance
and centered, then, for a.e. $\omega$, there is $C(\omega) > 0$ such
that
\begin{eqnarray} V_n(\omega) \geq C(\omega) \,
n^{\frac32} \, (\Log \Log \, n)^{-\frac12}. \label{minVn1}
\end{eqnarray}
\end{lem} \proof By the law of iterated logarithm, there
is a constant $c > 0$ such that, for a.e. $\omega$, the inequality
$|Z_n(\omega)| > c \, (n \ \Log \Log \, n)^\frac12$ is satisfied
only for finitely many values of $n$. This implies that, for a.e.
$\omega$, there is $N(\omega)$ such that $|Z_n(\omega)| \leq (c \, n
\ \Log \Log \, n)^\frac12$, for $ n \geq N(\omega)$. It follows, for
$N(\omega) \leq k < n$, $|Z_k(\omega) \, | \leq (c \, k \ \Log \Log
k)^\frac12 \leq (c \, n \ \Log \Log \, n)^\frac12$.

Therefore, if $\Cal R_n(\omega)$ is the set of points visited by the
random walk up to time $n$, we have $\Card (\Cal R_n(\omega)) \leq 2
(c_1(\omega) \, n \ \Log \Log \, n)^\frac12$, with an a.e. finite
constant $c_1(\omega)$.

We have $n = \sum_{\el \in \Z} \sum_{k=0}^{n-1} 1_{Z_k(\omega) \,  =
\el} \leq (\sum_{\el \in \Cal R_n(\omega)} (\sum_{k=0}^{n-1}
1_{Z_k(\omega) \,  = \el})^2)^\frac12 \, \Card(\Cal
R_n(\omega))^\frac12$ by Cauchy-Schwarz inequality; hence:
$V_n(\omega) \geq {n^2 / \Card(\Cal R_n(\omega))}$, which implies
(\ref{minVn1}). \eop

\vskip 3mm \begin{lem} For an aperiodic r.w. in dimension 1 with
finite variance and centered, we have:
\begin{eqnarray}
\sup_n \, |\sum_{j=1}^n [2\,\PP(Z_j = 0) - \PP(Z_j = \p) - \PP(Z_j =
-\p)]| < + \infty, \ \forall \p \in L(W). \label{boundPot}
\end{eqnarray}
\end{lem} \proof
Since $1 - \Psi(t)$ vanishes on the torus only at $t = \0$ with an
order 2, we have:
\begin{eqnarray*}
&&|\sum_{j=1}^{N-1} [2\PP(Z_j = 0) - \PP(Z_j = \p) - \PP(Z_j =
-\p)]| = 4 |\int_{\T} \sin^2 \pi \langle\p, t\rangle \, \Re e ({1
- \Psi^N(t) \over 1 - \Psi(t)}) \, dt|\\
&&\leq 4 \int_{\T} \sin^2 \pi \langle\p, t\rangle \, |{1 - \Psi^N(t)
\over 1 - \Psi(t)}| \, dt \leq 8 \int_{\T} { \sin^2 \pi \langle\p,
t\rangle \, \over |1 - \Psi(t)|} \, dt < + \infty.
\end{eqnarray*} \eop

{\bf Proof of Theorem \ref{recd12}} (for $d=1$: $\displaystyle
\lim_n {V_{n, p}(\omega) \over V_{n}(\omega)} = 1$, a.e. for every
$p \in L(W)$.)

Recall that $V_{n}(\omega) - V_{n, p}(\omega) \geq 0$ (cf. Ex.
\ref{exples2}). By (\ref{VnpErgSum}) we can bound $n^{-1} \, \E
[V_{n} - V_{n, p}]$ by
\begin{eqnarray*}
1 + n^{-1} \sum_{k = 1}^{n-1} \, |\sum_{j = 0}^{n-k-1} (\PP(Z_j =
\p) - \PP(Z_j = 0)) \, + \, \sum_{j = 0}^{n-k-1} (\PP(Z_j = -\p) -
\PP(Z_j = 0))|
\end{eqnarray*}
which is bounded according to  (\ref{boundPot}).

Let $\delta > 0$. The bound $\displaystyle \E \, ({V_{n}(\omega) -
V_{n, p} \over n^{\frac32}} ) = O(n^{-\frac12})$ implies by the
Borel-Cantelli lemma that $\displaystyle
({V_{n}(\omega) - V_{n, p} \over n^{\frac32}} )$ tends to 0 a.e.
along the sequence $k_n = [n^{2+\delta}], \, n \geq 1$.

Since $V_{n} \geq c(\omega) \, n^\frac32 \, (\Log \Log \,
n)^{-\frac12}$ by (\ref{minVn1}), we obtain
$$0 \leq  1 - {V_{n, \p}(\omega) \over V_{n}(\omega)} \leq {V_{n}
(\omega) - V_{n, p}(\omega) \over V_{n}(\omega)} < {V_{n}(\omega) -
V_{n, p}(\omega) \over c(\omega) \, n^\frac32 \, (\Log \Log \,
n)^{-\frac12}}.$$

Therefore $({V_{n, \p}(\omega) \over V_{n}(\omega)})_{n \geq 1}$
converges to 1 a.e. along the sequence $(k_n)$.

To complete the proof, it suffices to prove that a.s. $\lim_n
\max_{k_n \le j<k_{n+1}} |V_{j,p} / V_j - V_{k_n ,p} / V_{k_n}| =
0$. By monotonicity of $V_{j,p}$ and $V_j$ with respect to $j$, we
have
$$\frac{ V_{k_{n},p }}{V_{k_{n}}} \frac{V_{k_{n}}}{ V_{k_{n+1}}}
\le\frac{V_{j,p} }{ V_{j}}\le \frac{ V_{k_{n+1},p }}{V_{k_{n+1}}}
\frac{V_{k_{n+1}}}{ V_{k_n}}, \ k_n \le j <k_{n+1}.$$ Therefore,
since the first factors in the left and right terms of the
inequality tend to 1, it is enough to prove that
$\frac{V_{k_{n+1}}-V_{k_n}} {V_{k_n}}\to 0$ a.s.

By (\ref{IneqepsilRn2}), for $d= 1$, for each $\varepsilon
> 0$, there is an a.e. finite constant  $c_\varepsilon(\omega)$ such
that $\sup_{\el \in \Z} R_n(\omega,\el) = \sup_{\el \in \Z}
\sum_{k=0}^{n-1} 1_{Z_k(\omega) \, = \el} = c_\varepsilon(\omega) \,
n^{\frac12 + \varepsilon}$. This implies
\begin{eqnarray*}
&&V_{k_{n+1}}-V_{k_n}=\sum_{\ell, j=0}^{k_{n+1}} 1_{\{Z_\ell=Z_j\}}-
\sum_{\ell, j=0}^{k_{n}} 1_{\{Z_\ell=Z_j\}}\le
\sum_{\ell=k_n+1}^{k_{n+1}}\sum_{j=1}^{k_{n+1}} 1_{\{Z_\ell=Z_j\}}+
\sum_{j=k_n+1}^{k_{n+1}} \sum_{\ell=1}^{k_{n+1}} 1_{\{Z_\ell=Z_j\}}
\\ &&\le 2 (k_{n+1}- k_n) \, \sup_{p \in \Z} \sum_{j=1}^{k_{n+1}} 1_{\{Z_j =
p\}} \leq 2 c_\varepsilon(\omega) \, (k_{n+1}- k_n) \,
k_{n+1}^{\frac12 + \varepsilon} \leq K(\omega) \,  n^{1+\delta + (2
+ \delta)(\frac12 + \varepsilon)}.
\end{eqnarray*}

Therefore, $V_{k_{n+1}}-V_{k_n} \leq K(\omega) \,  n^{2+
\frac32\delta + 2\varepsilon + \varepsilon \delta}$ and
$\frac{V_{k_{n+1}}-V_{k_n}} {V_{k_n}}$ is a.s. bounded by
\begin{eqnarray*}
&&{V_{k_{n+1}}-V_{k_n} \over k_n^{3/2} \, (\Log\Log \,
k_n)^{-\frac12}} \leq 2 K(\omega) \,  {n^{2+ \frac32\delta +
2\varepsilon + \varepsilon \delta} \over (n^{2+\delta})^{3/2} \,
(\Log\Log \, (n^{2+\delta}))^{-\frac12}} \\
&& \leq 2 K(\omega) \, {n^{2+ \frac32\delta + 2\varepsilon +
\varepsilon \delta} \, n^{-(3+ \frac32 \delta)} \, ((\Log\Log \,
(n^{2+\delta}))^{\frac12}} = 2 K(\omega) \, {n^{-1 + 2\varepsilon +
\varepsilon \delta} \, (\Log\Log \, (n^{2+\delta}))^{\frac12}},
\end{eqnarray*}
which tends to 0, if $2\varepsilon + \varepsilon \delta < 1$. \eop

\vskip 3mm \subsection{\bf d=2: variance and SLLN for $V_{n,p}$}

\

For $d = 2$, in the centered case with finite variance, the a.s.
convergence $V_{n, \p} / V_{n, \0} \to 1$ for $p \in L(W)$ follows
from the strong law of large numbers (SLLN): $\lim_n V_{n, \p}/ \E
V_{n, \p} = 1$, a.s. We adapt the method of \cite{Lew93} to the case
$\p \not = \0$ in the estimation of $\Var(V_{n, \p})$. See also
\cite{DeUt11} for the computation of the variance of $V_{n, \0}$. We
need two auxiliary results.

\begin{lem} \label{extLewis}
There is $C$ such that, if $\p, \q$ are in $L$ and $n, k$ are such
that $n \el_1 = \p \mod D$ and $(n+k) \el_1 = \q \mod D$, then
\begin{eqnarray}
|\PP(Z_{n+k} = \q) - \PP(Z_n = \p)| \leq C ({1 \over(n+k)^{\frac32}}
+ { k\over n (n+k)}), \, \forall n, k \geq 1. \label{majDiffPPnk}
\end{eqnarray}
\end{lem}
\proof We have $\PP(Z_{n+r} = \q) - \PP(Z_n = \p) = \int_{\T^2}
G_{n,r}(\t) \, d\t$, with
\begin{eqnarray*}
G_{n,r}(\t) := \Re e \, [e^{-2\pi i \langle \q, \t \rangle} \,
\Psi(\t)^{n+r} -e^{-2\pi i \langle \p, \t \rangle} \, \Psi(\t)^{n}].
\end{eqnarray*}

The functions $e^{- 2\pi i \langle \q, \t \rangle} \,
\Psi(\t)^{n+r}$ and $e^{- 2\pi i \langle \p, \t \rangle} \,
\Psi(\t)^{n}$ are invariant by translation by the elements of
$\Gamma_1$ and have a modulus $< 1$, except for $\t \in \Gamma_1$
(cf. Lemma \ref{modulus1}). To bound the integral of $G_{n,r}$, it
suffices to bound its integral $I_n^0 := \int_{U_0} G_{n,r}(\t) _,
d\t$ restricted to a fundamental domain $U_0$ of $\Gamma_1$ acting
on $\T^2$.

Denote by $B(\eta)$ the ball with a small radius $\eta$ and center
$\0$ in $\T^2$. If $\eta >0$ is small enough, on $U_0 \setminus
B(\eta)$, we have $|\Psi(\t)| \leq \lambda(\eta)$ with
$\lambda(\eta) < 1$, which implies:
\begin{eqnarray*}
|I_n^0| \leq 2 C \lambda^n(\eta) + \int_{B(\eta)} |G_{n,r}(\t)| \,
d\t.
\end{eqnarray*}
and we have
\begin{eqnarray*}
|G_{n,r}(\t)| := |\Re e \, [e^{-2\pi i \langle \q, \t \rangle} \,
\Psi(\t)^{n+r} -e^{-2\pi i \langle \p, \t \rangle} \, \Psi(\t)^{n}]|
\leq |\Psi(\t)|^{n} \, |\Psi(\t)^{r} - e^{2\pi i \langle \q - \p, \t
\rangle} |.
\end{eqnarray*}
 Since the distribution $\nu$ of the centered r.w. $W$
is assumed to have a moment of order 2, by Lemma \ref{quadrF}, for
$\eta$ sufficiently small, there are constants $a, b>0$ such that
$$|\Psi(t)|  < 1-a \|t\|^2, \ |1-\Psi(t)|  <b \|t\|^2, \forall t \in B(\eta).$$

We distinguish two cases: \hfill \break if $\p = \q$, then
$|G_{n,r}(\t)| \leq  C(r) (1-a \|t\|^2)^n \|\t\|^2$; \hfill \break
if $\p \not = \q$, $|G_{n,r}(\t)| \leq C(r) |\Psi(\t)|^{n} \|\t\|
\leq C(r) (1-a \|t\|^2)^n \, \|\t\|$.

Now we bound the integral $\int_{B(\eta)}\ (1 - a {\|\t\|^2})^{n} \,
\|\t\|^2\, {dt_1 \, dt_2}$. By the change of variable $(t_1, t_2)
\to ({s_1, \over \sqrt n}, {s_2, \over \sqrt n}) $, then $(s_1, s_2)
\to (\rho \cos \theta, \rho \sin \theta)$, it becomes successively:
$${1 \over n^2} \, \int_{B(\eta)}\ (1 - a{\|\s\|^2 \over n})^n  \, \|\s\|^2 \, ds_1 \, ds_2
 \leq {C \over n^2} \, \int_\R e^{- a \rho^2} \, \rho^2 \, d\rho = O({1 \over n^2}).$$
So for $\p = \q$, we get the bound $O(n^2)$. Likewise,  if $\p \not
= \q$, we get the bound $O(n^{3/2})$.

Recall the $L / D$ is a cyclic group (Lemma \ref{DLatt}). If $n, k$
satisfy the condition of the lemma, we write: $k= r+ uv$, where $r$
and $v$ ($= |L/D|$) are the smallest positive integers such that
respectively: $r \el_1 = \q - \p \mod D$, $v \el_1 = \0 \mod D$.
Then, using the previous bounds and writing the difference
$\PP(Z_{n+k} = \q) - \PP(Z_n = \p)$ as
\begin{eqnarray*}
&&\PP(Z_{n+k} = \q) - \PP(Z_{n+r+(u-1)v} = \q)
+ \PP(Z_{n+r+(u-1)v} = \q) - \PP(Z_{n+r+(u-2)v} = \q) \\
&&+ ... + \PP(Z_{n+r+v} = \q) - \PP(Z_{n+r} = \q) + \PP(Z_{n+r} =
\q) - \PP(Z_n = \p),
\end{eqnarray*}
the telescoping argument gives (\ref{majDiffPPnk}). \eop

\begin{lem} \label{sumLambda}
For $d \geq 1$, let $M_n^{(d)} := \{(m_1, ..., m_d) \in \N^d: \,
\sum_i m_i = n\}$. Let $\tilde M_n^{(5)} := \{(m_1, ..., m_5) \, \in
M_n^{(5)}, m_3, \, m_4 > 0 \}$. If $|\lambda|, |\alpha|, |\beta|,
|\gamma| < 1$, then
\begin{eqnarray}
&&\sum_{n= 0}^\infty \lambda^{n} \sum_{(m_1, ..., m_5) \, \in \tilde
M_n^{(5)}} \alpha^{m_2} \, \beta^{m_3} \, \gamma^{m_4} = {1 \over (1
- \lambda)^2} {\lambda^2 \beta \gamma \over (1 - \lambda \alpha) (1
- \lambda \beta) (1 - \lambda \gamma)}. \label{Mn5}
\end{eqnarray}
\end{lem}
\Proof We have for $\lambda, \alpha_1, ..., \alpha_d$ such that
$|\lambda\alpha_1|, ..., |\lambda \alpha_d| \, \in [0, 1[$:
\begin{eqnarray}
\sum_{n= 0}^\infty \lambda^{n} \sum_{(m_1, ..., m_d) \, \in
M_n^{(d)}} \alpha_1^{m_1} \, \alpha_2^{m_2} \, ... \, \alpha_d^{m_d}
= \prod_{i=1}^d {1 \over (1 - \lambda \alpha_i)}.
\label{sumLambdaForm}
\end{eqnarray}

Indeed, the left hand of (\ref{sumLambdaForm}) is the sum over $\Z$
of the discrete convolution product of the functions $G_i$ defined
on $\Z$ by $G_i(k) = 1_{[0, \infty[}(k) (\lambda \alpha_i)^k$, hence
is equal to
\begin{eqnarray*}
&& \sum_{k \, \in \Z} \, (G_1 * ... * G_d)(k) = \prod_{i=1}^d \,
\bigr(\sum_{k \in \Z} \, G_i(k)\bigr) = \prod_{i=1}^d {1 \over (1 -
\lambda \alpha_i)}.
\end{eqnarray*}
For $d=5$ we find
\begin{eqnarray*}
\sum_{n= 0}^\infty \lambda^{n} \sum_{(m_1, ..., m_5) \in M_n^{(5)}}
\alpha^{m_2} \, \beta^{m_3} \gamma^{m_4} = {1 \over (1 - \lambda)^2}
{1 \over (1 - \lambda \alpha) (1 - \lambda \beta) (1 - \lambda
\gamma)}.
\end{eqnarray*}
The left hand of (\ref{Mn5}) reads
\begin{eqnarray*}
&&\sum_{n= 0}^\infty \lambda^{n} \sum_{(m_1, ..., m_5) \, \in
M_n^{(5)}} \alpha^{m_2} \, \beta^{m_3} \, \gamma^{m_4} - \sum_{n=
0}^\infty \lambda^{n} \sum_{(m_1, m_2, m_4, m_5) \, \in M_n^{(4)}}
\alpha^{m_2} \, \gamma^{m_4}\\
&& - \sum_{n= 0}^\infty \lambda^{n} \sum_{(m_1, m_2, m_3, m_5) \,
\in M_n^{(4)}} \alpha^{m_2} \, \beta^{m_3} + \sum_{n= 0}^\infty
\lambda^{n} \sum_{(m_1, m_2, m_5) \, \in M_n^{(3)}} \alpha^{m_2}\\
&& = {1 \over (1 - \lambda)^2} [{1 \over (1 - \lambda \alpha) (1 -
\lambda \beta) (1 - \lambda \gamma)} - {1 \over (1 - \lambda \alpha)
(1 - \lambda \gamma)}\\ &&-{1 \over (1 - \lambda \alpha) (1 -
\lambda \beta)} + {1 \over (1 - \lambda \alpha)}] = {1 \over (1 -
\lambda)^2} {\lambda^2 \beta \gamma \over (1 - \lambda \alpha) (1 -
\lambda \beta) (1 - \lambda \gamma)}. \eop
\end{eqnarray*}

\vskip 3mm {\bf Proof of Theorem \ref{recd12}} (Case $d=2$)

\vskip 3mm A) Below we consider:
\begin{eqnarray}
\sum_{0 \leq i_1 < j_1 < n; \ 0 \leq i_2 < j_2 < n} \, \PP(Z_{j_1} -
Z_{i_1} = \r, \, Z_{j_2} - Z_{i_2} = \s). \label{sumij}
\end{eqnarray}

The (disjoint) sets of possible configurations for $(i_1, j_1, i_2,
j_2)$ are the following: \hfill \break (1a) \ $i_1 \leq i_2 < j_1 <
j_2$, (1b) \ $i_2 \leq i_1 < j_2 < j_1$, \hfill \break (2a) \ $i_1 <
i_2 < j_2 \leq j_1$, (2b) \ $i_2 < i_1 < j_1 \leq j_2$, \hfill
\break (3a) \ $i_1 < j_1 \leq i_2 < j_2$, (3b) \ $i_2 < j_2 \leq i_1
< j_1$, \hfill \break (4) \ $i_2 = i_1 < j_1 = j_2$.

For the case 4), by the LLT the sum (\ref{sumij}) restricted to the
configurations (4) is less than $\sum_{0 \leq i_1 < j_1 < n} 1_{\r =
\s} \, \PP(Z_{j_1} - Z_{i_1} = \r) \leq \sum {C \over j_1 - i_1}
\leq C n \, \Log \, n$.

(1a) and (1b), (2a) and (2b), (3a) and (3b) are respectively the
same up to the exchange of indices $1$ and $2$. To bound
(\ref{sumij}), it suffices to consider the subsums corresponding to
1a), 2a), 3a).

\vskip 3mm {\bf 1a)} ($i_1 \leq i_2 < j_1 < j_2$)  Setting $m_1 =
i_1, \ m_2 = i_2 - i_1, \ m_3 = j_1 - i_2, \ m_4 = j_2 - j_1, \ m_5
= n - j_2$, we have: $Z_{j_1} - Z_{i_1} = \tau^{m_1} Z_{m_2 + m_3}$,
$Z_{j_2} - Z_{i_2} = \tau^{m_1+m_2} Z_{m_3 + m_4}$ with $m_1, m_2
\geq 0, m_3, m_4 > 0$; hence:
\begin{eqnarray}
&&\PP(Z_{j_1} - Z_{i_1} = \r, \, Z_{j_2} - Z_{i_2} = \s) =
\PP(Z_{m_2+ m_3} = \r, \, \tau^{m_2} Z_{m_3 + m_4} = \s) \nonumber\\
&&= \sum_\el [\PP(Z_{m_2} = \el) \, \PP(Z_{m_3} = \r - \el)
\, \PP(Z_{m_4} = \s - \r + \el)] \nonumber \\
&&= \int e^{-2\pi i(\langle \r, \u \rangle + \langle \s - \r, \v
\rangle)} \sum_\el e^{-2\pi i \langle \el, \, \t - \u + \v \rangle}
\ \Psi(\t)^{m_2} \, \Psi(\u)^{m_3} \, \Psi(\v)^{m_4} \, d\t \, d\u
\, d\v. \label{case1a}
\end{eqnarray}
The last equation can be shown by approximating the probability
vector $(p_j)_{j \in J}$ by probability vectors with finite support
and using the continuity of $\Psi$.

\vskip 3mm {\bf 2a)} ($i_1 < i_2 < j_2 \leq j_1$) Setting $m_1 =
i_1, \ m_2 = i_2 - i_1, \ m_3 = j_2 - i_2, \ m_4 = j_1 - j_2, \ m_5
= n - j_1$, we have: $Z_{j_1} - Z_{i_1} = \tau^{m_1} Z_{m_2 + m_3 +
m_4}$, $Z_{j_2} - Z_{i_2} = \tau^{m_1+m_2} Z_{m_3}$, with $m_1, m_4
\geq 0, m_2, m_3 > 0$.

Since $\PP(Z_{m_2} + \tau^{m_2+m_3} Z_{m_4} = \q) = \sum_\el
\PP(Z_{m_2} = \el, \tau^{m_2+m_3} Z_{m_4} = \q - \el) = \sum_\el
\PP(Z_{m_2} = \el) \, \PP(\tau^{m_2+m_3} Z_{m_4} = \q - \el) =
\sum_\el \PP(Z_{m_2} = \el) \, \PP(\tau^{m_2} Z_{m_4} = \q - \el) =
\PP(Z_{m_2+m_4} = \q)$, we have:
\begin{eqnarray*}
&&\PP(Z_{j_1} - Z_{i_1} = \r, \, Z_{j_2} - Z_{i_2} = \s)=
\PP(Z_{m_2+ m_3+m_4} = \r, \, \tau^{m_2}Z_{m_3} = \s) \nonumber \\
&&= \PP(Z_{m_2} + \tau^{m_2+m_3} Z_{m_4} = \r - \s, \,  \tau^{m_2}
Z_{m_3} = \s) \nonumber \\
&&= \PP(Z_{m_2} + \tau^{m_2+m_3} Z_{m_4} = \r - \s) \, \PP(Z_{m_3} =
\s) = \PP(Z_{m_2 + m_4} = \r - \s) \, \PP(Z_{m_3} = \s).
\end{eqnarray*}
Hence, in case 2a), we get:
\begin{eqnarray}
\PP(Z_{j_1} - Z_{i_1} = \r, \, Z_{j_2} - Z_{i_2} = \s) &=&
\PP(Z_{m_2 + m_4} = \r - \s) \, \PP(Z_{m_3} = \s). \label{1case2a}
\end{eqnarray}

{\bf 3a)} ($i_1 < j_1 \leq i_2 < j_2$) The events $Z_{j_1} - Z_{i_1}
= \r$ and $Z_{j_2} - Z_{i_2} = \s$ are independent.

\vskip 3mm Following the method of \cite{Lew93}, now we estimate
$\Var(V_{n, \p}) = \E(V_{n, \p}^2) - (\E V_{n, \p})^2$. Recall that
$V_{n, \p} = \sum_{0 \leq i, j < n} \, 1_{Z_{j} - Z_{i} = \p} =
\sum_{0 \leq i < j < n} \, (1_{Z_{j} - Z_{i} = \p} + 1_{Z_{j} -
Z_{i} = -\p}) + n \, 1_{\p = \0}$. We have
\begin{eqnarray*}
V_{n, \p}^2 &&= \sum_{0 \leq i_1 < j_1 < n; \ 0 \leq i_2 < j_2 < n}
\, (1_{Z_{j_1} - Z_{i_1} = \p} + 1_{Z_{j_1} - Z_{i_1} = -
\p}) \ (1_{Z_{j_2} - Z_{i_2} = \p} + 1_{Z_{j_2} - Z_{i_2} = - \p})\\
&&\ \ + 2n  \, 1_{p=\0} \sum_{0 \leq i < j < n} \, (1_{Z_{j} - Z_{i}
= \p} + 1_{Z_{j} - Z_{i} = - \p}) + n^2 \, 1_{p=\0}.
\end{eqnarray*}
The last term gives 0 in the computation of the variance, so that it
suffices to bound
$$\sum_{0 \leq i_1 < j_1 < n; \ 0 \leq
i_2 < j_2 < n} \, \E [ (1_{Z_{j_1} - Z_{i_1} = \p} + 1_{Z_{j_1} -
Z_{i_1} = - \p}) \ (1_{Z_{j_2} - Z_{i_2} = \p} + 1_{Z_{j_2} -
Z_{i_2} = - \p})]$$
$$-\sum_{0 \leq i_1 < j_1 < n; \ 0 \leq i_2 <
j_2 < n} \, [\PP(Z_{j_1} - Z_{i_1} = \p) + \PP(Z_{j_1} - Z_{i_1} = -
\p)] \ [\PP(Z_{j_2} - Z_{i_2} = \p) + \PP(Z_{j_2} - Z_{i_2} = -
\p)],$$ i.e., the sum over the finite set $(\r,\s) \in \{\pm \p\}$
of the sums
\begin{eqnarray*}
\sum_{0 \leq i_1 < j_1 < n; \ 0 \leq i_2 < j_2 < n} \, [\PP(Z_{j_1}
- Z_{i_1} = \r, \ Z_{j_2} - Z_{i_2} = \s) - \PP(Z_{j_1} - Z_{i_1} =
\r) \ \PP(Z_{j_2} - Z_{i_2} = \s)].
\end{eqnarray*}

B) By the previous analysis, we are reduced to cases (1a), (2a),
(3a). For (3a), by independence, we get 0. As $\E(V_{n, \p}^2) - (\E
V_{n, \p})^2 \geq 0$, for the bound of (1a), we can neglect the
corresponding centering terms, since they are subtracted and non
negative.

{\bf (1a)} For $\r, \s \in \Z^2$, let $a_{\r, \s}(n):= \sum_{0 \leq
i_1 \leq i_2 < j_1 < j_2 < n} \, \PP(Z_{j_1} - Z_{i_1} = \r, \,
Z_{j_2} - Z_{i_2} = \s),$ $a(\p, n):= \sum_{\r, \s \in \{\pm \p\}}
a_{\r, \s}(n)$.

Setting $m_1 = i_1, m_2= i_2 - i_1, m_3 = j_1 - i_2, m_4 = j_2 -
j_1$, $m_5 = n- j_2$ (so that $(m_1, ... , m_5)$ runs in the set
$\tilde M_n^{(5)}$ of Lemma \ref{sumLambda}), by (\ref{case1a}) we
have for the generating function $A_{\r, \s}(\lambda) := \sum_{n
\geq 0} \lambda^n a_{\r, \s}(n)$, $0 \leq \lambda < 1$:
\begin{eqnarray*}
A_{\r, \s}(\lambda) &=& \int_{\T^2}\int_{\T^2} \, e^{-2\pi i
(\langle\r, \t \rangle + \langle \s, \u \rangle)} \, \sum_n
\lambda^n \sum_{(m_1, ... , m_5) \, \in \, \tilde M_n^{(5)}} \,
\Psi(\t)^{m_2} \, \Psi(\t + \u)^{m_3} \, \Psi(\u)^{m_4} \, d\t \,
d\u.
\end{eqnarray*}
Using  $\sum_{\r, \s \in \{\pm \p\}} e^{-2\pi i (\langle\r, \t
\rangle + \langle \s, \u \rangle)} = \cos 2\pi\langle\p, \t \rangle
\, \cos 2 \pi \langle \p, \u \rangle$, for the generating function
$A_\p(\lambda) := \sum_{n \geq 0} \lambda^n a(\p, n)$, we have by
(\ref{sumLambdaForm}) with $\alpha = \Psi(\t), \beta = \Psi(\u),
\gamma = \Psi(\t+\u)$:
$$A_\p(\lambda) = 4 {\lambda^2 \over (1 - \lambda)^2}
\int \int {\cos 2\pi\langle\p, \t \rangle \, \cos 2 \pi \langle \p,
\u \rangle \, \Psi(\u) \Psi(\t+\u) \over (1 - \lambda \Psi(\t)) \,
(1- \lambda \Psi(\u)) \, (1 - \lambda \Psi(\t + \u))} d\t d\u.$$

For an aperiodic r.w., the bound obtained for $A_\p(\lambda)$ in
\cite{Lew93} for $\p = \0$ is valid. Indeed the bounds for the
domains $T_{i, \delta}$ and $E_\delta$ hold: for $T_{i, \delta}$
this is clear; for $E_\delta$ this is because $(1 - \Psi(\t)) \, (1-
\Psi(\u)) \, (1 - \Psi(\t + \u))$ does not vanish on $E_\delta$.
(See \cite{Lew93} p. 227 for the notation for the domains.) The main
contribution comes from a small neighborhood $U_\delta$ of diameter
$\delta > 0$ of $\0$ and the factor $\cos 2\pi\langle\p, \t \rangle
\, \cos 2 \pi \langle \p, \u \rangle$ plays no role.

It follows that $\displaystyle A_\p(\lambda) \leq {C \over (1 -
\lambda)^3}$. Since $a(\p, n)$ increases with $n$, we obtain $a(\p,
n) = O(n^2)$ by using the following elementary ``Tauberian"
argument:

Let $(u_k)$ be a non decreasing sequence of non negative numbers. If
there is $C$ such that $\displaystyle \sum_{k=0}^{+\infty} \lambda^k
u_k \leq {C \over (1 - \lambda)^3}, \ \forall \lambda \in [0, 1[$,
then $u_n \leq C' n^2, \ \forall n \geq 1$. Indeed we have:
$$\lambda^n u_n = (1 - \lambda)(\sum_{k=n}^{+\infty} \lambda^k) \, u_n \leq
(1 - \lambda) \sum_{k=n}^{+\infty} \lambda^k u_k \leq (1 - \lambda)
\sum_{k=0}^{+\infty} \lambda^k u_k  \leq {C \over (1 -
\lambda)^2}.$$ For $\lambda = 1 - {1\over n}$ in the previous
inequality, we get: $u_n \leq C n^2 \, (1 - {1\over n})^{-n} \leq C'
n^2$.

\vskip 3mm  {\bf (2a)} Let $b_{\r, \s}(n)$ be the sum
\begin{eqnarray*}
\sum_{0 \leq i_1 < i_2 < j_2 \leq j_1 < n} [\PP(Z_{j_1} - Z_{i_1} =
\r, \, Z_{j_2} - Z_{i_2} = \s) - \PP(Z_{j_1} - Z_{i_1} = \r) \, \PP(
Z_{j_2} - Z_{i_2} = \s)]
\end{eqnarray*}
and $b(\p, n) := \sum_{\r, \s \in \{\pm \p\}} b_{\r, \s}(n)$.

Putting $m_1 = i_1, m_2 = i_2 - i_1, m_3 = j_2 - i_2, m_4 = j_1 -
j_2$, we have by (\ref{1case2a}):
$$b(\p, n) = \sum_{\r, \s \in \{\pm \p\}} \sum_{{\underset{m_1, m_4 \geq 0, \, m_2, m_3 \geq 1} {m_1
+ m_2 + m_3 + m_4 < n}}} \PP(Z_{m_3} = \s) [\PP(Z_{m_2 + m_4} = \r -
\s) - \PP(Z_{m_2 + m_3+ m_4} = \r)].$$ We will use Lemma
\ref{extLewis} to bound $b(\p, n)$.

If $m_3 \el_1 = \s \mod D$, then either $(m_2 + m_4) \el_1 = \r - \s
\mod D$ and $(m_2+m_4 + m_3) \el_1 = \r \mod D$, or $(m_2 + m_4)
\el_1 \not = \r - \s \mod D$ and $(m_2+m_4 + m_3) \el_1 \not = \r
\mod D$. In the latter case, the corresponding probabilities are 0.
Moreover, $\PP(Z_{m_3} = \s) =0$ if $m_3 \el_1 \not = \s \mod D$.
Therefore the sum reduces to those indices such that $m_3 \el_1 = \s
\mod D$, $(m_2 + m_4) \, \el_1 = \r - \s \mod D$ and $(m_2+m_4 +
m_3) \, \el_1 = \r \mod D$ and the bound (\ref{majDiffPPnk}) applies
for the non zero terms. We obtain:
$$b(\p, n) \leq 4 \sum_{{\underset{m_1, m_4 \geq 0, \, m_2, m_3 \geq 1} {m_1
+ m_2 + m_3 + m_4 < n}}} [{1 \over m_3} {m_3 \over (m_2 + m_4) (m_2
+ m_3 + m_4)} + {1 \over m_3 \, (m_2 + m_4 + m_3)^\frac32}].$$ The
sum of the first term restricted to $\{m_i\ge 1, m_1+m_2+m_3+m_4\le
n\}$ reads:
\begin{eqnarray*}
&&\sum_{k=1}^n\sum_{m_1+m_2+m_3+m_4=k}\frac{1}{(m_2+m_4)(m_2+m_4+m_3)}
= \sum_{k=1}^n\sum_{\ell=1}^{k}\sum_{m_2+m_4=\ell,
m_1+m_3=k-\ell}\frac{1}{\ell(\ell+m_3)} \\ &&\le
\sum_{k=1}^n\sum_{\ell=1}^{k}\sum_{m_1+m_3=k-\ell}\ell
\cdot\frac{1}{\ell(\ell+m_3)} \le
\sum_{k=1}^n\sum_{\ell=1}^{k}(\frac{1}{\ell}+\cdots+\frac{1}{k}) =
\frac12 n(n+1).
\end{eqnarray*}

Likewise, we have:
\begin{eqnarray*}
&&\sum_{k=1}^n\sum_{m_1+m_2+m_3+m_4=k}\frac{1}{m_3(m_2+m_4+m_3)^\frac32}
\leq n \sum_{k=1}^n\sum_{m_2+m_3+m_4=k}\frac{1}{m_3 \, k^\frac32} \\
&&\le n\sum_{k=1}^n\frac{1}{k^{3/2}}\sum_{m_3=1}^k\frac{k}{m_3}
\le n\sum_{k=1}^n\frac{\Log k}{\sqrt{k}}\le C\, n^{3/2}\Log n.
\end{eqnarray*}
Therefore, we obtain $b(n, \p) = O(n^2)$.

The previous estimations imply $\Var(V_{n, \p}) = O(n^2)$. By
(\ref{EVnp}), for $\p \in L(W)$, $\lim_n V_{n, \p} / \E V_{n, \p} =
1$ a.e. follows as in \cite{Lew93} and therefore $\lim_n V_{n, \p} /
V_{n, \0} = 1$ a.e. \eop

\section{\bf Appendix II: mixing, moments and cumulants}
\label{cumulAppend}

 For the sake of completeness, we recall in this appendix some general results
on mixing of all orders, moments and cumulants (see \cite{Leo60c}
and the references given therein). Implicitly we assume existence of
moments of all orders when they are used.

For a real random variable $Y$ (or for a probability distribution on
$\R$), the cumulants (or semi-invariants) can be formally defined as
the coefficients $c^{(r)}(Y)$ of the cumulant generating function $t
\to \ln \E(e^{tY}) = \sum_{r=0}^\infty \, c^{(r)}(Y) \,{t^r \over
r!}$, i.e., $c^{(r)}(Y) = {\partial^r \over \partial^r t} \ln
\E(e^{tY})|_{t = 0}.$ Similarly the joint cumulant of a random
vector $(X_1, ... , X_r)$ is defined by
\begin{eqnarray}
c(X_1, ... , X_r) = {\partial^r \over \partial t_1 ... \partial t_r}
\ln \E(e^{\sum_{j=1}^r t_jX_j})|_{t_1 = ... = \, t_r = 0}.
\label{defCumul}
\end{eqnarray}
This definition can be given as well for a finite measure $\nu$ on
$\R^r$. Its cumulant is noted $c_\nu(x_1, ..., x_r)$.
The joint cumulant of $(Y, ...,Y)$ ($r$ copies of $Y$) is $c^{(r)}(Y)$.

 For any subset $I = \{i_1, ..., i_p\} \subset J_r:= \{1, ..., r\}$,
we put
$$m(I) = m(i_1, ..., i_p):= \E(X_{i_1} \cdots X_{i_p}),
\ s(I)= s(i_1, ..., i_p):= c(X_{i_1}, ... , X_{i_p}).$$ The
cumulants of a process $(X_j)_{j \in \Cal J}$, where $\Cal J$ is a
set of indexes, is the family
$$\{c(X_{i_1}, ..., X_{i_r}), ({i_1}, ..., {i_r}) \in \Cal J^r, r \geq 1\}.$$
The following formulas link moments and cumulants and vice-versa:
\begin{eqnarray}
c(X_1, ... , X_r) &=& s(J_r) = \sum_{\Cal P} (-1)^{p-1} (p - 1)!
\ m(I_1) \cdots m(I_p), \label{cumFormu1}\\
\E(X_{1} \cdots X_{r}) &=& m(J_r) = \sum_{\Cal P} s(I_1) \cdots s(I_p),
\label{cumFormu2}
\end{eqnarray}
where in both formulas, $\Cal P = \{I_1, I_2, ..., I_p\}$ runs
through the set of partitions of $J_r = \{1, ..., r\}$ into $p \le
r$ nonempty intervals, with $p$ varying from 1 to $r$.

Now let be given a random field of real random variables $(X_\k)_{\k
\in \Z^d}$ and a summable weight $R$ from $\Z^d$ to $\R^+$. For $Y
:= \sum_{\el \in \Z^d} \, R(\el) \, X_\el$ we obtain from the
definition (\ref{defCumul}):
\begin{eqnarray}
c^{(r)}(Y) = c(Y, ..., Y) = \sum_{(\el_1, ..., \el_r) \, \in
(\Z^d)^r} \, c(X_{\el_1}, ... , X_{\el_r}) \, R(\el_1) \cdots
R(\el_r). \label{cumLin0}
\end{eqnarray}

{\bf Limiting distribution and cumulants}

For our purpose, we state in terms of cumulants a particular case of
a theorem of M. Fr\'echet and J. Shohat, generalizing classical
results of A. Markov. Using the formulas linking moments and
cumulants, a special case of their ``generalized statement of the
second limit-theorem" can be expressed as follows:

\begin{thm} \cite{FreSho31} \label{FreSho} Let $(Z^n, n \geq 1)$ be a sequence
of centered real r.v. such that
\begin{eqnarray}
\lim_n c^{(2)}(Z^n) = \sigma^2, \ \lim_n c^{(r)}(Z^n) = 0, \forall r
\geq 3, \label{limCum}
\end{eqnarray}
then $(Z^n)$ tends in distribution to $\Cal N(0, \sigma^2)$. (If
$\sigma = 0$, then the limit is $\delta_0$).
\end{thm}

\begin{thm} \label{Leonv} (cf. Theorem 7 in \cite{Leo60b})
Let $(X_\k)_{\k\in \Z^d}$ be a random process and $(R_n)_{n \geq 1}$
a summation sequence on $\Z^d$. Let $(Y^n)$ be defined by $Y_n =
\sum_\el \, R_n(\el) \, X_\el, n \geq 1$. If
\begin{eqnarray}
\sum_{(\el_1, ..., \el_r) \, \in (\Z^d)^r} \, c(X_{\el_1}, ... ,
X_{\el_r}) \, R_n(\el_1) ... R_n(\el_r) = o(\|Y^n\|_2^r), \forall
r \geq 3, \label{smallCumul}
\end{eqnarray}
then ${Y^n \over \|Y^n\|_2}$ tends in distribution to $\Cal N(0, 1)$
when $n$ tends to $\infty$.
\end{thm}
\proof \ Let $\beta_n := \|Y^n\|_2 = \|\sum_\el \, R_n(\el) \,
X_\el\|_2$ and $Z^n = \beta_n^{-1} Y^n$.

Using (\ref{cumLin0}), we have $c^{(r)}(Z^n) = \beta_n^{-r}
\sum_{(\el_1, ..., \el_r) \, \in (\Z^d)^r} \, c(X_{\el_1}, ... ,
X_{\el_r}) \, R(\el_1) ... R(\el_r)$. The theorem follows then from
the assumption (\ref{smallCumul}) by Theorem \ref{FreSho} applied to
$(Z^n)$. \eop
\begin{defn} A measure preserving $\N^d$ (or $\Z^d$)-action $T: {\n} \to T^\n$
on a probability space $(X, \Cal A, \mu)$ is $r$-mixing, $r > 1$, if
for all sets $B_1, ..., B_r \in \Cal A$
$$\lim_{\min_{{1 \leq \ell < \ell' \leq r}} \|\n_\ell - \n_{\ell'}\| \to
\infty} \mu(\bigcap_{\ell = 1}^r T^{-\n_\ell} \, B_\ell)=\prod_{\ell
= 1}^r \mu(B_\ell).$$
\end{defn}
\begin{nota} {\rm For $f$ in the space $L_0^\infty(X)$ of measurable essentially
bounded functions on $(X, \mu)$ with $\int f \, d\mu = 0$, we apply the
definition of moments and cumulants to $(T^{\n_1}f,..., T^{\n_r}f)$
and put}
\begin{eqnarray}
m_f(\n_1,..., \n_r) = \int_X \, T^{\n_1} f \cdots T^{\n_r} f \, d\mu, \
\ s_f(\n_1,..., \n_r) := c(T^{\n_1} f,..., T^{\n_r} f).
\label{notasf}
\end{eqnarray}
\end{nota}
To use the property of mixing of all orders, we need the following lemma.
\begin{lem} \label{subSeqLem} For every sequence $(\n_1^k,..., \n_r^k)$ in
$(\Z^d)^r$, there are a subsequence with possibly a permutation of
indices (still written $(\n_1^k,..., \n_r^k)$), an integer
$\kappa(r) \in [1, r]$, a subdivision $1 = r_{1} < r_{2} < ... <
r_{{\kappa(r)-1}} < r_{{\kappa(r)}} \leq r$ of $\{1, ..., r\}$ and a
constant integral vector $\a_j$ such that
\begin{eqnarray}
&&\lim_k \min_{1 \leq s \not = s' \leq \kappa(r)}\|\n_{r_{s}}^k -
\n_{r_{s'}}^k\| = \infty, \label{minCond}\\
&&\n_j^k = \n_{r_{s}}^k + \a_j, \for r_{s} < j < r_{{s+1}}, \ s =
1, ..., \kappa(r)-1, \text{ and for } r_{\kappa(r)} < j \leq r.
\label{constant1}
\end{eqnarray}
If $(\n_1^k,..., \n_r^k)$ satisfies $\lim_k \max_{i \not = j}
\|\n_i^k - \n_j^k\| = \infty$, then $\kappa(r) > 1$.
\end{lem}
\proof \ Remark that if $\sup_k \max_{i \not = j} \|\n_i^k -
\n_j^k\| < \infty$, then $\kappa(r) = 1$ so that (\ref{minCond}) is
void and (\ref{constant1}) is void for the indexes such that
$r_{{s+1}} = r_{s} +1$. The proof of the lemma is by induction. The
result is clear for $r= 2$. Suppose that the subsequence for the
sequence of $(r-1)$-tuples $(\n_1^k,..., \n_{r-1}^k)$ is built.

Let $1 \leq r_{1} < r_{2} < ... < r_{{\kappa(r-1)}} \leq r-1$ be the
corresponding subdivision of $\{1, ..., r -1\}$, as stated above for
the sequence $(\n_1^k,..., \n_{r-1}^k)$. If
$(\n_1^k,..., \n_{r-1}^k)$ satisfies $\lim_k \max_{1 < i < j \leq
r-1} \|\n_i^k - \n_j^k\| = \infty$, then $\kappa(r-1) > 1$ by
construction in the induction process.

Now consider $(\n_1^k,..., \n_r^k)$. If $\lim_k \|\n_r^k - \n_i^k\|
= +\infty$, for $i=1, ..., r-1$, then we take $1 \leq r_{1} < r_{2}
< ... < r_{{\kappa(r-1)}} < r_{{\kappa(r)}} = r$ as new subdivision
of $\{1, ..., r\}$. If $\liminf_k \|\n_r^k - \n_{i_s}^k\| < +\infty$
for some $s \leq \kappa(r-1)$, then along a new subsequence (still
denoted $\n_r^k$) we have $\n_r^k = \n_{i_s}^k + \a_r$, where $\a_r$
is a constant integral vector. After changing the labels, we insert
$n_r$ in the subdivision of $\{1, ..., r-1\}$ and obtain the new
subdivision of $\{1, ..., r\}$.

For the last condition on $\kappa$, suppose that $\lim_k \max_{1 < i
< j \leq r} \|\n_i^k - \n_j^k\| = \infty$. Then, if $\liminf_k
\max_{1 < i < j \leq r-1} \|\n_i^k - \n_j^k\| < +\infty$,
necessarily, $\kappa(r) > 1$. If, on the contrary, the sequence
$(\n_1^k,..., \n_{r-1}^k)$ satisfies $\lim_k \max_{1 < i < j \leq
r-1} \|\n_i^k - \n_j^k\| = \infty$, then $\kappa(r-1) > 1$, so that
$\kappa(r) \geq \kappa(r-1) > 1$. \eop
\begin{lem} \label{skTozeroLem} If a $\Z^d$-dynamical system on
$(X, \Cal A, \mu)$ is mixing of order $r \geq 2$, then, for any $f
\in L_0^\infty(X)$, $\underset{\max_{i \not = j} \|\n_i - \n_j\| \to
\infty} \lim s_f(\n_1,..., \n_r) = 0$.
\end{lem}
{\proof} The notation $s_f$ was introduced in
(\ref{notasf}). Suppose that the above convergence does not hold.
Then there is $\varepsilon > 0$ and a sequence of $r$-tuples
$(\n_1^k =\0,..., \n_r^k)$ such that $|s_f(\n_1^k,..., \n_r^k)| \geq
\varepsilon$ and $\max_{i \not = j} \|\n_i^k - \n_j^k\| \to \infty$
(we use stationarity).

By taking a subsequence (but keeping the same notation), we can
assume that, for two fixed indexes $i, j$, $\lim_k \|\n_i^k -
\n_j^k\| = \infty$. By Lemma \ref{subSeqLem}, there is a subdivision
$1 = r_{1} < r_{2} < ... < r_{{\kappa(r)-1}} < r_{{\kappa(r)}} \leq
r$ and constant integer vectors $\a_j$ such that
\begin{eqnarray}
&&\lim_k \min_{1 \leq s \not = s' \leq \kappa(r)}\|\n_{r_{s}}^k -
\n_{r_{s'}}^k\| = \infty, \label{minCond2}\\
&&\n_j^k = \n_{r_{s}}^k + \a_j, \for r_{s} < j < r_{{s+1}}, \ s =
1, ..., \kappa(r)-1, \text{ and for } r_{\kappa(r)} < j \leq r.
\end{eqnarray}
Let $d\mu_k(x_1, ..., x_r)$ denote the probability measure on $\R^r$
defined by the distribution of the random vector $(T^{\n_1^k} f(.),
..., T^{\n_r^k}f(.))$. We can extract a converging subsequence from
the sequence $(\mu_k)$, as well as for the moments of order $\leq
r$. Let $\nu(x_1, ..., x_r)$ (resp. $\nu(x_{i_1}, ..., x_{i_p})$) be
the limit of $\mu_k(x_1, ..., x_r)$ (resp. of its marginal measures
$\mu_k(x_{i_1}, ..., x_{i_p})$ for $\{i_1, ..., i_p\} \subset \{1,
..., r\}$). Let $\varphi_i, i=1, ..., r$, be continuous functions
with compact support on $\R$. Mixing of order $r$ and condition
(\ref{minCond2}) imply
\begin{eqnarray*}
&&\nu(\varphi_1 \otimes \varphi_2 \otimes ... \otimes \varphi_r)
=\lim_k \int_{\R^d} \varphi_1 \otimes \varphi_2 \otimes ... \otimes
\varphi_r \, d\mu_k = \lim_k \int \prod_{i=1}^r
\varphi_i(f(T^{\n_i^k}x))\ d\mu(x) \\
&&=\lim_k \int \ \bigl[\prod_{s=1}^{\kappa(r)-1} \ \prod_{r_{s} \leq
j < r_{{s+1}}} \varphi_j(f(T^{\n_s^k + \a_j} x))\bigr] \,
\prod_{\kappa(r) \leq j \leq r} \varphi_j(f(T^{\n_{\kappa(r)}^k +
\a_j} x)) \ d\mu(x)
\end{eqnarray*}
\begin{eqnarray*}
= \bigl[\prod_{s=1}^{\kappa(r)-1} \, \bigl(\int \prod_{r_{s} \leq
j < r_{{s+1}}} \varphi_j(f(T^{\a_j}x)) \ d\mu(x) \bigr)\bigr] \,
\bigl[\int \prod_{\kappa(r) \leq j \leq r} \varphi_j(f(T^{\a_j} x))
\ d\mu(x)\bigr].
\end{eqnarray*}
Therefore $\nu$ is the product of marginal measures corresponding to
disjoint subsets: there are $I_1 = \{i_1, ..., i_p\}, I_2 = \{i_1',
..., i_{p'}'\}$ two nonempty subsets of $J_r = \{1, ..., r\}$, such
that $(I_1, I_2)$ is a partition of $J_r$ and $d\nu(x_1, ...,x_r) =
d\nu(x_{i_1}, ..., x_{i_p}) \times d\nu(x_{i_1'}, ...,
x_{i_{p'}'})$.

With $\Phi(t_1, ..., t_r) = \ln \int e^{\sum t_j x_j} \ d\nu(x_1,
...,x_r)$ and the analogous formulas for $\nu(x_{i_1}, ...,
x_{i_p})$ and $\nu(x_{i_1'}, ..., x_{i_p'})$, we get $\Phi(t_1, ...,
t_r) = \Phi(t_{i_1}, ..., t_{i_p}) + \Phi(t_{i_1'}, ..., t_{i_p'})$;
hence $c_\nu(x_1, ..., x_r) = {\partial^r \over \partial t_1 ...
\partial t_r}\Phi(t_1, ..., t_r)|_{t_1 = ... = \, t_r = 0} = 0$,
which contradicts $\liminf_k |s_f(\n_1^k,..., \n_r^k)| > 0$. \eop

\vskip 3mm {\bf Acknowledgements} This research was carried out
during visits of the first author to the IRMAR at the University of Rennes 1 and
of the second author to the Center for Advanced Studies in
Mathematics at Ben Gurion University.
The first author was partially supported by the ISF grant 1/12.
The authors are grateful to their hosts for their support.
They thank Y. Guivarc'h, S. Le Borgne and M. Lin for helpful discussions
and B. Weiss for the reference \cite{Lev13}.

\end{document}